\documentclass{amsbook}
\usepackage[centertags]{amsmath}
\usepackage{amsfonts}
\usepackage{amssymb}
\usepackage{amsthm}
\usepackage[ansinew]{inputenc}
\usepackage[cmtip,all]{xy}
\usepackage{graphicx}
\newtheorem{thm}[equation]{Theorem}
\newtheorem{cor}[equation]{Corollary}
\newtheorem{lem}[equation]{Lemma}
\newtheorem{prop}[equation]{Proposition}
\theoremstyle{definition}
\newtheorem{defn}[equation]{Definition}
\theoremstyle{remark}
\newtheorem{rem}[equation]{Remark}

\numberwithin{section}{chapter}
\numberwithin{equation}{section}

\newcommand{\set}[1]{\left\{#1\right\}}
\newcommand{\Real}{\mathbb R}

\newcommand{\To}{\longrightarrow}

\newcommand{\A}{\mathcal{A}}

\newcommand{\mo}{\mathrm{\mathbf{mod}}}
\newcommand{\chain}{\mathrm{\mathbf{chain}}}

\newcommand{\Ab}{\mathrm{\mathbf{Ab}}}
\newcommand{\ab}{\mathrm{\mathbf{ab}}}
\newcommand{\gr}{\mathrm{\mathbf{gr}}}
\newcommand{\nil}{\mathrm{\mathbf{nil}}}
\newcommand{\M}{\mathrm{\mathbf{M}}}

\renewcommand{\S}{\Sigma}

\def\fy{\varphi}
\def\ker{\operatorname{Ker}}
\def\op{{op}}
\def\ni{{nil}}
\def\abb{{ab}}
\def\im{\operatorname{Im}}
\def\coker{\operatorname{Coker}}
\def\card{\operatorname{card}}
\def\hom{\operatorname{Hom}}
\def\ext{\operatorname{Ext}}

\def\tor{\operatorname{Tor}}
\def\endo{\operatorname{End}}

\def\colim{\mathop{\operatorname{colim}}}
\newcommand{\grupo}[1]{\langle #1\rangle}
\def\F{\mathfrak{F}} 
\def\Gam{\Gamma} 

\def\Z{\mathbb{Z}}

\def\N{\mathbb{N}}

\def\st{\stackrel} 
\def\ul{\underline} 
\def\ol{\overline} 
\newcommand{\um}[1]{\mathcal{#1}} 
\newcommand{\uf}[1]{\mathbb{#1}} 
\newcommand{\mat}[1]{\mathsf{#1}}

\newcommand{\cc}{\mathcal{C}} 
\newcommand{\hh}{\mathcal{H}} 
\def\h{\mathcal{H}} 
\newdir{ >}{{}*!/-10pt/\dir{>}} 

\def\r{\rightarrow} 
\newcommand{\C}[1]{\mathbf{#1}} 
\def\ul{\underline}

\def\NN{\N_0}

\begin{document}
\frontmatter
\title{On the proper homotopy type of locally compact $A^2_n$-polyhedra}%
\author{Fernando Muro}%
\address{Max-Planck-Institut f\"ur Mathematik, Vivatsgasse 7, 53111 Bonn, Germany}
\email{muro@mpim-bonn.mpg.de}
\thanks{The author was partially supported
by the project MTM2004-01865 and the MEC postdoctoral fellowship
EX2004-0616.}



\begin{abstract}
In this paper we address the classification problem for locally compact $(n-1)$-connected $CW$-complexes with dimension $\leq n+2$ up to proper homotopy type. We obtain
complete classification theorems in terms of purely algebraic data in those cases where the representation type of the involved algebra is finite. For this we define
new quadratic functors in controlled algebra and new homotopy and cohomology invariants in proper homotopy theory.
\end{abstract}

\maketitle

\setcounter{page}{4} \tableofcontents

\chapter*{Introduction}

The program for the homotopy classification of polyhedra in terms of algebraic invariants was initiated by J. H. C. Whitehead in the late 1940s. He completely solved in a series of papers the homotopy classification problem for $(n-1)$-connected compact
$CW$-complexes with dimension $\leq n+2$, the so called
$A^2_n$-polyhedra, for $n\geq 2$, see \cite{htskp}, \cite{sc4p} and
\cite{aces}. Later progress was made by H.-J. Baues and his collaborators, see the survey article \cite{htsur} on the history and the development of Whitehead's program.

The proper homotopy type of a non-compact space is known to be a
stronger topological invariant than the ordinary homotopy type.
Proper maps, unlike ordinary continuous maps, keep track of the
infinity behavior which arises naturally in non-compact spaces.
Moreover, proper maps are still suitable to do homotopy theory with
them. There are two main axiomatic approachs to proper homotopy theory in the literature. The Edwards-Hastings embedding of the proper category into the Quillen model category of pro-spaces (\cite{edhas}) and the Baues-Quintero approach in terms of $I$-categories (\cite{iht}) which is the one we follow in this paper since it is intrinsic.

So far Whitehead's classification program is at a very preliminary stage in the proper homotopy setting. The only class of non-compact spaces for which a nice algebraic model is available is given by the $n$-dimensional spherical objects under a locally compact tree $T$. Such a spherical object is obtained by pasting $n$-spheres on $T$ in a locally finite way. As a consequence of results in \cite{iht} and \cite{projdim} the proper homotopy category of these objects is equivalent to the category of finitely generated projective modules over the ring $\Z(\F(T))$ for any $n\geq 2$, see Sections \ref{ch2} and \ref{ephai} below. This ring, up to isomorphism, only depends on the space of Freudenthal ends $\F(T)$ of the tree $T$. For $T$ a one-ended tree the ring $\Z(\F(T))$ is the well-known ring of row-column-finite matrices, or locally finite matrices, over the integers. This ring plays also a role in algebraic $K$-theory, it is the cone of $\Z$ in the sense of \cite{deloopw}, see also \cite{kloday}. In general $\Z(\F(T))$ is the endomorphism ring of a certain object in a category of free controlled modules, compare \cite{canc}.

Spherical objects, despite being very simple, play a crucial role since they are the building blocks of any locally compact finite-dimensional $CW$-complex $X$, therefore the algebraic models for spherical objects allow the construction of a cellular homology
$$\hh_*X$$
with values in $\Z(\F(T))$-modules, see \cite{iht}. This is the homology of the proper cellular chain complex $\cc_*X$ which is a bounded complex of finitely generated projective $\Z(\F(T))$-modules. The cohomology of this complex with coefficients in a $\Z(\F(T))$-module $\um{M}$ is the cohomology of $X$
$$H^*(X,\um{M}).$$
These are just abelian groups since $\Z(\F(T))$ is non-commutative.

In this paper we consider locally compact $A^2_n$-polyhedra in
proper homotopy theory. These spaces also arose early in the
development of proper homotopy theory for independent reasons. The
discovery of proper homology theories detecting proper homotopy
equivalences led to the question of the existence and uniqueness of
Moore spaces in proper homotopy, i. e. $1$-connected locally compact
$CW$-complexes with non-vanishing homology in a single dimension. It
was noticed in \cite{tsukuba95} that the proper homology of locally compact
$CW$-complexes
can have projective dimension $\geq 2$, and actually no more than $2$ by the main result in \cite{projdim}. As a consequence proper Moore spaces exist and they are
locally compact $A^2_n$-polyhedra, but they are not unique, i. e.
the proper homotopy type is not determined by the homology. The
first explicit example of two different proper Moore spaces with the same homology was obtained in \cite{tsukuba95}. At this
point it was reasonable to wonder how many different proper Moore
spaces are there with the same homology in the same dimension, or
more generally, about the proper homotopy classification of proper
Moore spaces in terms of purely algebraic data. When we addressed
this problem we noticed that there was not substantial difference
between considering only proper Moore spaces or more generally all
locally compact $A^2_n$-polyhedra. 

For the homotopy classification of compact $A^2_n$-polyhedra 
Whitehead used homology and the Pontrjagin-Steenrod invariant in cohomology. He also showed that one can also use as an alternative tool his ``certain'' exact sequence for the Hurewicz homomorphism, which can be deduced from the Pontrjagin-Steenrod invariant. A generalization of this sequence is also
available in proper homotopy theory. However it is not useful since
two proper Moore spaces with the same homology in the same degree have
always the same exact sequence but they need not have the same
proper homotopy type. We then turned to consider the existence of
proper Pontrjagin-Steenrod invariants. When one tries to construct
these invariants in proper cohomology by imitating the classical
methods one finds again an obstacle coming from the projective
dimension, which is not $\leq 1$ like for abelian groups.
Nevertheless the possibility of the existence of adequate proper
Pontrjaging-Steenrod invariants can not be discarded just because of
this. Actually we have been successful in answering the
question on the existence of proper Pontrjagin-Steenrod invariants
by placing the problem in the realm of cohomology of categories in
the sense of Hochschild-Mitchell \cite{rso} and Baues-Wirsching
\cite{csc}. For instance, the obstruction to the existence of a proper Pontrjagin-Steenrod invariant suitable to classify proper homotopy types of $A^2_n$-polyhedra, $n\geq 3$, is a class
\begin{equation*}\tag{a}
\set{\theta}\in H^1(\C{chain}_n(\Z(\F(T)))/\!\simeq\;,H^{n+2}(-,H_n\otimes\Z/2))
\end{equation*}
in the first cohomology of the homotopy category $\C{chain}_n(\Z(\F(T)))/\!\simeq$ of bounded chain complexes of finitely generated projective $\Z(\F(T))$-modules concentrated on dimensions $\geq n$, see Section \ref{epsi} below. The use of these techniques in obstruction theory and
homotopy classification problems goes back to \cite{ah}. Cohomology
of categories is strongly related to representation theory,
therefore one can not expect to be able to carry out all necessary
computations in the presence of a wild representation type. The
representation theory needed to deal with proper
Pontrjaging-Steenrod invariants was already developed in \cite{rca}.
Here we give a complete answer to the possible existence of proper
Pontrjaging-Steenrod invariants in the finite representation type
cases.

Chapter \ref{chI} is a basic introduction to the proper homotopy
category and controlled algebra following \cite{iht}. We include two
new results on controlled algebra, Propositions \ref{extigrande} and
\ref{esexacto}, which are proved in the two appendices of this
paper. We also recall from \cite{iht} the
definition of the proper homotopy and homology modules and proper
cohomology groups, as well as their basic properties. In Chapter \ref{qf} we define new quadratic
functors in controlled algebra. Among these new functors we have generalizations of the exterior square, the reduced tensor square, and Whitehead's quadratic functor 
$$\wedge^2_T,\hat{\otimes}^2_T,\Gamma_T\colon\mo(\Z(\F(T)))\To\mo(\Z(\F(T))).$$
These quadratic functors are applied in several
computations of proper homotopy invariants, which are
extensively used in Chapter \ref{chVI} for the definition of new
James-Hopf and cup-product invariants in proper homotopy theory, such as for example the reduced cup-product invariant of an $(n-1)$-connected locally compact finite-dimensional $CW$-complex $X$, $n\geq 3$,
$$\hat{\cup}_X\in H^{n+2}(X,\hat{\otimes}^2_T\hh_nX).$$ 
The cup-product cohomology invariants are well-behaved with respect to
obstruction theory. The necessary obstruction theory is recalled in
the last section of Chapter \ref{qf} following the approach given by the ``tower of
categories'' introduced in \cite{ah} and generalized in \cite{cfhh}.
Among the different cup-product invariants defined in Chapter
\ref{chVI} we have the so called chain cup-product. This is not the
cup-product invariant of a space but of a bounded chain complex $\cc_*$ of finitely generated projective $\Z(\F(T))$-modules concentrated in dimensions $\geq n$
$$\bar{\cup}_{\cc_*}\in H^{n+2}(\cc_*,\wedge^2_TH_n\cc_*),$$
as for example the proper
cellular chain complex $\cc_*X$ of a locally compact $A^2_n$-polyhedron $X$. The chain
cup-product $\cup_{\cc_*}$ of a chain complex $\cc_*$ concentrated in dimensions $2$, $3$ and $4$ is the obstruction to the existence of an $A^2_2$-polyhedron admitting a co-H-multiplication with proper cellular chain complex $\cc_*X=\cc_*$. Such an obstruction does not arise
in ordinary homotopy theory since ordinary Moore spaces are much
better behaved. The collection of all chain cup-products assembles to a $0$-dimensional class in cohomology of categories 
\begin{equation*}\tag{b}
\bar{\cup}\in H^0(\C{chain}_n(\Z(\F(T)))/\!\simeq\;,H^{n+2}(-,\wedge^2_TH_n)).
\end{equation*}
The computation of this cohomology class
is a major step towards the solution of the proper
Prontrjagin-Steenrod invariant problem. This computation is
addressed in Chapter \ref{chVII}, where we use techniques from
homotopical algebra developed in \cite{cosqc}, as well as the
localization theorem in cohomology of categories obtained in
\cite{fcc}. We also use in this paper the representation theory of $\Z(\F(T))\otimes\uf{F}_2$ studied in \cite{rca} and the homological computations in Chapter \ref{chIX}. As a consequence of the computations in Chapter \ref{chVII} we obtain in this paper
the first examples of proper Moore spaces in degree $2$ which are
not co-H-spaces, in contrast with the ordinary case, see Corollary \ref{kio} and Remark \ref{kio2} below.

In Chapter \ref{chVIII} we finally attack the problem of the
existence of proper Pontrjagin-Steenrod invariants. We show here that the cohomology class $\set{\theta}$ in (a) is the universal obstruction to the existence of proper Pontrjagin-Steenrod invariants, while the chain cup-product cohomology class (b) is the universal obstruction for the existence of proper Pontrjagin-Steenrod invariants which are in addition compatible with the cup-product invariants of $CW$-complexes in Chapter \ref{chVI}.
This proves, together with the computations carried out in Chapter \ref{chVII}, that for spaces with $3$ or more ends it is impossible to
obtain proper Pontrjagin-Steenrod invariants which are well related
to the cup-product invariants defined in Chapter \ref{chVI}. This
is a surprising result which emphasizes the contrast between
ordinary homotopy theory and proper homotopy theory. On the other hand we also
show in Chapter \ref{chVIII} that for spaces with no more than $3$ ends there exist
Pontrjagin-Steenrod invariants which are suitable to classify proper
homotopy types of locally compact $A^2_n$-polyhedra in the stable
range $n\geq 3$. For this we use the computations with finitely presented $\Z(\F(T))\otimes\uf{F}_2$-modules in Chapter \ref{chIX} in order to relate the cohomology classes (b) and (a) by a Bockstein exact sequence. In this way we deduce from the computation of (b) in Chapter \ref{chVII} that (a) vanishes whenever the tree $T$ has $3$ or less ends. This leads to one of the main results of this paper which asserts that for a locally compact $A^2_n$-polyhedron $X$ with less than $4$ ends and $n\geq 3$ the pair
$$(\cc_*X,\wp_n(X))$$
given by the proper cellular chain complex $\cc_*X$ and the proper Steenrod invariant
$$\wp_n(X)\in H^{n+2}(X,\hh_nX\otimes\Z/2)$$
constructed in Chapter \ref{chVIII} is an algebraic model for the proper homotopy type of $X$. See Theorem \ref{clasi1} in Section \ref{kia} for a precise statement of this result.
The implications of these results for proper Moore
spaces are considered in a separate section within Chapter \ref{chVIII}.

The unstable case $n=2$ and the case of spaces
with more than $4$ ends are out of any scope because of impediments
coming from representation theory, in the sense explained above. In
the remaining case, i. e. $4$-ended locally compact
$A^2_n$-polyhedra for $n\geq 3$, it should still be possible to
answer the question on proper Pontrjaging-Steenrod invariants and
maybe to extend the classification theorem proved here. Nevertheless
the representation theory in this case is connected to the
representation theory of the $4$-subspace quiver, see \cite{rca}.
This quiver has tame representation type, as it was shown by Nazarova in
\cite{nazarova}, but the classification of indecomposable
$4$-subspaces over the field with $2$ elements is quite intricate.
Therefore the algebraic computations carried out in the already
lengthy Chapter \ref{chIX} would have to be considerably expanded
to include this special case. For this reason we decided not to deal
with this case in this paper.

\mainmatter



\chapter{Proper homotopy theory and controlled algebra}\label{chI}

In this chapter we summarize the background on proper homotopy theory and controlled algebra which is necessary to understand the rest of this paper. We follow the axiomatic approach to proper homotopy theory in \cite{iht}. In Section \ref{ch2} we include two new results in controlled algebra, Propositions \ref{extigrande} and \ref{esexacto}, which are crucial in this paper in order to extend some results for a tree with a certain number of Freudenthal ends to a tree with more ends. We also establish connections with the results in \cite{projdim} on the projective dimension in controlled algebra.

\section{Basic aspects of proper homotopy theory}\label{icat}


A continuous map $f\colon X\r Y$ is \emph{proper} if it is closed
and $f^{-1}(y)$ is compact for all $y\in Y$. An alternative
definition says that $f$ is proper provided $f^{-1}(K)$ is compact
for any compact subspace $K\subset Y$. Both definitions are
equivalent on locally compact Hausdorff spaces, in particular over
locally compact $CW$-complexes, and this is the class of spaces on
which we will concentrate soon. All maps in this paper are supposed
to be proper, unless we explicitly state the contrary.

The ordinary \emph{cylinder} functor $X\mapsto IX=[0,1]\times X$
restricts to the category $\C{Topp}$ of spaces and proper maps.
Moreover, the structural maps of a cylinder
$$\xymatrix{X\ar@<.5ex>[r]^{i_0}\ar@<-.5ex>[r]_{i_1}&IX\ar[r]^p&X},$$
$i_k(x)=(k,x)$, $p(t,x)=x$, are proper, so one can define a
homotopy natural equivalence relation in this category.

In fact this cylinder gives rise to an $I$-category structure on
$\C{Topp}$, see \cite{iht} I.3.6. \emph{Cofibrations}
$X\rightarrowtail Y$ in this category are maps with the proper
homotopy extension property, in particular they are injective.
Examples of cofibrations are inclusions of subcomplexes into
finite-dimensional locally compact $CW$-complexes.

The axioms of an $I$-category, and the weaker ones of a cofibration
category, were introduced in \cite{ah}. There Baues develops the
homotopy theory of these kind of categories, with emphasis on
obstruction theory. He shows that avoiding the use of fibrations one
can still generalize most of ordinary homotopy theory. These axioms
are specially adequate for proper homotopy theory since $\C{Topp}$
lacks of enough inverse limits, even the product of two spaces does
not exist unless one of them is compact, so it does not fit into
Quillen's model category axioms.

A \emph{space $X$ under $A$} is a map $i_X\colon A\r X$ in
$\C{Topp}$. It is said to be \emph{cofibrant} if this map
$i_X\colon A\rightarrowtail X$ is a cofibration. 
The category $\C{Topp}^A_c$ of cofibrant spaces and maps under $A$
is also an $I$-category, see \cite{iht} I.3.12. A map under $A$ is a
cofibration if it is a cofibration in $\C{Topp}$. The cylinder
functor is the \emph{relative cylinder} $I_A$.
We shall write $$[X,Y]^A$$ for the set of homotopy classes of maps
$X\rightarrow Y$ in $\C{Topp}^A_c$.

A \emph{pair of spaces $(X,Y)$ under $A$} is a map $Y\rightarrow X$
in $\C{Topp}^A_c$. It is \emph{cofibrant} if $Y\rightarrowtail X$ is
a cofibration. The
category of cofibrant pairs under $A$ is also an $I$-category as a
consequence of \cite{iht} I.3.8 and I.3.6, see also the following
remark. The set of
homotopy classes of maps $(X,Y)\rightarrow (U,V)$ between cofibrant
pairs under $A$ is denoted by $$[(X,Y),(U,V)]^A.$$

\begin{rem}
In \cite{iht} I.3.8 it is claimed that the category of pairs of an
$I$-category with an extra condition satisfied by $\C{Topp}$ is an
$I$-category, but this is not correct since cofibrant pairs are
exactly those pairs for which the map from the initial object is a
cofibration, and this is a requirement of the $I$-category axioms.
However the proof of \cite{iht} I.3.8 is enough to check that the
category of cofibrant pairs is an $I$-category.
\end{rem}

A \emph{based space under $A$} is a cofibrant space $X$ under $A$
together with a retraction $0_X\colon X\r A$ of $i_X$. If $X$ is
based the \emph{trivial map} $0\colon X\r Y$ is the composite
$$X\st{0_X}\To A\st{i_Y}\To Y,$$ and $[X,Y]^A$ is a based set. A map $f\colon X\r Y$ under $A$
between based spaces is a \emph{based map} if $0_X=0_Yf$.
We say that $f$ is \emph{based up to homotopy} if $0_X$ is just homotopic to $0_Yf$. Based spaces and maps are used to
define in the usual way cones $CX$, cofibers $C_f$ of maps $f\colon X\r Y$,
suspensions $\Sigma X$,
and related constructions in $\C{Topp}^A_c$ such as the long
cofiber sequence of a based map
\begin{equation*}
X\st{f}\r Y\st{i_f}\rightarrowtail C_f\st{q}\r\Sigma X\st{\Sigma
f}\r\Sigma Y\st{\Sigma i_f}\rightarrowtail\Sigma C_f\st{\Sigma
q}\r\Sigma^2 X\st{\Sigma^2f}\r\cdots.
\end{equation*}
If the based map $f\colon X\rightarrowtail Y$ is already a
cofibration the cofiber $C_f$ is equivalent to the quotient
space $Y/X$.
One can check by using the ``gluing lemma'' \cite{ah} II.1.2 that up
to homotopy equivalence these constructions only depend on the
homotopy class of the retraction $0_X\colon X\r A$ in
$\C{Topp}^A_c$. Moreover, up to homotopy equivalence the cofiber of
a map $f$ only depends on the homotopy class of $f$ in
$\C{Topp}^A_c$. Furthermore, cylinders, cones and suspensions of
based spaces and cofibers of based maps are canonically based. See
\cite{iht} I.6 and I.7 for further details.

A map $f\colon A\r B$ induces a ``cobase change'' functor
$$f_*\colon \C{Topp}^A_c\To \C{Topp}^B_c$$ given by
$$\xymatrix{A\ar[r]^f\ar@{
>->}[d]\ar@{}[rd]|{\text{push}}&B\ar@{
>->}[d]\\X\ar[r]&f_*X}$$ This functor preserves relative cylinders
$f_*I_AX=I_Bf_*X$, hence it induces a functor between the respective
homotopy categories
\begin{equation}\label{cbase}
f_*\colon \C{Topp}^A_c/\!\simeq\;\To \C{Topp}^B_c/\!\simeq.
\end{equation}
Moreover, up to natural equivalence this last functor only depends
on the homotopy class of $f$ in $\C{Topp}$. In particular it is an
equivalence of categories if $f$ is a homotopy equivalence. This can
also be checked by using the ``gluing lemma''.


So far we have just considered the features of proper homotopy theory which are common to many other homotopy theories. The most basic purely proper homotopy invariant is the space of Freudenthal ends. The space of \emph{Freudenthal ends} $\F(X)$ of a connected locally compact
finite-dimensional $CW$-complex $X$ is the following
inverse limit
$$\F(X)=\lim_{\begin{scriptsize}\begin{array}{c}K\subset X\\\text{compact}\end{array}\end{scriptsize}}\pi_0(X-K).$$
Here $\pi_0(X-K)$ denotes the finite set of connected
components of $X-K$ with the discrete topology. This space is
known to be homeomorphic to a closed subspace of the Cantor set.
In fact any closed subspace of the Cantor set is homeomorphic to
the space of Freudenthal ends of some tree.
Recall that \emph{trees} are contractible $1$-dimensional
$CW$-complexes. In addition in this paper we shall suppose that
all trees are locally compact.

Taking space of Freudenthal ends defines a functor $\F$ from the
proper homotopy category of connected locally compact
finite-dimensional $CW$-complexes. This functor determines an
equivalence between the proper homotopy category of trees and the
category of closed subspaces of the Cantor set, see \cite{iht}
II.1.10. Therefore, given a tree $T$, up to equivalence of
categories $\C{Topp}^T_c$ only depends on the space $\F(T)$.

The \emph{Freudenthal compactification} $\hat{X}$ of a connected locally compact
finite-dimensional $CW$-complex $X$
is, as a set, the disjoint union $\hat{X}=X\sqcup \F(X)$. Given an
open subset $U\subset X$ we define $U^\F$ as the set of
Freudenthal ends $\varepsilon\in\F(X)$ such that $U$ is the
coordinate of $\varepsilon$ in $\pi_0(X-K)$ for some compact
subset $K\subset X$. The sets $U\sqcup U^\F$ form a basis of the
topology of $\hat{X}$. The space $X$ is an open dense subspace of
$\hat{X}$. Moreover, the Freudenthal compactification is a functor
and the inclusions $X\subset \hat{X}\supset\F(X)$ are natural.


Following \cite{iht}, we now describe the reduced
and normalized models of connected locally compact $CW$-complexes in
proper homotopy theory. The building blocks of these models are the
so called spherical objects.

The $n$-dimensional \emph{spherical object} $S^n_\alpha$ under a
tree $T$ associated to a proper map $\alpha \colon A\r T$ with $A$ a
discrete set is the space obtained by pasting an $n$-sphere $S^n$ to
$\alpha(a)$ by the base-point $*\in S^n$ for every $a\in A$
$$\xymatrix{A\ar[r]^\alpha\ar@{ >->}[d]_{(1,*)}\ar@{}[rd]|{\text{push}}&T\ar@{ >->}[d]\\A\times S^n\ar[r]&**[r]
T\cup_\alpha A\times S^n=S^n_\alpha}$$ Spherical objects are based
by the map $S^n_\alpha\r T$ which collapses each sphere to the point
pasted to $T$. One can readily check that $\Sigma
S^n_\alpha=S^{n+1}_\alpha$ $(n\geq 0)$.

A \emph{$T$-complex} is a finite-dimensional $CW$-complex $X$
whose $1$-skeleton is a $1$-dimensional spherical object
$X=S^1_{\alpha_1}$ and such that the $(n+1)$-skeleton $X^{n+1}$
$(n\geq 1)$ is the cofiber of a map $f_{n+1}\colon
S^n_{\alpha_{n+1}}\r X^n$ in $\C{Topp}^T_c$ which is called
\emph{attaching map} of $(n+1)$-cells. This is what Baues and
Quintero call in \cite{iht} a finite-dimensional, reduced and
normalized $CW$-complex relative to $T$.

$T$-complexes $X$ are connected, locally compact, and the inclusion
$T\subset X$ induces a homeomorphism between the spaces of
Freudenthal ends $\F(T)\cong\F(X)$ which will be regarded as an
identification. Conversely for every connected finite-dimensional
locally compact $CW$-complex $Y$ there exists a subtree $T\subset Y$
which induces a homeomorphism on Freudenthal ends $\F(T)\cong\F(Y)$,
and for any such a tree there exists a $T$-complex $X$ and a
homotopy equivalence $X\simeq Y$ in $\C{Topp}^T_c$.

Any $T$-complex $X$ can be regarded as a based space in a unique
way up to homotopy since the set $[X,T]^T$ is a singleton, see
\cite{iht} IV.6.4. In particular all maps under $T$ between
$T$-complexes are based up to homotopy.

The elementary properties of $T$-complexes are similar to those of
connected pointed $CW$-complexes, see \cite{iht} for details. In
particular there is a proper cellular approximation theorem, so if
$\C{CW}^T$ is the category of $T$-complexes and cellular maps under
$T$ the homotopy category $\C{CW}^T/\!\simeq$ is a full subcategory
of $\C{Topp}^T_c/\!\simeq$. Also, as it happened with
$\C{Topp}^T_c/\!\simeq$, up to equivalence of categories
$\C{CW}^T/\!\simeq$ only depends on the proper homotopy type of the
tree $T$, or equivalently on its space of Freudenthal ends $\F(T)$.
Therefore the first step in the study of proper homotopy types of
connected finite-dimensional locally compact $CW$-complexes is to
separate them attending to their most basic homotopy invariant, the
space of Freudenthal ends, and then to concentrate on the study of
the categories of those with a fixed space of Freudenthal ends, or
equivalently on the categories $\C{CW}^T$.

As it is well-known to experts, there are drastic differences
between the homotopy theories of the categories $\C{CW}^T$ when
$T$ is or not compact. In fact if $T$ is any compact tree
$\F(T)=\emptyset$ is the empty set and $\C{CW}^T/\!\simeq$ is
equivalent to the ordinary homotopy category of connected compact
pointed $CW$-complexes. In this paper we are concerned with
the non-compact case so, unless we state the contrary, all trees
considered will be non-compact.

\section{Controlled algebra}\label{ch2}\label{modules}\label{freecontrol}


Recall the category $\mo(\C{A})$ of (right) \emph{modules} over a
small additive category $\C{A}$ is the abelian category of
additive functors $\um{M}\colon\C{A}^\op\r\C{Ab}$ and natural
transformations. Here $\C{Ab}$ is the category of abelian groups.
There is a full Yoneda inclusion of categories
\begin{equation}\label{yoneda}
\mathrm{Yoneda}\colon\C{A}\hookrightarrow\mo(\C{A})
\end{equation} which sends an object $X$ to the representable
functor $\hom_\C{A}(-,X)$. These modules are said to be
\emph{finitely generated (f. g.) free} and constitute a set of small
projective generators of $\mo(\C{A})$. For simplicity we shall
identify $X$ with $\hom_\C{A}(-,X)$. An $\C{A}$-module is
\emph{finitely generated (f. g.)} if it is the image of a f. g. free
$\C{A}$-module, and it is \emph{finitely presented (f. p.)} if it is
the cokernel of a morphism between f. g. free $\C{A}$-modules.

If $\C{B}$ is another small additive category an additive functor
$\uf{F}\colon\C{A}\r\C{B}$ induces two ``change of coefficient''
functors
$$\uf{F}^*\colon\mo(\C{B})\To\mo(\C{A}),$$
$$\uf{F}_*\colon\mo(\C{A})\To\mo(\C{B}).$$
The first one is exact and preserves colimits, it is given by
right composition with $\uf{F}$, $\uf{F}^*\um{M}=\um{M}\uf{F}$,
and the second one is characterized by being left adjoint to
$\uf{F}^*$, in particular $\uf{F}_*$ is right-exact and more
generally colimit-preserving, see \cite{hilton} II.7.7. Moreover,
$\uf{F}_*$ fits into the following commutative diagram
\begin{equation}\label{yonedacom}
\xymatrix{\C{A}\ar[r]^{\uf{F}}\ar@{^{(}->}[d]_{\text{Yoneda}}&
\C{B}\ar@{^{(}->}[d]^{\text{Yoneda}}\\
\mo(\C{A})\ar[r]^{\uf{F}_*}&\mo(\C{B})}
\end{equation}
hence $\uf{F}_*$ restricts to the full subcategories of f. p.
modules. Furthermore, if $\uf{F}$ is full and faithful $\uf{F}_*$
is too, and in this case $\uf{F}^*\uf{F}_*$ is naturally
equivalent to the identity, see \cite{borceux1} 3.4.1. This can be
directly checked by using that any module has a projective
resolution by (arbitrary) direct sums of f. g. free ones.

It is interesting to know when the category of modules over a small
additive category $\C{A}$ is equivalent as an abelian category to
the category $\mo(R)$ of modules over a ring $R$, i. e. when $\C{A}$
is Morita equivalent to a ring. The typical case in this paper will
be the following: if there exists an object $X$ of $\C{A}$ such that
any other object is a retract of $X$ then $\C{A}$ is Morita
equivalent to the endomorphism ring of $X$, that is
$\endo_\C{A}(X)=\hom_\C{A}(X,X)$. This follows from \cite{rso} 8.1.
The equivalence of categories is given by the evaluation functor
$$ev_X\colon\mo(\C{A})\To\mo(\endo_\C{A}(X))\colon \um{M}\mapsto\um{M}(X).$$

If $Y$ is an object of $\C{B}$ the morphism set
$\hom_\C{B}(Y,\uf{F}X)$ is a
left-$\endo_\C{A}(X)$-right-$\endo_\C{B}(Y)$-bimodule and one can
check that the next diagram commutes up to natural equivalence
\begin{equation}\label{conmotra}
\xymatrix@C=100pt{\mo(\C{A})\ar[r]^{\uf{F}_*}\ar[d]_{ev_X}&\mo(\C{B})\ar[d]^{ev_Y}\\
\mo(\endo_\C{A}(X))\ar[r]_{-\mathop{\otimes}\limits_{\endo_\C{A}(X)}\hom_\C{B}(Y,\uf{F}X)}&
\mo(\endo_\C{B}(Y))}
\end{equation}


Let $T$ be a tree and $R$ a commutative ring. A \emph{free
$T$-controlled $R$-module} is a pair $R\grupo{A}_\alpha$ formed by a
free $R$ module $R\grupo{A}$ with basis a discrete set $A$ together
with a proper map $\alpha\colon A\r T$ which is called the
\emph{height function}. If we drop the properness condition for
$\alpha$ we say that $R\grupo{A}_\alpha$ is a \emph{big free
$T$-controlled $R$-module}. A \emph{controlled homomorphism}
$\varphi\colon R\grupo{A}_\alpha\r R\grupo{B}_\beta$ is a
homomorphism between the underlying $R$-modules such that for any
neighbourhood $U$ of $\varepsilon\in\F(T)$ in $\hat{T}$ there exists
another one $V\subset U$ with $\varphi(\alpha^{-1}(V))\subset
R\grupo{\beta^{-1}(U)}$. The category $\C{M}^b_R(T)$ of possibly big
free $T$-controlled $R$-modules and controlled homomorphisms is an
additive category, and the full additive subcategory $\C{M}_R(T)$ of
(non-big) free $T$-controlled $R$-modules is small.

\begin{rem}
The category $\C{M}_R(T)$ was introduced in \cite{iht}. Here we give
a slightly different but equivalent definition. For $R=\Z$ it is
also isomorphic to the category of finitely generated free $T$-trees
of abelian groups considered in \cite{projdim}, this observation can
be regarded as a particular case of \cite{iht} VIII.3.5, see also
\cite{iht} V.3.10. In \cite{projdim} we only consider trees $T$ with
a base vertex $v_0$ and no leaves apart from possibly $v_0$, however
this is not a restriction since any tree is proper homotopy
equivalent to such a kind of tree. The representation theory of f.
p. $\C{M}_R(T)$-modules has been considered by us in \cite{rca}, see
Section \ref{rrt} for a review on some of its results. The category
$\C{M}_R(T)$ corresponds to $\C{M}_R(\bar{T})$ in \cite{rca} for
$\bar{T}=(\hat{T},T,\F(T))$.
\end{rem}

The \emph{support} of a free $T$-controlled $R$-module
$R\grupo{A}_\alpha$ is the derived set of $\alpha(A)$ in $\hat{T}$,
that is the subset $\alpha(A)'\subset\hat{T}$ formed by those points
which contain infinitely many points of $\alpha(A)$ in any
neighborhood. One can readily check that the support is always
formed by Freudenthal ends $\alpha(A)'\subset\F(T)$. Moreover, the
support of a finite direct sum is the union of the supports of the
factors.

The isomorphism class of a free $T$-controlled $R$-module is
completely determined by its support and the rank of the underlying
$R$-module, see \cite{rca} 3.2 and notice that we here assume that
$R$ is commutative. Notice also that the rank of the underlying
$R$-module is $\aleph_0$ if and only if the support is non-empty.

Any free $T$-controlled $R$-module is a retract of any other one
with support $\F(T)$, see \cite{rca} 3.5. The canonical free
$T$-controlled $R$-module with support $\F(T)$ is
$R\grupo{T^0}_\delta$ with $\delta\colon T^0\hookrightarrow T$
the inclusion of the vertex set. In particular, as we observed in
above, $\C{M}_R(T)$ is Morita equivalent to the
endomorphism ring of a free $T$-controlled $R$-module with support
$\F(T)$. Up to isomorphism this ring (in fact $R$-algebra) only
depends on $\F(T)$ and will be denoted $R(\F(T))$, see \cite{rca}
3.4. An explicit equivalence of module categories is given by the
evaluation functor
\begin{equation}\label{expli}
ev_{R\grupo{T^0}_\delta}\colon\mo(\C{M}_R(T))\To\mo(R(\F(T))).
\end{equation}

A curious property of the rings $R(\F(T))$ is stated in the
following proposition.

\begin{prop}\label{solo1}
Any finitely generated free $R(\F(T))$-module is isomorphic to
just one copy of $R(\F(T))$.
\end{prop}

The proposition follows from the previous equivalence and the fact
that the support of a finite direct sum of free $T$-controlled
$R$-modules with support $\F(T)$ has support $\F(T)$ as well.

\begin{prop}\label{induce}
A proper map between trees $f\colon T\rightarrow T'$ induces a
``change of tree'' additive functor
$$\uf{F}^f\colon\C{M}^b_R(T)\To\C{M}^b_R(T'),\;\;\uf{F}^fR\grupo{A}_\alpha=R\grupo{A}_{f\alpha}.$$
This functor restricts to the full subcategories of non-big objects.
\end{prop}

This proposition follows from the following
sequential characterization of controlled homomorphisms. 

\begin{lem}\label{critter}
Given two possibly big free $T$-controlled $R$-modules
$R\grupo{A}_\alpha$ and $ R\grupo{B}_\beta$, a homomorphism
$\varphi\colon R\grupo{A}\r R\grupo{B}$ induces a controlled
homomorphism $\varphi\colon R\grupo{A}_\alpha\r R\grupo{B}_\beta$ if
and only if 
given a sequence $\set{a_n}_{n\geq 0}$ in $A$ such that
$\lim\limits_{n\r\infty}\alpha(a_n)=\varepsilon\in\F(T)$ in
$\hat{T}$, if $b_n$ appears with non-trivial coefficient in the
linear expansion of $\varphi(a_n)$ then
$\lim\limits_{n\r\infty}\beta(b_n)=\varepsilon$.
\end{lem}

This lemma is an
extension of \cite{iht} III.4.14 to the possibly big case. The proof is analogous.

Moreover, as we show in the next proposition, up to natural
equivalence $\uf{F}^f$ only depends on the proper homotopy class of
$f$.

\begin{prop}\label{homono}
If two proper maps between trees $f,f'\colon T\r T'$ are homotopic
then there is a natural equivalence $\uf{F}^f\simeq\uf{F}^{f'}$.
\end{prop}

\begin{proof}
For any possibly big free
$T$-controlled $R$-module $R\grupo{A}_\alpha$ the identity in
$R\grupo{A}$ induces a controlled isomorphism
$R\grupo{A}_{f\alpha}\simeq R\grupo{A}_{f'\alpha}$. This can be checked by using Lemma \ref{critter} and the fact that $\F(f)=\F(f')$.
\end{proof}

This proposition shows that the categories $\C{M}^b_R(T)$ and
$\C{M}_R(T)$ only depend, up to equivalence of additive categories,
on the proper homotopy type of $T$, or equivalently on $\F(T)$.

The functor $\uf{F}^f$ is obviously faithful. It is not full in
general, but it is when $f$ induces an injection in Freudenthal
ends.

\begin{prop}\label{fullis}
If the proper map between trees $f\colon T\r T'$ induces an
injection $\F(f)\colon\F(T)\hookrightarrow\F(T')$ then $\uf{F}^f$
is full.
\end{prop}

This proposition follows from Lemma \ref{critter}.

The fully faithful functor
$\C{M}_R(T)\hookrightarrow\mo(\C{M}_R(T))$ in (\ref{yoneda}) can be
extended to the whole category of possibly big free $T$-controlled
$R$-modules
\begin{equation}\label{yonedabig}
\C{M}^b_R(T)\To\mo(\C{M}_R(T))\colon R\grupo{A}_\alpha\mapsto
\hom_{\C{M}_R^b(T)}(-,R\grupo{A}_\alpha)_{|_{\C{M}_R^\op(T)}}.
\end{equation}
In this way we can regard a big free $T$-controlled $R$-module
$R\grupo{A}_\alpha$ as an $\C{M}_R(T)$-module. This extension is
easily seen to be faithful but not necessarily full, and also
satisfies the following property.

\begin{prop}\label{extigrande}\label{EXTIGRANDE}
Given a proper map between two trees $f\colon T\r T'$ the following
diagram of functors is commutative up to natural equivalence
$$\xymatrix{\C{M}_R^b(T)\ar[r]^{\uf{F}^f}\ar@{^{(}->}[d]&\C{M}_R^b(T')\ar@{^{(}->}[d]\\
\mo(\C{M}_R(T))\ar[r]^{\uf{F}^f_*}&\mo(\C{M}_R(T'))}$$
\end{prop}

A proof of this proposition will be given in Appendix
\ref{proofofextigrande}.

Another crucial property of ``change of tree'' functors is stated in
the next proposition.

\begin{prop}\label{esexacto}\label{ESEXACTO}
Given a proper map between two trees $f\colon T\r T'$ such that
$\F(f)\colon\F(T)\hookrightarrow\F(T')$ is injective, the change of
coefficients
$$\uf{F}^f_*\colon\mo(\C{M}_R(T))\To\mo(\M_R(T'))$$
is exact.
\end{prop}

The proof of this proposition is in Appendix
\ref{proofofesexacto}.

We will be specially interested in the case $R=\Z$ the integers.
In order to keep the notation close to \cite{iht} we will write
$\ab(T)$ and $\ab^b(T)$ for $\C{M}_\Z(T)$ and $\C{M}^b_\Z(T)$
respectively.

A crucial property of the category of $\ab(T)$-modules is the
following.

\begin{thm}\label{kernel}
The kernel of a morphism between f. g. free $\ab(T)$-modules is f.
g. free.
\end{thm}

This follows from \cite{projdim} 3.2. A translation of the result
in \cite{projdim} to the language used here for the special case
$T=\Real_+=[0,+\infty)$ the half line appears in \cite{iht} V.5.3.
Some consequences of this theorem are the following.

\begin{cor}\label{pd2}
F. p. $\ab(T)$-modules have projective dimension $\leq 2$.
\end{cor}

\begin{cor}\label{extcol}
If $\um{M}$ is a f. p. $\ab(T)$-module the functors
$\ext^n_{\ab(T)}(\um{M},-)$ preserve colimits $(n\geq 0)$.
\end{cor}

\begin{cor}
Every f. g. projective $\ab(T)$-module is f. g. free.
\end{cor}

This is a consequence of the standard fact that a f. g. projective
$\ab(T)$-module is the kernel of an endomorphism of a f. g. free
$\ab(T)$-module.

\begin{cor}
$\ab(T)$ is equivalent to the category of f. g. projective
$\Z(\F(T))$-modules.
\end{cor}

This is a consequence of the previous corollary and the
equivalence of abelian categories in (\ref{expli}).

Another relevant ring for us will be the field with two elements
$\uf{F}_2$. For simplicity we write $\C{vect}^b(T)$ and
$\C{vect}(T)$ for the categories $\C{M}^b_{\uf{F}_2}(T)$ and
$\C{M}_{\uf{F}_2}(T)$, respectively. The tensor product by $\Z/2$
induces a full functor
$$-\otimes\Z/2\colon\ab^b(T)\To\C{vect}^g(T)$$ which  restricts to
the full subcategories of non-big objects. Moreover, if we denote
$\ab^b(T)\otimes\Z/2$ and $\ab(T)\otimes\Z/2$ to the categories
obtained from $\ab^b(T)$ and $\ab(T)$ by tensoring their morphism
abelian groups by $\Z/2$, respectively, it is easy to check that
the following proposition holds.
\begin{prop}\label{yoqse}
The functor $-\otimes\Z/2$ gives rise to isomorphisms
$\ab^b(T)\otimes\Z/2\simeq\C{vect}^b(T)$ and
$\ab(T)\otimes\Z/2\simeq\C{vect}(T)$.
\end{prop}

The functor $-\otimes\Z/2\colon\ab(T)\r\C{vect}(T)$ induces two
``change of coefficient'' functors
\begin{equation}\label{amod2}
\begin{array}{c}
(-\otimes\Z/2)^*\colon\mo(\C{vect}(T))\To\mo(\ab(T)),\\{}\\
(-\otimes\Z/2)_*\colon\mo(\ab(T))\To\mo(\C{vect}(T)).
\end{array}
\end{equation}

The elementary properties of these functors are summarized in the following proposition.

\begin{prop}\label{acabalo}
The composite $(-\otimes\Z/2)_*(-\otimes\Z/2)^*$ is naturally equivalent to the identity, and $(-\otimes\Z/2)^*(-\otimes\Z/2)_*$ is naturally
equivalent to the functor
$$-\otimes\Z/2\colon\mo(\ab(T))\To\mo(\ab(T))$$ given by
left-composition with $-\otimes\Z/2\colon\Ab\r\Ab$.
\end{prop}

Notice that now we can identify $\mo(\C{vect}(T))$ with the full
subcategory of $\mo(\ab(T))$ formed by those $\ab(T)$-modules which
take values in abelian groups with exponent $2$. Moreover, as a
consequence of Propositions \ref{yoqse} and \ref{acabalo} we obtain the
following result.

\begin{cor}\label{calculete}
For any possibly big free $T$-controlled $\Z$-module
$\Z\grupo{A}_\alpha$ we have a natural identification
$\Z\grupo{A}_\alpha\otimes\Z/2=\uf{F}_2\grupo{A}_\alpha$ as
$\ab(T)$-modules.
\end{cor}

\section{Elementary proper homotopy algebraic
invariants}\label{ephai}\label{ctpht}


Let $\C{S}^n(T)$ be the proper homotopy category of $n$-dimensional
spherical objects under $T$. One of the main tools to construct
algebraic invariants in proper homotopy theory is the following
result.

\begin{prop}\label{fund}
There are isomorphisms of categories $(n\geq 2)$
$$\C{S}^n(T)\st{\simeq}\To\ab(T)\colon
S^n_\alpha\mapsto\Z\grupo{A}_\alpha$$ which are compatible with
the suspension functors
$\Sigma\colon\C{S}^n(T)\rightarrow\C{S}^{n+1}(T)$.
\end{prop}

These isomorphisms are given by the ordinary homotopy groups
$\pi_n(S^n_\alpha,T)=\Z\grupo{A}$ in the level of underlying
abelian groups. See \cite{iht} II.4.15 and II.4.20 for a proof.
From now on we shall use these isomorphisms as identifications.

\begin{defn}
The \emph{proper homotopy $\ab(T)$-modules} of a cofibrant space $X$
under $T$ are $$\Pi_nX=[-,X]^T\colon \C{S}^n(T)^\op\To\Ab,\;\; n\geq
2.$$ The proper homotopy modules of a cofibrant pair $(X,Y)$ under
$T$ are also defined
$$\Pi_{n+1}(X,Y)=[(C-,-),(X,Y)]^T\colon\C{S}^n(T)^\op\To\Ab,\;\;
n\geq 2.$$

These proper homotopy modules are functors in the corresponding
homotopy categories. Moreover, there is a natural long exact
sequence
$$\cdots\r\Pi_{n+1}X\st{j}\r\Pi_{n+1}(X,Y)\st{\partial}\r\Pi_nY\r\Pi_nX\r\cdots.$$

A cofibrant space $X$ under $T$ is \emph{$0$-connected} if the map
$T\rightarrowtail X$ induces a surjection
$\F(T)\twoheadrightarrow\F(X)$ between the spaces of Freudenthal
ends, and \emph{$n$-connected} $(n\geq 1)$ if in addition
$[S^1_\alpha,X]^T=0$ is a singleton for every $1$-dimensional
spherical object $S^1_\alpha$ and $\Pi_kX=0$ for all $2\leq k\leq
n$.

Similarly, a cofibrant pair $(X,Y)$ is \emph{$0$-connected} if the
map $Y\rightarrowtail X$ induces an epimorphism
$\F(Y)\twoheadrightarrow\F(X)$, \emph{$1$-connected} if
$Y\rightarrowtail X$ also induces an epimorphism
$[S^1_\alpha,Y]^T\twoheadrightarrow[S^1_\alpha,X]^T$ for all
$1$-dimensional spherical object $S^1_\alpha$, and $n$-connected
$(n\geq 2)$ if in addition $[(CS^1_\alpha,S^1_\alpha),(X,Y)]^T=0$ is
always a singleton and $\Pi_k(X,Y)=0$ for all $3\leq k\leq n$.

The \emph{Whitehead modules} of a $T$-complex $X$ are
$$\Gam_nX=\im[\Pi_nX^{n-1}\rightarrow\Pi_nX^n],\;\; n\geq 2.$$
\end{defn}


\begin{defn}
The \emph{cellular chain complex} $\cc_*X$ of a $T$-complex $X$
with $X^1=S^1_{\alpha_1}$ and attaching maps $f_{n+1}\colon
S^n_{\alpha_{n+1}}\rightarrow X^n$ $(n\geq 1)$ is the positive
chain complex in $\ab(T)$ defined as
$\cc_nX=\Z\grupo{A_n}_{\alpha_n}$ with lower differential
$d_2\colon\cc_2X\rightarrow\cc_1X$ given by $\Sigma f_2$ and
higher differentials $d_{n+1}\colon\cc_{n+1}X\rightarrow\cc_nX$
given by the following composites $(n\geq 2)$
$$S^n_{\alpha_{n+1}}\st{f_{n+1}}\To X^n\To X^n/X^{n-1}=S^n_{\alpha_n}.$$

The cellular chain complex is a functor from the category of
$T$-complexes to the category of bounded chain complexes in
$\ab(T)$ concentrated in positive degrees
$$\cc_*\colon\C{CW}^T\To\C{chain}(\ab(T)).$$ The category
$\C{chain}(\ab(T))$ is also an $I$-category with the usual
cylinder of chain complexes, compare \cite{cfhh} III.9.2, and
$\cc_*$ preserves all the structure. In particular $\cc_*$ factors
trough the respective homotopy categories.

The homology of a chain complex in $\ab(T)$ is defined via the
Yoneda inclusion (\ref{yoneda}), in particular one can define the
\emph{proper homology $\ab(T)$-modules} of $X$
$$\hh_nX=H_n\cc_*X.$$
The \emph{proper cohomology groups} of $X$ with coefficients in an
$\ab(T)$-module $\um{M}$ are
$$H^n(X,\um{M})=H^n\hom_{\ab(T)}(\cc_*X,\um{M}).$$
\end{defn}

One of the most relevant consequences of Theorem \ref{kernel} in
proper homotopy theory is the following.

\begin{prop}
The proper homology modules of a $T$-complex are finitely presented.
\end{prop}


Several classical theorems generalize to the proper homotopy invariants defined above. In particular there is a homological Whitehead theorem (\cite{iht} VI.6) and a Blakers-Massey's excision theorem (\cite{iht} VI.7.3). This Blakers-Massey theorem can be used to prove a Freudenthal suspension theorem for proper homotopy modules. Nevertheless its
main consequence is the fact that the category $\C{Topp}^T$ is a homological cofibration category under $\C{S}^1(T)$ in the sense of \cite{cfhh} V.1.1., see \cite{cfhh} V.1.2. Homological cofibration categories are cofibration categories with a
certain class of objects, satisfying certain axioms, which are used,
together with their suspensions, to construct the analogue of
$CW$-complexes. This extra structure allows to develop a rich
obstruction theory for these complexes in an abstract setting. We
will review some of this obstruction theory in Section \ref{hs} for
the particular example $\C{Topp}^T$, see \cite{cfhh} for further
details.

A $T$-complex $X$ is \emph{$n$-reduced} $(n\geq 1)$ if $X^n=T$, the
relation of this concept with the connectivity is established in the
next proposition.

\begin{prop}\label{conexred}
Given an $n$-connected $T$-complex $X$ there exists a homotopy
equivalent one $Y$ which is $n$-reduced. Moreover, one can assume
that the dimension of $Y$ does not exceed the maximum between the
dimension of $X$ and $n+2$.
\end{prop}

The first part of the statement is a particular case of \cite{iht}
IV.7.2. In fact the techniques used there are go back to classical
results by J. H. C. Whitehead adapted to proper homotopy theory.
Following carefully the proof of \cite{iht} IV.7.2 one notices that
the second part of the proposition is also proven there.

\begin{cor}\label{gammaes0}
If $X$ is an $(n-1)$-connected $T$-complex $\Gam_kX=0$ for all
$2\leq k\leq n$.
\end{cor}

For any $T$-complex $X$ there is a sequence of $\ab(T)$-modules
\begin{flushleft}$\cdots\st{j}\r\Pi_{n+1}(X^{n+1},X^n)\st{\partial}\r\Pi_nX^n\st{j}\r\Pi_n(X^n,X^{n-1})\st{\partial}\r\cdots$\end{flushleft}
\begin{flushright}$\cdots\st{j}\r\Pi_{3}(X^{3},X^2)\st{\partial}\r\Pi_2X^2\st{=}\r\Pi_2X^2\r 0,$\end{flushright}
which is exact in the proper homotopy modules of all pairs and also
in the second $\Pi_2X^2$. If $X$ is $1$-reduced the cellular chain
complex $\cc_*X$ coincides with the following chain complex of
$\ab(T)$-modules, see \cite{iht} VI.5.5,
$$\cdots\r\Pi_{n+1}(X^{n+1},X^n)\st{j\partial}\r\Pi_n(X^n,X^{n-1})\r\cdots
\r\Pi_{3}(X^{3},X^2)\st{\partial}\r\Pi_2X^2\rightarrow 0,$$ since
in this case
$\Pi_n(X^n,X^{n-1})=\Pi_nX^n/X^{n-1}=\Pi_nS^n_{\alpha_n}=\cc_nX$ for
all $n\geq 2$ so one can construct a long exact sequence involving
proper homotopy, homology and Whitehead modules as J. H. C.
Whitehead did in \cite{aces}, see \cite{iht} VI.5.4,
$$\cdots\r\hh_{n+1}X\st{b_{n+1}}\r\Gam_nX\st{i_n}\r\Pi_nX\st{h_n}\r\hh_nX\r\cdots\r\Pi_2X\st{h_2}\twoheadrightarrow\hh_2X.$$
The \emph{Hurewicz morphisms} in proper homotopy theory are defined
to be $h_n$ $(n\geq 2)$. Moreover, following Whitehead's
terminology we call the $b_n$ $(n\geq 3)$ \emph{secondary boundary
operators}.

The -r6-er Hurewicz theorem follows fron Whitehead's long exact sequence in proper homotopy theory,
Proposition \ref{conexred} and Corollary \ref{gammaes0}.

\section{$1$-dimensional spherical objects}

The proper homotopy category of $1$-dimensional spherical objects
$\C{S}^1(T)$ has a purely algebraic model, as it happens in higher
dimensions, see Proposition \ref{fund}. It is given by the following
non-abelian version of the free $T$-controlled modules introduced in
Section \ref{freecontrol}.

\begin{defn}\label{fcnag}
A \emph{free $T$-controlled group} is a pair $\grupo{A}_\alpha$
where $\grupo{A}$ is the free group with basis the discrete set $A$
and $\alpha\colon A\r T$ is a proper map, the \emph{height
function}. A \emph{controlled homomorphism}
$\alpha\colon\grupo{A}_\alpha\r\grupo{B}_\beta$ is a homomorphism
between the underlying groups such that for any neighborhood $U$ of
$\varepsilon\in\F(T)$ in $\hat{T}$ there exists another one
$V\subset U$ with
$\varphi(\alpha^{-1}(V))\subset\grupo{\beta^{-1}(U)}$. Let $\gr(T)$
be the category of free $T$-controlled groups and controlled
homomorphisms. There is an obvious abelianization full functor
$$ab\colon\gr(T)\To\ab(T)\colon\grupo{A}_\alpha\mapsto\Z\grupo{A}_\alpha.$$
\end{defn}

\begin{prop}\label{fund2}
There is an isomorphism of categories
$$\C{S}^1(T)\st{\simeq}\To\gr(T)\colon S^1_\alpha\mapsto\grupo{A}_\alpha$$
which makes commutative the following diagram
$$\xymatrix{\C{S}^1(T)\ar[r]^\S&\C{S}^2(T)\\\gr(T)\ar[u]^\simeq\ar[r]^{ab}&\ab(T)\ar[u]_\simeq}$$
\end{prop}

The isomorphism is induced by the ordinary homotopy groups
$\pi_1(S^1_\alpha,T)=\grupo{A}$ in the underlying groups, and the
isomorphism in the left of the square is in Proposition \ref{fund}.
See \cite{iht} 4.20 for a proof.

\chapter[Quadratic functors in controlled algebra\dots]{Quadratic functors in controlled algebra and proper homotopy
theory}\label{qf}

In this chapter we first recall the basic definitions and elemetary examples of quadratic functors. In the Section \ref{eqf} we construct quadratic functors in controlled algebra which play an important role in this paper. These functors are applied in Section \ref{slwm} to the computation of proper homotopy invariants. In the Section \ref{f2} we use the controlled exterior square to describe the category of free controlled groups of nilpotency class $2$ as a linear extension. Finally in Section \ref{chV} we give an account of the obstruction theory in proper homotopy theory developed in \cite{iht} from the perspective of the tower of categories in \cite{cfhh} and we establish the links with the cohomology of categories. For this we use the controlled quadratic functors defined in the Section \ref{eqf} and the computations in Section \ref{slwm}.

\section{General quadratic functors}\label{gqf}

\begin{defn}\label{dgqf}
A function between abelian groups $f\colon A\rightarrow B$ is
\emph{quadratic} if the function $$[-,-]\colon A\times A\rightarrow
B\colon (a,b)\mapsto [a,b]=f(a+b)-f(a)-f(b)$$ is bilinear.

A functor $F\colon\C{A}\rightarrow\C{B}$ between additive categories
is a \emph{quadratic functor} if it induces quadratic functions
between morphism abelian groups.

Suppose that $\C{B}$ is abelian. Let $X$ and $Y$ be two objects in
$\C{A}$ and
\begin{equation}\label{dsd}
\xymatrix{X\ar@<-.5ex>[r]_<(.3){i_1}&X\oplus
Y\ar@<-.5ex>[l]_<(.3){p_1}\ar@<.5ex>[r]^<(.37){p_2}&Y
\ar@<.5ex>[l]^<(.25){i_2}}
\end{equation}
the projections and inclusions of the factors of the direct sum.
The \emph{quadratic crossed effect} $F(X|Y)$ is defined as
$$F(X|Y)=\im[[i_1p_1,i_2p_2]\colon F(X\oplus Y)\rightarrow
F(X\oplus Y)].$$ In particular $[i_1p_1,i_2p_2]$ factors as
follows
$$[i_1p_1,i_2p_2]\colon F(X\oplus Y)\st{r_{12}}\twoheadrightarrow
F(X|Y)\st{i_{12}}\hookrightarrow F(X\oplus Y),$$ and there is an
isomorphism $$F(X)\oplus F(X|Y)\oplus F(Y)\simeq F(X\oplus Y).$$
Since $F$ is quadratic $F(-|-)$ is a biadditive functor. Moreover,
$F(-|-)$ vanishes if and only if $F$ is additive.
\end{defn}

We are specially interested in the following classical quadratic
functors from abelian groups to abelian groups:
\begin{itemize}
\item Whitehead's functor $\Gamma$,

\item exterior square $\wedge^2$,

\item tensor square $\otimes^2$,

\item and reduced tensor square $\hat{\otimes}^2$.
\end{itemize}

The functor $\Gamma$ is characterized by the existence of a
natural function $\gamma\colon A\rightarrow\Gamma A$ which is
universal among all quadratic functions $f\colon A\rightarrow B$
with $f(a)=f(-a)$ for all $a\in A$. The exterior square $\wedge^2$
and the tensor square $\otimes^2$ are fairly well-known, and the
reduced tensor square is determined by the push-out square in the
following commutative diagram with exact rows and column
\begin{equation}\label{d1}
\xymatrix{\otimes^2A\ar[d]^{[-,-]}&&\\\Gamma
A\ar[r]^\tau\ar@{->>}[d]^\sigma\ar@{}[dr]|{\mathrm{push}}&\otimes^2A\ar@{->>}[d]^{\bar{\sigma}}\ar@{->>}[r]^q&\wedge^2A\ar@{=}[d]\\
A\otimes\Z_2\ar@{^{(}->}[r]^{\bar{\tau}}&\hat{\otimes}^2A\ar@{->>}[r]^{\bar{q}}&\wedge^2A}
\end{equation}
Here $\tau(\gamma(a))=a\otimes a$, $q(a\otimes b)=a\wedge b$,
$\sigma(\gamma(a))=a\otimes 1$ and $\bar{\sigma}(a\otimes
b)=a\hat{\otimes}b$. One can check that $\bar{\tau}$ is always a
monomorphism by using its naturality, the natural projection
$A\twoheadrightarrow A\otimes\Z/2$, and the fact that all vector
spaces (over the field with $2$ elements) have a basis, but it is
not completely trivial, i. e. it is a not so easy consequence of its
definition.

The quadratic crossed effect of $\Gamma$, $\wedge^2$ and
$\hat{\otimes}^2$ is the tensor product $\otimes$ and the structure
homomorphisms $i_{12}$ and $r_{12}$ are given as follows,
$$\begin{array}{c} \Gamma(A\oplus B)\st{r_{12}}\twoheadrightarrow
A\otimes B\st{i_{12}}\hookrightarrow\Gamma(A\oplus B),\\{}\\
r_{12}(\gamma(i_1a+i_2b))=a\otimes b,\;\;\;i_{12}(a\otimes
b)=[i_1a,i_2,b];\\{}\\{}\\ \wedge^2(A\oplus
B)\st{r_{12}}\twoheadrightarrow A\otimes
B\st{i_{12}}\hookrightarrow\wedge^2(A\oplus B),\\{}\\ r_{12}((i_1
a+i_2b)\wedge(i_1a'+i_2b'))= i_1a\otimes i_2b'-i_1a'\otimes
i_2b,\;\;\; i_{12}(a\otimes b)=i_1a\wedge i_2b;\\{}\\{}\\
\hat{\otimes}^2(A\oplus B)\st{r_{12}}\twoheadrightarrow A\otimes
B\st{i_{12}}\hookrightarrow\hat{\otimes}^2(A\oplus B),\\{}\\
r_{12}((i_1 a+i_2b)\hat{\otimes}(i_1a'+i_2b'))=
i_1a{\otimes}i_2b'-i_1a'{\otimes}i_2b,\;\;\; i_{12}(a\otimes
b)=i_1a\hat{\otimes}i_2b.
\end{array}$$
The tensor square $\otimes^2$ is the square of a biadditive functor,
the tensor product. All quadratic functors arising in this way have
a quite obvious quadratic crossed effect, in particular for
$\otimes^2$ we have
$$\otimes^2(A|B)=A\otimes B\oplus B\otimes A.$$
One can notice that $[-,-]=\Gamma(1,1)i_{12}$ and
$\tau=r_{12}\Gamma(i_1+i_2)$, therefore the whole diagram (\ref{d1})
is determined by the functor $\Gamma$.

The functor $\otimes^2$ sends free abelian groups to free abelian
groups, since the tensor product of two free abelian groups is
free abelian $\Z\grupo{A}\otimes\Z\grupo{B}\simeq\Z\grupo{A\times
B}$. The functors $\wedge^2$ and $\Gamma$ also preserve free abelian
groups. In order to describe bases of $\wedge^2\Z\grupo{E}$ and
$\Gamma\Z\grupo{E}$ we choose a total ordering $\preceq$ in $E$. The
bases of $\wedge^2\Z\grupo{E}$ and $\Gamma\Z\grupo{E}$ are the sets
$\wedge^2E$ and $\Gamma E$, respectively, defined as follows,
\begin{itemize}
\item $\Gamma E=\set{\gamma(e),[e_1,e_2]\,;\, e, e_1, e_2\in
E, e_1\succ e_2}$,

\item $\wedge^2E=\set{e_1\wedge e_2\,;\,e_1\prec e_2\in E}$.
\end{itemize}

\section[Extension of classical quadratic functors\dots]{Extension of classical quadratic functors to controlled
algebra}\label{eqf}

In this section we define controlled generalizations of the classical quadratic functors considered in the previous section. The consistency of these new constructions depends on topological properties of trees that we now recall.

\begin{defn}\label{topgen}
Recall that an \emph{arc} in a topological space is a subspace
homeomorphic to a compact interval. Any two points in a tree
$u,v\in T$ can be joined by a unique arc denoted by $\ol{uv}$.
Moreover, if $w\in\ol{uv}$ the union of $\ol{uw}$ and $\ol{wv}$ by
their common point is $\ol{uv}$.

It is well known that there is a unique metric $d$ in $T$ such
that de distance between two adjacent vertices is $1$ and every
arc $\ol{uv}$ is isometric to the interval
$[0,d(u,v)]\subset\Real$.

Fixed a base-point $v_0\in T$ we define the map
$$\ell\colon T\times T\To T$$ as the unique one which satisfies $\overline{uv}\cap\overline{v_0\ell(u,v)}=\set{\ell(u,v)}$ for any $u,v\in T$.
This map is obviously symmetric $\ell(u,v)=\ell(v,u)$.

\end{defn}

\begin{lem}
The map $\ell$ is well defined.
\end{lem}

\begin{proof}
Given $u,v\in T$ there exists at least one $w\in T$ such that
$\overline{uv}\cap\overline{v_0w}=\set{w}$. We can take $w$ to be
a point in the compact space $\overline{uv}$ where the function
$\overline{uv}\r\Real\colon x\mapsto d(v_0,x)$ reaches the minimum
value. Moreover, $w$ is the unique point with this property.
Suppose by the contrary that there where two points $w_1,w_2\in T$
with $\overline{uv}\cap\overline{v_0w_i}=\set{w_i}$ $(i=1,2)$.
Then
$\overline{v_0w_1}\cup\overline{w_1w_2}\cup\overline{w_2v_0}\subset
T$ would be a subspace homeomorphic to a circle, so $T$ would not
be a tree. 
\end{proof}

\begin{defn}\label{quad1}
Given a tree $T$ we fix a base-point $v_0\in T$ and a total
ordering in all discrete sets so that the next definitions make
sense.

For any map from a discrete set $\alpha\colon A\r T$ we define the
maps
\begin{itemize}
\item $\Gamma\alpha\colon\Gamma A\r T$,

\item $\wedge^2\alpha\colon\wedge^2 A\r T$,
\end{itemize}
as follows $(a,a_1,a_2\in A; a_1\succ a_2)$
\begin{itemize}
\item $(\Gamma\alpha)(\gamma(a))=\alpha(a)$,

\item $(\Gamma\alpha)([a_1,a_2])=\ell(\alpha(a_1),\alpha(a_2))$,

\item $(\wedge^2\alpha)(a_2\wedge
a_1)=\ell(\alpha(a_2),\alpha(a_1))$.
\end{itemize}
If $\beta\colon B\r T$ is another map we also define
\begin{itemize}
\item $\alpha\otimes\beta=\ell(\alpha\times\beta)\colon A\times B\r T.$
\end{itemize}

The \emph{$T$-controlled Whitehead functor}
$$\Gamma_T\colon\ab^b(T)\To\ab^b(T)$$ is defined as the unique functor with
$$\Gamma_T\Z\grupo{A}_\alpha=\Z\grupo{\Gamma A}_{\Gamma \alpha}$$
which coincides with the ordinary Whitehead functor $\Gamma$ in
underlying abelian groups.

Similarly the \emph{$T$-controlled exterior square}
$$\wedge^2_T\colon\ab^b(T)\To\ab^b(T)$$ is the unique
functor with
$$\wedge^2_T\Z\grupo{A}_\alpha=\Z\grupo{\wedge^2A}_{\wedge^2\alpha}$$
which coincides with the ordinary exterior square $\wedge^2$ in
underlying abelian groups.

Finally the \emph{$T$-controlled tensor product} is the functor
$$-\otimes_T-\colon\ab^b(T)\times\ab^b(T)\To\ab^b(T)$$ defined as
$$\Z\grupo{A}_\alpha\otimes_T\Z\grupo{B}_\beta=\Z\grupo{A\times B}_{\alpha\otimes\beta}$$
which coincides with the ordinary tensor product in underlying
abelian groups.
\end{defn}

It is not completely trivial that the functors $\Gamma_T$,
$\wedge^2_T$ and $\otimes_T$ are well defined. One needs to check,
for example for $\wedge^2_T$, that given a controlled homomorphism
$\varphi\colon\Z\grupo{A}_\alpha\r\Z\grupo{B}_\beta$ the abelian
group homomorphism
$\wedge^2\varphi\colon\wedge^2\Z\grupo{A}\r\wedge^2\Z\grupo{B}$
induces a controlled homomorphism
$\wedge^2\varphi=\wedge^2_T\varphi\colon\wedge^2_T\Z\grupo{A}_\alpha\r\wedge^2_T\Z\grupo{B}_\beta$.

\begin{prop}\label{consistente}
The functors $\Gamma_T$, $\wedge^2_T$ and $\otimes_T$ are well
defined.
\end{prop}

In the proof of this proposition we shall use the following nice
property of the map $\ell$.

\begin{lem}\label{critterell}
Given two sequences $\set{x_n}_{n\in\N}$ and $\set{y_n}_{n\in\N}$ in
$T$, if
$\lim\limits_{n\r\infty}x_n=\lim\limits_{n\r\infty}y_n=x\in\hat{T}$
then $\lim\limits_{n\r\infty}\ell(x_n,y_n)=x$. Moreover, if
$x\in\F(T)$ the converse also holds.
\end{lem}

\begin{rem}\label{nice2}
It is well known that there is a basis of neighborhoods of the space
of Freudenthal ends of a tree $\F(T)$ in the Freudenthal
compactification $\hat{T}$ given by the sets $T_v\sqcup T_v^\F$
$(v\in T)$ where $T_v \subset T$ is the subtree formed by those
points $u\in T$ such that $v\in \overline{v_0u}$, see \cite{iht}
II.1.14, III.1.2 and III.1.3. These subtrees have the following
property.
\end{rem}

\begin{lem}\label{nice}
$x,y\in T_w$ if and only if $\ell(x,y)\in T_w$.
\end{lem}

\begin{proof}
Clearly $\Rightarrow$ holds because $\ell(x,y)\in\overline{xy}$
and $\overline{xy}\subset T_w$ since $T_w$ is a tree. On the other
hand by Definition \ref{topgen} $\ell(x,y)\in
\overline{v_0x}\cap\overline{v_0y}$ hence if $\ell(x,y)\in T_w$
then
$w\in\overline{v_0\ell(x,y)}\subset\overline{v_0x}\cap\overline{v_0y}$
so $x,y\in T_w$.
\end{proof}

\begin{proof}[Proof of Lemma \ref{critterell}]
For $x\in\F(T)$ Lemma \ref{critterell} is an immediate consequence
of Lemma \ref{nice}. If $x\in T$ then
$\lim\limits_{n\r\infty}d(x_n,y_n)=0$ and
$d(x_n,y_n)=d(x_n,\ell(x_n,y_n))+d(\ell(x_n,y_n),y_n)$ $(n\in\N)$ so
$$\lim\limits_{n\r\infty}d(x_n,\ell(x_n,y_n))=0=\lim\limits_{n\r\infty}d(y_n,\ell(x_n,y_n))$$
and therefore $\lim\limits_{n\r\infty}\ell(x_n,y_n)=x$.
\end{proof}

\begin{proof}[Proof of Proposition \ref{consistente}]
We will make the proof for $\wedge^2_T$ and leave to the reader
the cases $\Gamma_T$ and $\otimes_T$.



Consider sequences $\set{a_n^1\wedge a_n^2}_{n\geq 0}$ and
$\set{b_n^1\wedge b_n^2}_{n\geq 0}$ of elements in $\wedge^2 A$ and
$\wedge^2B$ respectively such that
$\lim\limits_{n\r\infty}\ell(\alpha(a_n^1),\alpha(a_n^2))=\varepsilon\in\F(T)$
in $\hat{T}$ and $b_n^1\wedge b_n^2$ appears with non-trivial
coefficient in the linear expansion of $\wedge^2\varphi(a_n^1\wedge
a_n^2)$. We can therefore take permutations $o_n$ of $\set{1,2}$
$(n\geq 0)$ such that $b_n^{o_n(i)}$ appears with non-trivial
coefficient in the linear expansion of $\varphi(a^i_n)$. By Lemma
\ref{critterell} $\lim\limits_{n\r\infty}\alpha(a_n^i)=\varepsilon$
$(i=1,2)$ so $\lim\limits_{n\r\infty}\beta(b_n^{o(i)})=\varepsilon$
because $\varphi$ is controlled, see Lemma \ref{critter}. This
implies that also $\lim\limits_{n\r\infty}\beta(b_n^i)=\varepsilon$
$(i=1,2)$, hence by Lemma \ref{critterell}
$\lim\limits_{n\r\infty}\ell(\beta(b_n^1),\beta(b_n^2))=\varepsilon$,
so $\wedge^2\varphi$ is controlled by Lemma \ref{critter}.
\end{proof}

The following proposition considers the compatibility of the
functors $\Gamma_T$, $\wedge^2_T$ and $\otimes_T$ with the
``change of tree'' functors determined by a proper map.

\begin{prop}\label{comp}
Given a proper map between two trees $f\colon T\r T'$ there are
natural transformations
\begin{itemize}
\item $\uf{F}^f\Gamma_T\r\Gamma_{T'}\uf{F}^f$,

\item $\uf{F}^f\wedge^2_T\r\wedge^2_{T'}\uf{F}^f$,

\item $\uf{F}^f(-\otimes_T-)\r\uf{F}^f(-)\otimes_{T'}\uf{F}^f(-)$,
\end{itemize}
which are natural equivalences provided $\F(f)\colon
\F(T)\hookrightarrow \F(T')$ is injective. Moreover, in any case
\begin{itemize}
\item $\uf{F}^f(-\otimes\Z/2)\simeq(\uf{F}^f-)\otimes\Z/2$.
\end{itemize}
\end{prop}

\begin{proof}
As in the proof of Proposition \ref{consistente} we shall leave the
cases $\Gamma_T$ and $\otimes_T$ to the reader. The case
$-\otimes\Z/2$ follows from Propositions \ref{extigrande} and \ref{acabalo} and the right-exactness of the additive
functor $\uf{F}^f_*$, since $-\otimes\Z/2$ is the cokernel of the
multiplication by $2$.

Let $\ell'\colon T'\times T'\r T'$ be the function associated to a
base-point $v_0'\in T'$ as in Definition \ref{topgen}. Notice that
for the definition of $\wedge^2_{T'}$ we choose total orderings
$\leq $ in discrete sets $A$ which are possibly different from
those orders $\preceq$ used to define $\wedge^2_T$. In order to
distinguish the two possible definitions of the set $\wedge^2A$ we
will write $\wedge^2_\leq A$ and $\wedge^2_\preceq A$. Moreover,
notice that the definition of $\wedge^2\alpha$
depends on whether the target of $\alpha$ is $T$ or $T'$.

We claim that for any possibly big free $T$-controlled $\Z$-module
$\Z\grupo{A}_\alpha$ the identity in $\wedge^2\Z\grupo{A}$ induces
a controlled homomorphism
$$\varphi\colon\uf{F}^f\wedge^2_T\Z\grupo{A}_\alpha=\Z\grupo{\wedge^2_\preceq A}_{f(\wedge^2\alpha)}\To
\Z\grupo{\wedge^2_\leq
A}_{\wedge^2(f\alpha)}=\wedge^2_{T'}\uf{F}^f\Z\grupo{A}_\alpha.$$
This morphism sends $a_1\wedge a_2\in\wedge^2_\preceq A$ to
$a_1\wedge a_2$ if $a_1<a_2$ or to $-a_2\wedge a_1$ in other case.




Suppose that $\set{a^1_n\wedge q^2_n}_{n\geq 0}$ is a sequence of
elements in $\wedge^2_\preceq A$ such that in
$\hat{T}'$ 
$$\lim\limits_{n\r\infty}
f\ell(\alpha(a^1_n),\alpha(a^2_n))=\varepsilon\in\F(T').$$ If $o_n$ is the permutation of $\set{1,2}$ such that
$a^{o(1)}_n\wedge a^{o(2)}_n$ belongs to $\wedge^2_\leq A$ then we
have to prove that
$\lim\limits_{n\r\infty}\ell'((f\alpha)(a^{o(1)}_n),(f\alpha)(a^{o(2)}_n))=\varepsilon$.
By the symmetry of $\ell'$ this is the same as showing that
$\lim\limits_{n\r\infty}\ell'((f\alpha)(a^1_n),(f\alpha)(a^2_n))=\varepsilon$.
We have that the accumulation points of
$\set{\ell(\alpha(a^1_n),\alpha(a^2_n))}_{n\geq 0}$ are contained in
$\F(f)^{-1}(\varepsilon)$, therefore by Lemma \ref{critterell} the
accumulation points of $\set{\alpha(a^1_n)}_{n\geq 0}$ and
$\set{\alpha(a^2_n)}_{n\geq 0}$ are all in
$\F(f)^{-1}(\varepsilon)$, so
$\lim\limits_{n\r\infty}(f\alpha)(a^i_n)=\varepsilon$ $(i=1,2)$.
Again by Lemma \ref{critterell}
$\lim\limits_{n\r\infty}\ell'((f\alpha)(a^1_n),(f\alpha)(a^2_n))=\varepsilon$
so $\varphi$ is controlled by Lemma \ref{critter}.

If $\F(f)$ is injective we claim that the identity in
$\wedge^2\Z\grupo{A}$ also induces a controlled homomorphism
$$\psi\colon\wedge^2_{T'}\uf{F}^f\Z\grupo{A}_\alpha=\Z\grupo{\wedge^2_\leq A}_{\wedge^2(f\alpha)}\To
\Z\grupo{\wedge^2_\preceq
A}_{f(\wedge^2\alpha)}=\uf{F}^f\wedge^2_T\Z\grupo{A}_\alpha.$$





Given a sequence $\set{a^1_n\wedge a^2_n}$ of elements in
$\wedge^2_\geq A$ such that
in $\hat{T}'$ 
$$\lim\limits_{n\r\infty}\ell'((f\alpha)(a^1_n),(f\alpha)(a^2_n))=\varepsilon\in\F(T')$$
we now have to prove that
$\lim\limits_{n\r\infty}f\ell(\alpha(a^1_n),\alpha(a^2_n))=\varepsilon$
as well. By Lemma \ref{critterell}
$\lim\limits_{n\r\infty}(f\alpha)(a^i_n)=\varepsilon$ $(i=1,2)$.
Moreover, since $\F(f)$ is injective
$\lim\limits_{n\r\infty}\alpha(a^i_n)=\varepsilon'$ is the unique
$\varepsilon'\in\F(T)$ such that $\F(\varepsilon')=\varepsilon$ so
again by Lemma \ref{critterell}
$$\lim\limits_{n\r\infty}\ell(\alpha(a^1_n),\alpha(a^2_n))=\varepsilon',$$
in particular $\lim\limits_{n\r\infty}
f\ell(\alpha(a^1_n),\alpha(a^2_n))=\F(f)(\varepsilon')=\varepsilon$
and $\psi$ is controlled by Lemma \ref{critter}.
\end{proof}

In Definition \ref{quad1} we have made several arbitrary choices
to define three functors, however up to natural equivalence the
functors do not depend on these choices, as the following
corollary shows.

\begin{cor}\label{indepe}
Up to natural equivalence the functors $\Gamma_T$, $\wedge^2_T$
and $\otimes_T$ do not depend on the choice of the base-point
$v_0\in T$ and the total orderings on discrete sets.
\end{cor}

This corollary can be regarded as a special case of Proposition
\ref{comp} for the identity map $f=1\colon T\r T$.

The functors $\Gamma_T$ and $\wedge^2_T$ are quadratic because
they are quadratic in the underlying abelian groups. Similarly
$\otimes_T$ is biadditive. We shall be interested in the composition
of these functors with the full inclusion of non-big objects
$\ab(T)\subset\ab^b(T)$ and the faithful functor
$\ab^b(T)\hookrightarrow\mo(\ab(T))$ in (\ref{yonedabig}) which we
will regard as the inclusion of a subcategory. We shall use the
same names for these composites
\begin{equation}\label{quad2}
\begin{array}{c}
\Gamma_T\colon\ab(T)\To\mo(\ab(T)),\\
\wedge^2_T\colon\ab(T)\To\mo(\ab(T)),\\
-\otimes_T-\colon\ab(T)\times\ab(T)\To\mo(\ab(T)).
\end{array}
\end{equation}

These functors satisfy the properties stated in the following
proposition where $\otimes^2_T$ denotes the square of the
bifunctor $\otimes_T$.

\begin{prop}\label{quadpropri}
The functor $\otimes_T$ is the quadratic crossed effect of
$\Gamma_T$ and $\wedge^2_T$ in (\ref{quad2}). Moreover, there are
two exact sequences of natural transformations
\begin{enumerate}
\item
$\Gamma_T\st{\tau_T}\hookrightarrow\otimes_T^2\st{q_T}\twoheadrightarrow\wedge^2_T$,

\item $\otimes^2_T\st{[-,-]_T}\To\Gamma_T\st{\sigma_T}\twoheadrightarrow -\otimes\Z/2$.
\end{enumerate}
Indeed (1) splits non-naturally; $\tau_T$ verifies
$\tau_T=r_{12}\Gamma_T(i_1+i_2)$, where $i_1$ and $i_2$ are the
inclusions of the factors of the coproduct; and
$[-,-]_T=\Gamma_T(1,1)i_{12}$. Furthermore, the functors $\Gamma_T$,
$\wedge^2_T$, $\otimes_T$ and $-\otimes\Z/2$ are compatible with the
``change of tree'' functors in the sense of Proposition \ref{comp},
as well as the natural transformations $\tau_T$, $q_T$, $[-,-]_T$
and $\sigma_T$.
\end{prop}

\begin{proof}
The last part of the statement follows from Propositions
\ref{extigrande} and \ref{comp}.

It is immediate to check that the structure maps of the quadratic
crossed effects of $\Gamma_T$ and $\wedge^2_T$ coincide with those
of the ordinary functors $\Gamma$ and $\wedge^2$ in the underlying
abelian groups. 


The exact sequence (1) is given by the upper row of (\ref{d1}) in
underlying abelian groups. In particular the formula
$\tau_T=r_{12}\Gamma_T(i_1,i_2)$ is satisfied. The splitting
$s\colon\wedge^2_T\Z\grupo{A}_\alpha\hookrightarrow\otimes^2_T\Z\grupo{A}_\alpha$
of $q_T$ and the retraction
$r\colon\otimes^2_T\Z\grupo{A}_\alpha\twoheadrightarrow\Gamma_T\Z\grupo{A}_\alpha$
of $\tau_T$ depend on a given total order $\preceq$ in $A$, more
precisely $s(a_1\wedge a_2)=a_1\otimes a_2$, $r(a_1\otimes a_2)=0$,
$r(a_2\otimes a_1)=[a_2,a_1]$ and $r(a\otimes a)=\gamma(a)$
$(a,a_1,a_2\in A,a_1\prec a_2)$.

Moreover, $[-,-]_T$ is also given by the homomorphism $[-,-]$ in
(\ref{d1}) in underlying abelian groups, in particular the equality
$[-,-]_T=\Gamma_T(1,1)i_{12}$ holds. In order to define $\sigma_T$
we consider the non-natural controlled homomorphisms
$\varsigma\colon\Gamma_T\Z\grupo{A}_\alpha\r\Z\grupo{A}_\alpha$
defined as $\varsigma(\gamma(a))=a$ and $\varsigma([a_1,a_2])=0$
$(a,a_1,a_2\in A, a_1\succ a_2)$. Now we define the $\ab(T)$-module
morphism
$\sigma_T\colon\Gamma_T\Z\grupo{A}_\alpha\r\Z\grupo{A}_\alpha\otimes\Z/2$
as
$$\sigma_T=\hat{p}\varsigma_*\colon\hom_{\ab^b(T)}(-,\Gamma_T\Z\grupo{A}_\alpha)\r\hom_{\ab^b(T)}(-,\Z\grupo{A}_\alpha)\otimes\Z/2,$$
where $\hat{p}$ is the natural projection $\hat{p}\colon
1\twoheadrightarrow -\otimes\Z/2$ in the category of abelian groups.
This is an epimorphism since $\varsigma$ is an epimorphism which
admits a section
$j\colon\Z\grupo{A}_\alpha\hookrightarrow\Gamma_T\Z\grupo{A}_\alpha$
defined as $j(a)=\gamma(a)$ $(a\in A)$. The naturality of $\sigma_T$
can be easily checked by using the naturality of $\sigma$ in
(\ref{d1}), therefore it is only left to check the exactness of (2)
in the middle.

Given a controlled homomorphism
$\varphi\colon\Z\grupo{B}_\beta\r\Gamma_T\Z\grupo{A}_\alpha$,
$\sigma_T\varphi=0$ if and only if
$(\varsigma\varphi)\otimes\Z/2=0$. This means that the coefficient
of $\gamma(a)$ in the linear expansion of $\varphi(b)$ is even for
all $a\in A$ and $b\in B$, so we can define a controlled
homomorphism
$\psi\colon\Z\grupo{B}_\beta\r\Gamma_T\Z\grupo{A}_\alpha$ in such a
way that the coefficient of $\psi(b)$ in $\gamma(a)$ is half the
corresponding coefficient in $\varphi(b)$ and the coefficients of
$[a_1,a_2]$ in $\psi(b)$ and $\varphi(b)$ coincide for all $a_1\succ
a_2$. In fact by construction $\varphi$ is controlled if and only if
$\psi$ is. The controlled homomorphism
$k\colon\Gamma_T\Z\grupo{A}_\alpha\rightarrow\otimes^2_T\Z\grupo{A}_\alpha$
defined by $k(\gamma(a))=a\otimes a$ and $k([a_1,a_2])=a_1\otimes
a_2$ $(a,a_1,a_2\in A; a_1\succ a_2)$ satisfies
$([-,-]_Tk)(\gamma(a))=2\gamma(a)$ and
$([-,-]_Tk)([a_1,a_2])=[a_1,a_2]$, so $\varphi=[-,-]_Tk\psi$ as
required.
\end{proof}

\begin{defn}\label{qre}
A quadratic functor $F\colon\C{A}\r\C{B}$ between abelian categories
is \emph{right-exact} if for any exact sequence in $\C{A}$
$$X\st{f}\r Y\st{g}\twoheadrightarrow Z$$ the following sequence in
$\C{B}$ is also exact
$$F(X)\oplus F(X|Y)\st{\zeta}\r F(Y)\st{F(g)}\twoheadrightarrow F(Z),\;\; \zeta=(F(f), F(f,1)i_{12}).$$
If $F$ is right-exact its quadratic crossed effect $F(-|-)$ is
exact in each variable in the usual additive sense.
\end{defn}

\begin{prop}\label{extender}
If $F\colon\C{A}\r\C{B}$ is a quadratic functor from a small
additive category $\C{A}$ to an abelian category $\C{B}$ with exact
filtered colimits, up to natural equivalence there exists a unique
extension $F\colon\mo(\C{A})\r\C{B}$ of $F$ through the Yoneda
inclusion $\C{A}\hookrightarrow\mo(\C{A})$ in (\ref{yoneda}) which
is right-exact, quadratic, and preserves filtered colimits.
Moreover, if $G\colon\C{A}\r\C{B}$ is another functor and $\mu\colon
F\r G$ is a natural transformation $\mu$ extends uniquely to
$\mo(\C{A})$.
\end{prop}

\begin{rem}\label{sketch}
It is not difficult to obtain a direct proof of this proposition,
however it is a tedious work to carry out all necessary
verifications and details are not relevant for the rest of this
paper, for this reason we leave them to the interested reader. The
idea is to extend the functor $F\colon\C{A}\r\C{B}$ to f. p.
$\C{A}$-modules by right-exactness and to all $\C{A}$-modules by
using the fact that any $\C{A}$-module is a filtered colimit of f.
p. ones. The uniqueness up to natural equivalence follows
immediately. The extended functor $F\colon\mo(\C{A})\r\C{B}$ is
quadratic since it is quadratic over f. g. free ones and f. g. free
$\C{A}$-modules are small projective generators of $\mo(\C{A})$. The
extension of the natural transformation $\mu$ can be constructed at
the same time as the extensions of the functors $F$ and $G$ by using
the same arguments.

Alternatively one can check that the extension of $F$ is the Andr\'e
homology functor $H_0(-,F)\colon\mo(\C{A})\r\C{B}$ in the sense of
\cite{msahac}. In fact a quadratic functor
$M\colon\mo(\C{A})\r\C{B}$ is right-exact in the sense of Definition
\ref{qre} if and only if its derived functor $L_0M(-,0)$ in the
sense of Dold and Puppe \cite{hnaf} coincides with $M$. On the other
hand by \cite{msahac} 1.1 and 3.3
$L_0H_0(-,F)(-,0)=H_0(-,H_0(-,F)_{|_{\C{A}}})=H_0(-,F)$ therefore
$H_0(-,F)$ is right-exact. The functor $H_0(-,F)$ preserves filtered
colimits by \cite{caccha} 14 (c). As an immediate consequence of
this characterization of extensions in terms of Andr\'e homology we
observe that whenever we have an exact sequence of natural
transformations
$$F\r G\twoheadrightarrow N$$
between functors $\C{A}\r\C{B}$ the extended sequence of functors
$\mo(\C{A})\r\C{B}$ is exact in the same way.
\end{rem}

We use Proposition \ref{extender} to extend the quadratic functors
$\Gamma_T$ and $\wedge^2_T$ in (\ref{quad2}) to the whole category
of $\ab(T)$-modules.
\begin{equation}\label{quad3}
\wedge^2_T,\Gamma_T\colon\mo(\ab(T))\To\mo(\ab(T)).
\end{equation}
Now we can define the $T$-controlled tensor product of
$\ab(T)$-modules
\begin{equation}\label{quad4}
-\otimes_T-\colon\mo(\ab(T))\times\mo(\ab(T))\To\mo(\ab(T)).
\end{equation}
as the quadratic crossed effect of either $\Gamma_T$ or
$\wedge^2_T$. These two quadratic functors have the same quadratic
crossed effect since their crossed effects are also right-exact and
preserve filtered colimits in each variable, and both coincide with
the $T$-controlled tensor product in (\ref{quad2}) over f. g. free
$\ab(T)$-modules by Proposition \ref{quadpropri}. In particular
$\otimes_T$ in (\ref{quad4}) is an extension of $\otimes_T$ in
(\ref{quad2}). 

\begin{prop}\label{quadpropri2}
There are two exact sequences of natural transformations involving
the functors in (\ref{quad3}) and (\ref{quad4})
\begin{equation*}
\begin{array}{c}
\Gamma_T\st{\tau_T}\r\otimes_T^2\st{q_T}\twoheadrightarrow\wedge^2_T,\\
\otimes^2_T\st{[-,-]_T}\To\Gamma_T\st{\sigma_T}\twoheadrightarrow
-\otimes\Z/2,
\end{array}
\end{equation*}
such that $\tau_T=r_{12}\Gamma_T(i_1+i_2)$ and
$[-,-]_T=\Gamma_T(1,1)i_{12}$. Moreover, the functors $\Gamma_T$,
$\wedge^2_T$, $\otimes_T$ and $-\otimes\Z/2$ are compatible with the
``change of tree'' functors in the sense of Proposition \ref{comp},
as well as the natural transformations $\tau_T$, $q_T$, $[-,-]_T$
and $\sigma_T$.
\end{prop}

This proposition follows from Propositions \ref{quadpropri} and
\ref{extender} and 
Remark \ref{sketch}. 

\begin{defn}\label{quad5}
The \emph{$T$-controlled reduced tensor square}
$$\hat{\otimes}^2_T\colon\mo(\ab(T))\To\mo(\ab(T))$$ is defined by
the following natural push-out
$$\xymatrix{\Gamma_T\ar[d]_{\sigma_T}\ar[r]^{\tau_T}\ar@{}[dr]|{\mathrm{push}}&\otimes^2_T\ar[d]^{\bar{\sigma}_T}\\-\otimes\Z/2\ar[r]^<(.35){\bar{\tau}_T}&\hat{\otimes}^2_T}$$
In particular $\hat{\otimes}^2_T$ is quadratic, right-exact and
preserves filtered colimits.
\end{defn}

This definition completes the extension of (\ref{d1}) to the
category of $\ab(T)$-modules, that is, we have a commutative
diagram of natural transformations with exact rows and column
\begin{equation}\label{d2}
\xymatrix{\otimes^2_T\ar[d]^{[-,-]_T}&&\\\Gamma_T
\ar[r]^{\tau_T}\ar@{->>}[d]^{\sigma_T}\ar@{}[dr]|{\mathrm{push}}&\otimes^2_T\ar@{->>}[d]^{\bar{\sigma}_T}\ar@{->>}[r]^{q_T}&\wedge^2_T\ar@{=}[d]\\
-\otimes\Z_2\ar[r]^{\bar{\tau}_T}&\hat{\otimes}^2_T\ar@{->>}[r]^{\bar{q}_T}&\wedge^2_T}
\end{equation}
Here $\tau_T$ and $\sigma_T$, and therefore the whole diagram, are
determined by the functor $\Gamma_T$, see Proposition
\ref{quadpropri2}, as it happened in the ordinary case (\ref{d1}).


\begin{prop}\label{conprts}
Given a proper map between two trees $f\colon T\r T'$ there is a
natural transformation
$$\uf{F}^f\hat{\otimes}^2_T\To\hat{\otimes}^2_{T'}\uf{F}^f$$
which is a natural equivalence provided $\F(f)\colon\F(T)\r\F(T')$
is injective. The natural transformations $\bar{\tau}_T$,
$\bar{\sigma}_T$ and $\bar{q}_T$ are also compatible with ``change
of tree'' functors. Moreover, $\otimes_T$ is the quadratic crossed
effect of $\hat{\otimes}^2_T$.
\end{prop}

\begin{proof}
The part concerning ``change of tree'' functors follows from
Definition \ref{quad5} and Proposition \ref{quadpropri2}. The
natural transformation $\bar{q}_T$ in (\ref{d2}) induces an
equivalence between the crossed effects because $-\otimes\Z/2$ is
additive.
\end{proof}


The following proposition shows that $\otimes_T$, $\wedge^2_T$ and
$\hat{\otimes}^2_T$ restrict to the full subcategory of
$\C{vect}(T)$-modules.



\begin{prop}\label{conZ2}
There are natural equivalences
\begin{enumerate}
\item
$(\Gamma_T-)\otimes\Z/2\simeq(\Gamma_T(-\otimes\Z_2))\otimes\Z/2$,

\item $(-\otimes_T-)\otimes\Z/2\simeq(-)\otimes_T(-\otimes\Z/2)$,

\item $(\wedge^2_T-)\otimes\Z/2\simeq\wedge^2_T(-\otimes\Z/2)$,

\item
$(\hat{\otimes}^2_T-)\otimes\Z/2\simeq\hat{\otimes}^2_T(-\otimes\Z/2)$.
\end{enumerate}
In particular the functors $(\Gamma_T-)\otimes\Z/2$, $\otimes_T$,
$\wedge^2_T$ and $\hat{\otimes}^2_T$ restrict to the category of
$\C{vect}(T)$-modules.
\end{prop}


\begin{proof}
All these isomorphisms are induced by the natural projection
$\hat{p}\colon 1\twoheadrightarrow-\otimes\Z/2$. A free
$T$-controlled $\uf{F}_2$-module
$\uf{F}_2\grupo{A}_\alpha=\Z\grupo{A}_\alpha\otimes\Z/2$ is the
cokernel of the multiplication-by-$2$ controlled homomorphism
$2\colon\Z\grupo{A}_\alpha\r\Z\grupo{A}_\alpha$. Since $\Gamma_T$ is
right-exact and coincides on $\ab(T)$ with the ordinary Whitehead
functor $\Gamma$ in underlying abelian groups we obtain an exact
sequence
$$\Gamma_T\Z\grupo{A}_\alpha\oplus\Z\grupo{A}_\alpha\otimes_T\Z\grupo{A}_\alpha\st{(4,2[-,-]_T)}\To
\Gamma_T\Z\grupo{A}_\alpha\st{\Gamma_T\hat{p}}\twoheadrightarrow\Gamma_T\uf{F}_2\grupo{A}_\alpha,$$
therefore we clearly obtain the desired isomorphism by applying
$-\otimes\Z/2$
$$(\Gamma_T\hat{p})\otimes\Z/2\colon(\Gamma_T\Z\grupo{A}_\alpha)\otimes\Z/2\st{\simeq}\To(\Gamma_T\uf{F}_2
\grupo{A}_\alpha)\otimes\Z/2.$$ The isomorphism (1) for f. p.
$\ab(T)$-modules follows now from right-exactness, and for arbitrary
$\ab(T)$-modules from the fact that they are filtered colimits of f.
p. ones.

The isomorphism (2) follows easily from the biadditivity and right
exactness properties of $\otimes_T$. The rest of isomorphisms follow
from (1), (2) and the properties of (\ref{d2}).
\end{proof}

The constructions of $\otimes_T$, $\wedge^2_T$ and $\Gamma_T$ allow
us to compute these functors from projective resolutions in a very
convenient way. We do not have in general such a method for
$\hat{\otimes}^2_T$, however as a consequence of the following
proposition we can obtain an analogous method when we restrict to
$\C{vect}(T)$-modules. In its statement we use the following
notation. Given a total ordering $\preceq$ on a discrete set $A$ we
define the set
$$\hat{\otimes}^2A=\set{a_1\hat{\otimes}a_2\,;\,a_1\preceq a_2},$$
and given a function $\alpha\colon A\r T$ we define another one
$$\hat{\otimes}^2\alpha\colon\hat{\otimes}^2A\r T,\;\;
(\hat{\otimes}^2\alpha)(a_1\hat{\otimes}a_2)=\ell(\alpha(a_1),\alpha(a_2)).$$

\begin{prop}\label{redten}
The restriction of $\hat{\otimes}^2_T$ to $\C{vect}(T)$ factors
through the functor
$\hat{\otimes}^2_T\colon\C{vect}(T)\r\C{vect}^b(T)$ which is
$\hat{\otimes}^2_T\uf{F}_2\grupo{A}_\alpha=\uf{F}_2\grupo{\hat{\otimes}^2A}_{\hat{\otimes}^2\alpha}$
on objects and coincides with the ordinary reduced tensor square on
underlying $\uf{F}_2$-vector spaces. Moreover, the natural
transformations $\bar{\tau}_T$, $\bar{\sigma}_T$, $\bar{q}_T$ and
the structure morphisms of the quadratic crossed effect $\otimes_T$
over $\C{vect}(T)$ also coincides with those of the ordinary reduced
tensor square on underlying $\uf{F}_2$-vector spaces.
\end{prop}

\begin{proof}
Since $-\otimes\Z/2$ is right-exact and
$\uf{F}_2\grupo{A}_\alpha=\Z\grupo{A}_\alpha\otimes\Z/2$, by
Proposition \ref{conZ2} $\hat{\otimes}^2_T\uf{F}_2\grupo{A}_\alpha$
is also the following push-out
$$\xymatrix{(\Gamma_T\Z\grupo{A}_\alpha)\otimes\Z/2\ar[d]_{\sigma_T\otimes\Z/2}\ar[r]^{\tau_T\otimes\Z/2}\ar@{}[dr]|{\mathrm{push}}&
(\otimes^2_T\Z\grupo{A}_\alpha)\otimes\Z/2\ar[d]^{\bar{\sigma}_T}\\
\Z\grupo{A}_\alpha\otimes\Z/2\ar[r]^{\bar{\tau}_T}&\hat{\otimes}^2_T\uf{F}_2\grupo{A}_\alpha}$$
By Propositions \ref{calculete} and \ref{quadpropri} we have that $(a,a_1,a_2\in A;a_1\succ a_2)$
\begin{eqnarray*}
(\Gamma_T\Z\grupo{A}_\alpha)\otimes\Z/2&=&\uf{F}_2\grupo{\Gamma
A}_{\Gamma\alpha},\\
(\otimes^2_T\Z\grupo{A}_\alpha)\otimes\Z/2&=&\uf{F}_2\grupo{A\times
A}_{\alpha\otimes\alpha},\\
\Z\grupo{A}_\alpha\otimes\Z/2&=&\uf{F}_2\grupo{A}_\alpha,\\
(\tau_T\otimes\Z/2)(\gamma(a))&=&a\otimes a,\\
(\tau_T\otimes\Z/2)([a_1,a_2])&=&a_1\otimes a_2+a_2\otimes a_1,\\
(\sigma_T\otimes\Z/2)(\gamma(a))&=&a,\\ 
(\sigma_T\otimes\Z/2)(
[a_1,a_2])&=&0.
\end{eqnarray*}  
Now it is easy to see
that the homomorphisms
$\bar{\tau}_T\colon\uf{F}_2\grupo{A}_\alpha\r\uf{F}_2\grupo{\hat{\otimes}^2A}_{\hat{\otimes}^2\alpha}$,
$\bar{\tau}_T(a)=a\hat{\otimes}a$ $(a\in A)$ and
$\bar{\sigma}_T\colon\uf{F}_2\grupo{A\times
A}_{\alpha\otimes\alpha}\r\uf{F}_2\grupo{\hat{\otimes}^2A}_{\hat{\otimes}^2\alpha}$
$\bar{\sigma}_T(a_1\otimes a_2)=a_1\hat{\otimes} a_2$ $(a_1,a_2\in
A)$ are controlled homomorphisms and
$$\xymatrix{(\Gamma_T\Z\grupo{A}_\alpha)\otimes\Z/2\ar[d]_{\sigma_T\otimes\Z/2}\ar[r]^{\tau_T\otimes\Z/2}\ar@{}[dr]|{\mathrm{push}}&
(\otimes^2_T\Z\grupo{A}_\alpha)\otimes\Z/2\ar[d]^{\bar{\sigma}_T}\\
\Z\grupo{A}_\alpha\otimes\Z/2\ar[r]^{\bar{\tau}_T}&\uf{F}_2\grupo{\hat{\otimes}^2A}_{\hat{\otimes}^2\alpha}}$$
is a push-out.


The functors $\Gamma_T$, $\wedge^2_T$, $\otimes^2_T$ as well as the
structure morphisms of their crossed effects and the natural
transformations $[-,-]_T$, $\tau_T$, $q_T$ are known to agree with
their ordinary analogues on underlying abelian groups when evaluated
in free $T$-controlled $\Z$-modules, see Proposition
\ref{quadpropri}. Moreover, we observe that in this case
$\sigma_T\otimes\Z/2$ also coincides with $\sigma\otimes\Z/2$ on
underlying $\uf{F}_2$-vector spaces. Furthermore, we also notice
above that $\bar{\tau}_T$ and $\bar{\sigma}_T$ coincide with
$\bar{\tau}$ and $\bar{\sigma}$ on underlying $\uf{F}_2$-vector
spaces, therefore the properties of diagram (\ref{d2}) also imply
this fact for $\bar{q}_T$, and hence for the structure morphisms of
the quadratic crossed effect of $\hat{\otimes}_T$ because
$\bar{q}_T$ induces an isomorphism between the quadratic crossed
effects of $\hat{\otimes}^2_T$ and $\wedge^2_T$, see Proposition
\ref{conprts}.
\end{proof}

\section{The lower Whitehead module}\label{slwm}

The following theorem shows that the lower Whitehead module of a $T$-complex can be obtained from the homology by using some of the functors defined in the previous section.

\begin{thm}\label{lwm}
For any $(n-1)$-connected $T$-complex $X$ there is a natural
isomorphism $$\Gam_{n+1}X\simeq\left\{%
\begin{array}{ll}
    \Gamma_T\hh_2X, & \hbox{if $n=2$;} \\ & \\
    \hh_nX\otimes\Z_2, & \hbox{if $n\geq 3$.} \\
\end{array}%
\right.    $$ These isomorphisms are compatible with the
suspension morphisms in proper homology and Whitehead modules,
that is, for $n=2$ there is a commutative diagram
$$\xymatrix{\Gamma_T\hh_2X\ar@{->>}[r]^<(.2){\sigma_T}&
\hh_2X\otimes\Z_2\ar[r]_<(.2)\simeq^<(.19){\Sigma\otimes\Z_2}&\hh_3\Sigma X\otimes\Z_2\\
\Gam_3X\ar[rr]^\Sigma\ar[u]^\simeq&&\Gam_4\Sigma X\ar[u]_\simeq}$$
and for $n\geq 3$
$$\xymatrix{\hh_nX\otimes\Z_2\ar[r]_<(.17)\simeq^<(.18){\Sigma\otimes\Z_2}&\hh_{n+1}\Sigma
X\otimes\Z_2\\\Gam_{n+1}X\ar[r]^\Sigma\ar[u]^\simeq&\Gam_{n+2}\Sigma
X\ar[u]_\simeq }$$
\end{thm}

\begin{proof}
Let us check that this proposition holds in the homotopy category
category of $n$-dimensional spherical objects $\C{S}^n(T)$, which
are obviously $(n-1)$-connected. Notice that the functor
$\Gamma_{n+1}$ coincides with $\Pi_{n+1}$ over this category.

By \cite{iht} II.4.21 the isomorphism in the stable case $n\geq 3$
holds.

In \cite{iht} II.2.15 the sets of proper homotopy classes
$[S^3_\beta,S^2_\alpha]^T$ are computed as subsets of the sets of
ordinary homotopy classes of pointed maps
$[S^3_\alpha/T,S^2_\beta/T]^*_{\mathrm{ord}}$ by using the
concepts of carrier and $(T,\alpha,\beta)$-proper element in
\cite{iht} II.2.13 and II.2.14, respectively, however the
equivalent description in terms of
$(\bar{T},\alpha,\beta)$-controlled elements given \cite{iht}
III.2.7 is better for our purposes, see \cite{iht} III.2.8. It is
well-known that $\pi_3(S^2_\beta/T)=\Gamma\Z\grupo{B}$ where $B$
is the source of $\beta$. Given a total ordering $\preceq$ in $B$
the carrier of an element in the basis $\Gamma B$ of
$\Gamma\Z\grupo{B}=[S^3,S^2_\alpha/T]_{\mathrm{ord}}^*$ is either
$\mathrm{carr}(\gamma(b))=\set{\beta(b)}$ or
$\mathrm{carr}([b_1,b_2])=\set{\beta(b_1),\beta(b_2)}$ for
$b,b_1,b_2\in B$ with $b_1\succ b_2$. Now it is easy to see by using
Lemma \ref{nice} that an element $\varphi\in
[S^3_\alpha/T,S^2_\beta/T]^*_{\mathrm{ord}}=\hom(\Z\grupo{A},\Gamma\Z\grupo{B})$,
where $A$ is the source of $\alpha$, is
$(\bar{T},\alpha,\beta)$-controlled if and only if
$\varphi\colon\Z\grupo{A}_\alpha\r\Gamma_T\Z\grupo{B}_\beta$ is a
controlled homomorphism, so the isomorphism also holds in the
unstable case $n=2$.

It is known that the suspension
$\Sigma\colon[S^3_\alpha/T,S^2_\beta/T]^*_{\mathrm{ord}}\r[S^4_\alpha/T,S^3_\beta/T]^*_{\mathrm{ord}}$
is equivalent to the following homomorphism
$$\hom(\Z\grupo{A},\Gamma\Z\grupo{B})\st{\hom(1,\sigma)}\To\hom(\Z\grupo{A},\Z\grupo{B}\otimes\Z/2)=
\hom(\Z\grupo{A},\Z\grupo{B})\otimes\Z/2.$$ Notice also that
$\hom(1,\sigma)=\hat{p}\hom(1,\varsigma)$, where, $\varsigma$ is the
homomorphism used to define $\sigma_T$ in the proof of Proposition
\ref{quadpropri} and $\hat{p}$ is the natural projection
$\hat{p}\colon 1\twoheadrightarrow-\otimes\Z/2$, hence the first
diagram of the statement commutes for spherical objects.

Once we have proved that the statement is true for spherical
objects the general case can be checked as in \cite{ah} IX.4.5
since a proper homotopy version of \cite{ah} V.7.6 holds as a
consequence of the proper Blakers-Massey theorem. 
\end{proof}

We shall often use the isomorphisms in Theorem \ref{lwm} as
identifications. Below we give some propositions whose proofs follow by diagram chasing from Theorem \ref{lwm},
the properties of the quadratic functors defined in Section
\ref{eqf} and the elementary proper homotopy theory reviewed in
Sections \ref{icat} and \ref{ephai}.

\begin{prop}\label{piaditiv}
Given two $1$-connected $T$-complexes $X$ and $Y$ the following
sequence is splitting short exact
$$\xymatrix{\h_2X\otimes_T\h_2Y\ar@{^{(}->}[r]^<(.22){i_3i_{12}}&\Pi_3(X\vee
Y)\ar@{->>}[r]^<(.2)\xi&\Pi_3X\oplus\Pi_3Y\ar@/^12pt/@{-->}[l]^{\zeta}}.$$ Here $\xi=\left(%
\begin{array}{c}
  {p_1}_* \\
  {p_2}_*
\end{array}%
\right)$ and $\zeta=({i_1}_*,{i_2}_*)$, where $i_k$ and $p_k$
$(k=1,2)$ are the inclusions and retractions of the first and
second factor of the coproduct $X\vee Y$; $i_3$ is one of the
morphisms in Whitehead's long exact sequence, and $i_{12}$ is the
inclusion of the quadratic crossed effect of $\Gamma_T$.
\end{prop}

\begin{proof}
One only has to consider the following commutative diagram with
exact row and columns
$$\xymatrix{&\Pi_4(X\vee Y)\ar@{->>}[r]^\xi\ar[d]_{h_4}&\Pi_4X\oplus \Pi_4Y\ar[d]^{h_4\oplus h_4}\ar@/^12pt/@{-->}[l]^\zeta
            \\&\h_4(X\vee Y)\ar@{=}[r]\ar[d]_{b_4}&\h_4 X\oplus\h_4
            Y\ar[d]^{b_4\oplus b_4}
            \\\h_2X\otimes_T\h_2Y\ar@{^{(}->}[r]^<(.2){i_{12}}&\Gamma_T(\h_2X\oplus\h_2Y)\ar@{->>}[r]\ar[d]_{i_3}&
            \Gamma_T\h_2X\oplus\Gamma_T\h_2Y\ar[d]^{i_3\oplus i_3}
            \\&\Pi_3(X\vee Y)\ar@{->>}[r]^{\xi}\ar@{->>}[d]_{h_3}&
            \Pi_3X\oplus\Pi_3Y\ar@{->>}[d]^{h_3\oplus h_3}\ar@/^12pt/@{-->}[l]^\zeta
            \\&\h_3(X\vee Y)\ar@{=}[r]&\h_3 X\oplus\h_3 Y}$$
\end{proof}

\begin{prop}\label{nucleosusp}
Given a $1$-connected $T$-complex $X$ the following sequence is
exact
$$\otimes^2_T\h_2X\st{i_3[-,-]_T}\To\Pi_3X\st{\Sigma}\twoheadrightarrow\Pi_4\Sigma X.$$
\end{prop}

\begin{proof}
It is enough to follow the next commutative diagram with exact rows
and column where we use the suspension isomorphism in proper
homology $\h_2X\simeq\h_3\Sigma X$ as an identification
$$\xymatrix{&\otimes^2_T\h_2X\ar[d]_{[-,-]_T}&&
            \\\h_4X\ar[d]_{\Sigma}^\simeq\ar[r]^<(.3){b_4}&\Gamma_T\h_2X\ar@{->>}[d]_{\sigma_T}\ar[r]^<(.35){i_3}&
              \Pi_3X\ar[d]_{\Sigma}\ar@{->>}[r]^<(.3){h_3}&\h_3X\ar[d]_{\Sigma}^\simeq
            \\\h_5\Sigma X\ar[r]^<(.25){b_5}&\h_2X\otimes\Z_2\ar[r]^<(.25){i_4}&\Pi_4\Sigma X\ar@{->>}[r]^<(.25){h_4}&\h_4\Sigma X}$$
\end{proof}

\begin{prop}\label{yanose}
Given two $T$-complexes $X^2$ and $Z$ with $\dim X^2\leq 2$ there
are natural central extensions $(n\geq 1)$
$$H^{n+2}(\Sigma^nX^2,\Pi_{n+2}Z)\st{j}\hookrightarrow[\Sigma^nX^2,Z]^T\twoheadrightarrow H^{n+1}(\Sigma^nX^2,\Pi_{n+1}Z).$$
\end{prop}

\begin{proof}
The long cofiber sequence associated to the pasting map of $2$-cells
$f\colon S^1_{\alpha_2}\rightarrow S^1_{\alpha_1}$ in $X^2$ gives
rise to exact sequences (in the left column) which fit into the
following commutative diagrams $(n\geq 1)$
$$\xymatrix{[S^{n+2}_{\alpha_1},Z]^T\ar[d]_{(\Sigma^{n+1}f)^*}\ar@{=}[r]&(\Pi_{n+2}Z)(\mathcal{C}_{n+1}\Sigma^nX^2)\ar[d]^{d_{n+2}^*}
\\[S^{n+2}_{\alpha_2},Z]^T\ar[d]_{(\Sigma^nq)^*}\ar@{=}[r]&(\Pi_{n+2}Z)(\mathcal{C}_{n+2}\Sigma^nX^2)\\[\Sigma^nX^2,Z]^T\ar[d]&\\
[S^{n+1}_{\alpha_1},Z]^T\ar[d]_{(\Sigma^nf)^*}\ar@{=}[r]&(\Pi_{n+1}Z)(\mathcal{C}_{n+1}\Sigma^nX^2)\ar[d]^{d_{n+2}^*}\\[S^{n+1}_{\alpha_2},Z]^T\ar@{=}[r]&
(\Pi_{n+1}Z)(\mathcal{C}_{n+2}\Sigma^nX^2)}$$ Here $(\Sigma^nq)^*$
is always central, see \cite{ah} II.8.26.
\end{proof}

\begin{prop}\label{portutatis}
Given three $T$-complexes $X^2$, $Y$ and $Z$ with $\dim X^2\leq 2$
and $Y$ and $Z$ $1$-connected, we have a natural central extension
$$\xymatrix@C=32pt{H^3(\Sigma X^2,\h_2Y\otimes_T\h_2Z)\ar@{^{(}->}[r]^<(.25){j(i_3i_{12})_*}&
[\Sigma X^2,Y\vee
Z]^T\ar@{->>}[r]^<(.20){({p_1}_*,{p_2}_*)}&[\Sigma X^2,Y]^T\times
[\Sigma X^2,Z]^T}.$$
\end{prop}

\begin{proof}
By evaluating the central extension in Proposition \ref{yanose}
for $n=1$ on the retractions of the coproduct $p_1\colon Y\vee
Z\rightarrow Y$, $p_2\colon Y\vee Z\rightarrow Z$ and using
Hurewicz's isomorphism as an identification we obtain the
following commutative diagram whose columns are central extensions
$$\xymatrix{H^3(\Sigma X^2,\Pi_3(Y\vee
Z))\ar[r]\ar@{^{(}->}[d]^j&H^3(\Sigma X^2,\Pi_3Y)\oplus H^3(\Sigma
X^2,\Pi_3Y)\ar@{^{(}->}[d]^j\\[\Sigma X^2,Y\vee Z]^T\ar[r]^<(.25){({p_1}_*,{p_2}_*)}\ar@{->>}[d]&[\Sigma X^2,Y]^T\times [\Sigma X^2,Z]^T\ar@{->>}[d]\\H^2(\Sigma
X^2,\h_2(Y\vee Z))\ar@{=}[r]&H^2(\Sigma X^2,\h_2Y)\oplus
H^2(\Sigma X^2,\h_2Z)}$$ Moreover, 
if we evaluate this functor in the splitting short exact sequence in
Proposition \ref{piaditiv} we obtain another one
$$H^3(\Sigma X^2,\h_2Y\otimes_T\h_2Z)\st{(i_3i_{12})_*}\hookrightarrow  H^3(\Sigma
X,\Pi_3(Y\vee Z))\twoheadrightarrow H^3(\Sigma X,\Pi_3Y\oplus
\Pi_3Z).$$ Now the proposition follows from the diagram obtained by
fitting this short exact sequence in the previous diagram.
\end{proof}

\begin{prop}\label{porcelanosa}
Given two $T$-complexes $X^2$ and $Z$ with $\dim X^2\leq 2$ and
$Z$ $1$-connected, the following natural sequence is exact
$$H^3(\Sigma X^2,\otimes^2_T\h_2Z)\st{j(i_3[-,-]_T)_*}\To[\Sigma
X^2,Z]^T\st{\Sigma}\twoheadrightarrow[\Sigma^2 X^2,\Sigma Z]^T.$$
Moreover, the suspension operator induces isomorphisms $(n\geq 2)$
$$\Sigma\colon[\Sigma^nX^2,\Sigma^{n-1}Z]^T\simeq[\Sigma^{n+1}X^2,\Sigma^nZ]^T.$$
\end{prop}

\begin{proof}
The central extensions in Proposition \ref{yanose}, the suspension
homomorphisms in proper homotopy modules and (co)homology, and the
Hurewicz isomorphism give rise to the following commutative diagram
whose columns are central extensions
$$\xymatrix{H^{n+2}(\Sigma^n
X^2,\Pi_{n+2}\Sigma^{n-1}Z)\ar@{^{(}->}[d]^j\ar[r]&H^{n+3}(\Sigma^{n+1}
X^2,\Pi_{n+3}\Sigma^nZ)\ar@{^{(}->}[d]^j\\[\Sigma^nX^2,\Sigma^{n-1}Z]^T\ar@{->>}[d]\ar[r]&
[\Sigma^{n+1}X^2,\Sigma^nZ]^T\ar@{->>}[d]
\\H^{n+1}(\Sigma^nX^2,\h_{n+1}\Sigma^{n-1}Z)\ar[r]^\simeq&H^{n+2}(\Sigma^{n+1}X^2,\h_{n+2}\Sigma^nZ)}$$
Moreover, by Freudenthal's suspension theorem in proper homotopy
theory the homomorphism
$\Sigma\colon\Pi_{n+2}\Sigma^{n-1}Z\rightarrow\Pi_{n+3}\Sigma^nZ$ is
an isomorphism for $n\geq 2$, hence the upper row in the diagram is
an isomorphism within this range, and the second part of the
statement follows. On the other hand for $n=1$ the functor
$H^3(\Sigma X^2,-)$ is right-exact since $\dim\Sigma X^2\leq 3$, in
particular by applying it to the exact sequence in Proposition
\ref{nucleosusp} we obtain another exact sequence
$$H^3(\Sigma X^2,\otimes^2_T\h_2Z)\st{(i_3[-,-]_T)_*}\To H^3(\Sigma
X^2,\Pi_3Z)\twoheadrightarrow H^3(\Sigma X^2,\Pi_4\Sigma Z).$$ The
first part of the proposition follows from the diagram obtained by
inserting this exact sequence into the previous diagram for $n=1$.
\end{proof}

\section{Free controlled groups of nilpotency class $2$}\label{f2}

Recall that a free group with basis $A$ in the variety of groups
of nilpotency class $2$ is the quotient
$$\grupo{A}^\ni=\frac{\grupo{A}}{\Gamma_3\grupo{A}},$$ where
$\Gamma_3\grupo{A}\subset\grupo{A}$ is the third term of the lower
central series of the free group $\grupo{A}$, that is the subgroup
generated by commutators of weight three, i. e. $[x,[y,z]]$
$(x,y,z\in \grupo{A})$. Here $[x,y]=-x-y+x+y$.

\begin{defn}
A \emph{free $T$-controlled group of nilpotency class $2$} is a pair
$\grupo{A}_\alpha^\ni$ formed by a free group of nilpotency class
$2$ $\grupo{A}^\ni$ with basis the discrete set $A$ and a proper map
$\alpha\colon A\r T$. A \emph{controlled homomorphism}
$\varphi\colon\grupo{A}_\alpha^\ni\r\grupo{B}_\beta^\ni$ is a
homomorphism between the underlying groups such that for any
neighbourhood $U$ of $\varepsilon\in\F(T)$ in $\hat{T}$ there exists
another one $V\subset U$ such that
$\varphi(\alpha^{-1}(V))\subset\grupo{\beta^{-1}(U)}^\ni$.
\end{defn}

Notice that this definition is completely analogous to that of free
$T$-controlled $R$-module in Section \ref{freecontrol} and free
$T$-controlled group in Definition \ref{fcnag}. We will denote
$\nil(T)$ to the category of free $T$-controlled groups of
nilpotency class $2$ and controlled homomorphisms. This category
admitis finite coproducts. More precisely, the coproduct of two
objects is
$$\grupo{A}^\ni_\alpha\vee\grupo{B}^\ni_\beta=\grupo{A\sqcup B}^\ni_{(\alpha,\beta)}.$$

There are full abelianization and nilization functors
$$\abb\colon\nil(T)\To\ab(T)\colon\grupo{A}^\ni_\alpha\mapsto\Z\grupo{A}_\alpha,$$
$$\ni\colon\gr(T)\To\nil(T)\colon\grupo{A}_\alpha\mapsto\grupo{A}_\alpha^\ni.$$

Recall that there exists a natural central extension
$$\wedge^2\Z\grupo{A}\st{i}\hookrightarrow\grupo{A}^\ni\st{p}\twoheadrightarrow\Z\grupo{A},\;\; i(a\wedge b)=[a,b], \;\; (a,b\in A).$$
By using this central extension we can define a (right) action of
the abelian group
$\hom_{\ab(T)}(\Z\grupo{B}_\beta,\wedge^2_T\Z\grupo{A}_\alpha)$ on
the set $\hom_{\nil(T)}(\grupo{B}_\beta^\ni,\grupo{A}_\alpha^\ni)$
as
$$\varphi+\zeta\colon\grupo{B}_\beta^\ni\To\grupo{A}_\alpha^\ni\colon
x\mapsto \varphi(x)+i\zeta p(x).$$ It is easy to check that
$\varphi+\zeta$ is controlled by using that $\varphi$ and $\zeta$
are.

\begin{prop}\label{leab}
The previous action determines a linear extension of categories in
the sense of Definition \ref{esf}
$$\hom_{\ab(T)}(-,\wedge^2_T)\st{+}\To\nil(T)\st{\abb}\To\ab(T).$$
\end{prop}

\begin{proof}
We only need to check that given two controlled homomorphisms
$\varphi,\psi\colon\grupo{B}_\beta^\ni\r\grupo{A}_\alpha^\ni$ with
$\varphi^\abb=\psi^\abb$ the unique homomorphism
$\zeta\colon\Z\grupo{B}_\beta\r\Z\grupo{A}_\alpha$ with
$\varphi+\zeta=\psi$ is indeed controlled.

Let $T_v\sqcup T_v^\F$, as defined in the paragraph preceding Lemma
\ref{nice}, be a neighbourhood of $\varepsilon\in\F(T)$ in
$\hat{T}$. Since $\varphi$ and $\psi$ are controlled we can take
$T_w\subset T_v$ such that $T_w\sqcup T_w^\F$ is also another
neighborhood of $\varepsilon$ in $\hat{T}$ and
$\varphi(\beta^{-1}(T_w)),\psi(\beta^{-1}(T_w))\subset\grupo{\alpha^{-1}(T_v)}^\ni$.
Given $b\in B$, $\zeta(b)$ is a linear combination of elements
$a_1\wedge a_2$ $(a_1,a_2\in A, a_1\prec a_2)$ such that $a_i$
appears with non-trivial coefficient in the linear expansion of
either $\varphi(b)$ or $\psi(b)$ $(i=1,2)$, in particular by Lemma
\ref{nice} $\zeta(b)\in\Z\grupo{(\wedge^2\alpha)^{-1}(T_v)}$
provided $\beta(b)\in T_w$, so
$\zeta(\beta^{-1}(T_w))\subset\Z\grupo{(\wedge^2\alpha)^{-1}(T_v)}$
and therefore $\zeta$ is a controlled homomorphism.
\end{proof}


\section{Obstruction theory and cohomology of categories}\label{chV}\label{cohcat}\label{hs}


Given a small category $\C{C}$ a \emph{$\C{C}$-bimodule} $D$ is a functor
$$D\colon\C{C}^\op\times\C{C}\To\Ab.$$ For simplicity given two
morphisms $\varphi$ and $\psi$ in $\C{C}$ we write
$D(\varphi,1)=\varphi^*$ and $D(1,\psi)=\psi_*$. The cohomology groups $H^*(\C{C},D)$ were first defined in \cite{rso}. They are also the cohomology of the cochain complex $F^*(\C{C},D)$ in \cite{ah} IV.5.

In order to describe this cohomology as a functor we
consider the category $\C{Nat}$ whose objects are pairs
$(\C{C},D)$ given by a small category $\C{C}$ and a
$\C{C}$-bimodule $D$. The \emph{pull-back} of a $\C{C}$-bimodule $D$ along a functor
$\lambda\colon\C{B}\r\C{C}$ is
$\lambda^*D=D(\lambda^\op\times\lambda)$. Similarly one defines the
pull-back of a natural transformation $t\colon D\r D'$ between
$\C{C}$-bimodules as $\lambda^*t=t(\lambda^\op\times\lambda)$. A morphism $(\lambda,t)\colon (\C{C},D)\r
(\C{C}',D')$ is a functor $\lambda\colon\C{C}'\r\C{C}$ together
with a natural transformation $t\colon\lambda^*D\r D'$. The
composition in $\C{Nat}$ is given by the formula
$(\mu,u)(\lambda,t)=(\lambda\mu,u(\mu^*t))$, in particular $(\lambda,t)=(1_{\C{C}'},t)(\lambda,1_{\lambda^*D})$.

The cohomology of categories induces functors $(n\geq 0)$
$$H^n\colon\C{Nat}\To \Ab,\;\; H^n(\lambda,t)=t_*\lambda^*.$$
Here we denote $t_*=H^n(1,t)$ and $\lambda^*(\lambda,1)$.

As usual a short exact sequence of $\C{C}$-bimodules
$$D\st{t}\hookrightarrow D'\st{u}\twoheadrightarrow D''$$ induces
a Bockstein long exact sequence in cohomology
$$\cdots\r H^n(\C{C},D)\st{t_*}\To
H^n(\C{C},D')\st{u_*}\To H^n(\C{C},D'')\st{\beta}\To
H^{n+1}(\C{C},D)\r\cdots.$$

A \emph{$0$-cochain} of $\C{C}$ with coefficients in $D$ is a
function $c$ which sends an object $A$ of $\C{C}$ to an element
$c_A\in D(A,A)$. It is a $0$-cocycle if for any $\sigma\colon A\r B$ in $\C{C}$ the equation $\sigma_*c_A=\sigma^*c_B$ holds.

Suppose that $\C{A}$ is a small additive category and $D$ an
$\C{A}$-bimodule which is additive in the first variable, that is,
given a direct sum diagram as in (\ref{dsd})
$$\xymatrix{X\ar@<-.5ex>[r]_<(.3){i_1}&X\oplus
Y\ar@<-.5ex>[l]_<(.3){p_1}\ar@<.5ex>[r]^<(.37){p_2}&Y
\ar@<.5ex>[l]^<(.25){i_2}}$$ the following homomorphism is an
isomorphism $$(p_1^*,p_2^*)\colon D(X,Z)\oplus D(Y,Z)\st{\simeq}\To
D(X\oplus Y,Z).$$ We say that a $0$-cochain $f$ is \emph{additive}
if the following equality always holds
$$f(X\oplus Y)=p_1^*{i_1}_*f(X)+p_2^*{i_2}_*f(Y)\in D(X\oplus
Y,X\oplus Y).$$ It is not difficult to check that all $0$-cocycles
are additive. This is also a consequence of much more general
results proved in \cite{sncc} which point out a strong connection
between cohomology of categories and representation theory.

\begin{defn}\label{esf}
Given two $\C{C}$-bimodules $H$ and $D$, an \emph{exact
sequence for a functor} $\lambda\colon\C{B}\r\C{C}$
$$H\st{+}\To\C{B}\st{\lambda}\To\C{C}\st{\theta}\To D$$ consists
of the following:
\begin{enumerate}
\item For any morphism $f\colon X\r Y$ in $\C{B}$ the abelian
group $H(\lambda X,\lambda Y)$ acts transitively on the right of
the set $\lambda^{-1}(\lambda f)\subset\hom_{\C{B}}(X,Y)$.

\item The previous action satisfies the formula $(f+\alpha)(g+\beta)=fg+(\lambda f)_*\beta+(\lambda
g)^*\alpha$.

\item Given two objects $X$ and $Y$ in $\C{B}$ and a morphism
$f\colon\lambda X\r \lambda Y$ there is a well-defined element
$\theta_{X,Y}(f)\in D(\lambda X,\lambda Y)$ which vanishes if and
only if $f=\lambda g$ for some $g\colon X\r Y$.

\item The obstruction operator $\theta$ is a derivation, that is,
the following formula holds
$$\theta_{X,Z}(gf)=g_*\theta_{X,Y}(f)+f^*\theta_{Y,Z}(g).$$

\item Given an object $X$ in $\C{B}$ and $\alpha\in D(\lambda
X,\lambda X)$ there exists another one $Y$ such that $\lambda
X=\lambda Y$ and $\theta_{X,Y}(1_{\lambda X})=\alpha$. One can
check by using the previous axioms that the object $Y$ is well
defined up to isomorphism by $X$ and $\alpha$, and also $X$ is
determined up to isomorphism by $Y$ and $\alpha$. In these
conditions we write $X=Y+\alpha$.
\end{enumerate}
Any functor fitting into an exact sequence reflects isomorphisms,
see \cite{ah} IV.4.11.

A \emph{linear extension of categories}
$$H\st{+}\To\C{B}\st{\lambda}\To\C{C}$$ is an exact sequence $H\st{+}\r\C{B}\st{\lambda}\r\C{C}\r
0$ such that $\lambda$ induces a bijection between the sets of
objects and the action in (1) is effective.
\end{defn}

\begin{rem}\label{coh1}
Suppose that the functor $\lambda$ fitting into an exact sequence as
above induces a surjection between object sets $\lambda\colon
Ob\C{B}\twoheadrightarrow Ob\C{C}$. In this case the obstruction
operator $\theta$ determines a well defined cohomology class
$$\set{\theta}\in H^1(\C{C},D)$$ in the following way, compare \cite{ah} IV.7. Take a
splitting function $s\colon Ob\C{C}\hookrightarrow Ob\C{B}$ of
$\lambda$. The function $\tilde{s}$ sending a morphism $\sigma\colon A\r B$ in $\C{C}$ to
$\tilde{s}(\sigma)=\theta_{sA,sB}(\sigma)$ is a $1$-cocycle representing
$\set{\theta}$. Moreover any representative cocycle of this
cohomology class can be obtained as $\tilde{s}$ for an adequate
splitting $s$. Furthermore, $\set{\theta}=0$ if and only if there
exists a $0$-cocycle $c$ of $\C{B}$ with coefficients in
$\lambda^*D$ such that $\theta_{X,Y}(f)=f_*(c_X)-f^*(c_Y)$. Indeed
if such a $0$-cocycle $c$ exists the $1$-cocycle $\tilde{s}$
representing $\set{\theta}$ is the coboundary of the $0$-cochain
sending an object $A$ in $\C{C}$ to $c_{sA}$. On the other hand if
$\set{\theta}=0$ there is a splitting function $s$ such that
$\tilde{s}=0$ and we can define the $0$-cocycle $c$ as
$c_X=\theta_{X,s\lambda X}(1_{\lambda
X})$. 
\end{rem}


\begin{defn}\label{htpysys}
An $(n-1)$-reduced \emph{$T$-homotopy system} $(n\geq 2)$ is a pair
$(\cc_*, f_{n+2})$ given by a bounded chain complex $\cc_*$ in
$\ab(T)$ concentrated in degrees $\geq n$ and an $\ab(T)$-module
morphism $f_{n+2}\colon\cc_{n+2}\r\Pi_{n+1}C_{d_{n+1}}$ such that
$f_{n+2}d_{n+3}=0$ and $d_{n+2}$ coincides with the following
composite
$$\cc_{n+2}\st{f_{n+2}}\To\Pi_{n+1}C_{d_{n+1}}\st{j}\To\Pi_{n+1}(C_{d_{n+1}},S^n_{\alpha_n})=\cc_{n+1}.$$ Here we
identify the differential $d_{n+1}\colon\cc_{n+1}\r\cc_n$ with a
map between $n$-dimensional spherical objects $d_{n+1}\colon
S^n_{\alpha_{n+1}}\r S^n_{\alpha_n}$ by using Proposition
\ref{fund} so that we can take its cofiber $C_{d_{n+1}}$.

A \emph{morphism} $(\xi,\eta)\colon (\cc_*,f_{n+2})\r(\cc_*',g_{n+2})$ of $(n-1)$-reduced $T$-homotopy systems
is given by a
chain complex morphism $\xi\colon\cc_*\r\cc_*'$ and a map
$\eta\colon C_{d_{n+1}}\r C_{d_{n+1}'}$ such that $\cc_*\eta$
coincides with $\xi$ in degrees $n$ and $n+1$, and the following
square commutes
\begin{equation}\label{morfismoes}
\xymatrix{\cc_{n+2}\ar[r]^{\xi_{n+2}}\ar[d]_{f_{n+2}}&\cc_{n+2}'\ar[d]^{g_{n+2}}\\\Pi_{n+1}
C_{d_{n+1}}\ar[r]_{\eta_*}&\Pi_{n+1} C_{d_{n+1}'}}
\end{equation}
that is
$$\eta f_{n+2}=g_{n+2}\xi_{n+2}\in[S^{n+1}_{\alpha_{n+2}},C_{d_{n+1}}']^T.$$
Here $S^{n+1}_{\alpha_{n+2}}$ is the $(n+1)$-dimensional spherical
object corresponding to $\cc_{n+2}$. Let $\C{H}_n(T)$ be the
category of $(n-1)$-reduced $T$-homotopy systems.

Two morphisms $(\xi,\eta),(\xi',\eta')\colon (\cc_*,f_{n+2})\r(\cc_*',g_{n+2})$ of $(n-1)$-reduced $T$-homotopy systems
are \emph{homotopic} if there exists a chain homotopy
$\xi\simeq\xi'$ given by morphisms
$\omega_m\colon\cc_m\r\cc_{m+1}'$ such that
$$\eta+g_{n+2}\omega_{n+1}=\eta'\in[C_{d_{n+1}},C_{d'_{n+1}}]^T.$$
Here we use the action of the group
$[S^{n+1}_{\alpha_{n+1}},C_{d'_{n+1}}]^T=(\Pi_{n+1}C_{d'_{n+1}})(\cc_{n+1})$
on the set $[C_{d_{n+1}},C_{d'_{n+1}}]^T$ induced by the 
long cofiber sequence
$$S^n_{\alpha_{n+1}}\st{d_{n+1}}\To S^n_{\alpha_n}\rightarrowtail C_{d_{n+1}}\r S^{n+1}_{\alpha_{n+1}}\r\cdots,$$
see \cite{ah} II.8.8.

The homotopy relation is a natural equivalence relation in he
category of $(n-1)$-reduced $T$-homotopy systems so the quotient
category $\C{H}_n(T)/\simeq$ is well defined.
\end{defn}

The general definition of homotopy systems in homological
cofibration categories appears in \cite{cfhh} VI.1. Using the
terminology of \cite{cfhh} we have defined above $(n-1)$-reduced
twisted homotopy systems of order $n+2$ satisfying the cocycle
condition in $\C{Topp}^T$. The homotopy systems we will use here are
closer to those used in \cite{ah} IX, since we will always be
in the $1$-reduced case.

Let $\C{CW}^T_n$ be the category of $(n-1)$-reduced $T$-complexes.
There is a canonical functor
$$r\colon\C{CW}^T_n\To\C{H}_n(T)$$ which sends a $T$-complex $X$
to the pair $rX=(\cc_*X,f_{n+2})$ where $f_{n+2}\colon
\cc_{n+2}X\r\Pi_{n+1}X^{n+1}$ is induced by the pasting map of
$(n+2)$-cells $S^{n+1}_{\alpha_{n+2}}\r X^{n+1}$. Notice that
$C_{d_{n+1}}=X^{n+1}$. We will not discuss here the cylinders in
$\C{H}_n(T)$ which give rise to the homotopy relation defined above,
however we point out that the functor $r$ preserves cylinders so it
factors through the homotopy categories
\begin{equation}\label{funtoR}
r\colon\C{CW}^T_n/\!\simeq\;\To\C{H}_n(T)/\!\simeq,
\end{equation}
see \cite{cfhh} VI.1.7 for further details.

Let us write $\C{A}^2_n(T)$ for the homotopy category of
$(n-1)$-reduced $T$-complexes of dimension $\leq n+2$ and
$\C{H}_n^2(T)$ for the category $(n-1)$-reduced $T$-homotopy
systems whose chain complex is concentrated in dimensions $\leq
n+2$.

\begin{prop}\label{isoob}
The functor $r$ restricts to a full functor
$$r\colon\C{A}^2_n(T)\To\C{H}_n^2(T)/\!\simeq$$ which induces a
bijection between the sets of isomorphisms classes of objects in
both categories.
\end{prop}

This proposition can be easily checked by using principal maps in
the sense of \cite{ah} V.2 and the proper homological Whitehead
theorem.

If $\chain_n(\ab(T))$ is the category of bounded chain complexes
in $\ab(T)$ concentrated in degrees $\geq n$ there is a forgetful
functor
$$\lambda\colon\C{H}_n(T)\To\chain_n(\ab(T)),\;\;\lambda(\cc_*,f_{n+2})=\cc_*.$$
This functor also factors through the homotopy categories
$$\lambda\colon\C{H}_n(T)/\!\simeq\;\To\chain_n(\ab(T))/\!\simeq.$$
Moreover $\lambda r=\cc_*$ is the proper cellular chain complex
functor.

In the rest of this section we will concentrate in describing the
properties the functor $\lambda$. We refer the reader to \cite{cfhh}
VI for the details.

There are exact sequences for the functors $\lambda$ in the sense of
Definition \ref{esf} as follows $(n\geq 3)$
\begin{equation}\label{thetacoh}
\begin{array}{c}
H^3(-,\Gamma_TH_2)\st{+}\To\C{H}_2(T)/\!\simeq\;\st{\lambda}\To
\chain_2(\ab(T))/\!\simeq\;\st{\theta^2}\To
H^4(-,\Gamma_TH_2),\\{}\\
H^{n+1}(-,H_n\otimes\Z/2)\st{+}\r\C{H}_n(T)/\!\simeq\;\st{\lambda}\r
\chain_n(\ab(T))/\!\simeq\;\st{\theta^n}\r
H^{n+2}(-,H_n\otimes\Z/2).
\end{array}
\end{equation}

These exact sequences are defined in \cite{cfhh} VI in a much more
general setting, here we will translate some of their properties
to proper homotopy theory, mainly the definition of the
obstruction operators $\theta^n$. Here we use Theorem \ref{lwm} to
give a convenient description of the bimodules involved, otherwise
these bimodules would be defined directly in terms of Whitehead
modules.

Given two $(n-1)$-reduced homotopy systems $(\cc_*,f_{n+2})$,
$(\cc_*',g_{n+2})$ and a chain homotopy class
$\xi\colon\cc_*\r\cc_*'$ we can take a principal map $\eta\colon
C_{d_{n+1}}\r C_{d_{n+1}'}$ in the sense of \cite{ah} V.2 such that
$\cc_*\eta$ coincides with $\xi$ in dimensions $n$ and $n+1$, but
the commutativity of (\ref{morfismoes}) is not guaranteed, that is
the pair $(\xi,\eta)$ need not be a morphism of homotopy systems.
However the lower square and the two triangles in the following
diagram do commute by the properties of $\eta$ and the definition of
homotopy systems, respectively,
$$\xymatrix@C=40pt{&\cc_{n+2}\ar[r]^{\xi_{n+2}}\ar[d]^{f_{n+2}}\ar[ddl]_{d_{n+2}}&\cc_{n+2}'\ar[d]_{g_{n+2}}
\ar[ddr]^{d_{n+2}}&\\
&\Pi_{n+1}C_{d_{n+1}}\ar[r]^{\eta_*}\ar[d]&\Pi_{n+1}C_{d_{n+1}'}\ar[d]^j&\\
\cc_{n+1}\ar@{=}[r]&\Pi_{n+1}(C_{d_{n+1}'},S^n_{\alpha_n})\ar[r]_{\xi_{n+1}}&\Pi_{n+1}(C_{d_{n+1}'},S^n_{\alpha_n'})\ar@{=}[r]&\cc_{n+1}'}$$
so the failure $-g_{n+2}\xi_{n+2}+\eta f_{n+2}$ in the commutativity
of the upper square, which is (\ref{morfismoes}), factors through
the kernel of $j$ which coincides with the Whitehead module
$\Gam_{n+1}C_{d_{n+1}'}$ by the exactness of the sequence of the
pair $(C_{d_{n+1}'},S^n_{\alpha_n'})$ in proper homotopy modules,
\begin{equation*}
\xymatrix{\cc_{n+2}\ar[rr]^{-g_{n+2}\xi_{n+2}+\eta f_{n+2}}\ar[rd]_\beta&&\Pi_{n+1}C_{d_{n+1}'}\\
&\Gam_{n+1}C_{d_{n+1}'}\ar@{^{(}->}[ru]&}
\end{equation*}
By Theorem \ref{lwm} $\Gam_{n+1}C_{d_{n+1}'}$ is either
$\Gamma_TH_2\cc_*'$ if $n=2$ or $H_n\cc_*'\otimes\Z/2$ if $n\geq
3$, and finally the cohomology class
$$\theta^n_{(\cc_*,f_{n+2}),(\cc_*',g_{n+2})}(\xi)\in\left\{%
\begin{array}{ll}
    H^4(\cc_*,\Gamma_TH_2\cc_*'), & \hbox{if $n=2$;} \\ & \\
    H^{n+2}(\cc_*,H_n\cc_*'\otimes\Z_2), & \hbox{if $n\geq 3$;} \\
\end{array}\right.$$
is represented by the cocycle $\beta$.

\begin{prop}\label{so}
The functors $\lambda$ are surjective on objects.
\end{prop}

\begin{proof}
Let $\cc_*$ be a bounded chain complex concentrated in dimensions
$\geq n$ for some $n\geq 2$. Since $C_{d_{n+1}}$ is $1$-connected
one can check by using the exact sequence of the pair
$(C_{d_{n+1}},S^n_{\alpha_n})$ in proper homotopy modules that the
image of
$j\colon\Pi_{n+1}C_{d_{n+1}}\r\Pi_{n+1}(C_{d_{n+1}},S^n_{\alpha_n})=\cc_{n+1}$
is $\ker d_{n+1}$ which is a f. g. free $\ab(T)$-module by Theorem
\ref{kernel}, so there exists a section $s\colon\ker
d_{n+1}\hookrightarrow\Pi_{n+1}C_{d_{n+1}}$ of the natural
projection onto the image
$p\colon\Pi_{n+1}C_{d_{n+1}}\twoheadrightarrow\ker d_{n+1}$. The
differential $d_{n+2}$ factors through the inclusion $i\colon\ker
d_{n+1}\hookrightarrow\cc_{n+1}$, that is
$d_{n+2}=i\bar{d}_{n+2}$. We claim that $(\cc_*,s\bar{d}_{n+2})$
is an $(n-1)$-reduced homotopy system. On one hand
$js\bar{d}_{n+2}=ips\bar{d}_{n+2}=i\bar{d}_{n+2}=d_{n+2}$, as
required, and on the other hand $\bar{d}_{n+2}d_{n+3}=0$ since
$i\bar{d}_{n+2}d_{n+3}=d_{n+2}d_{n+3}=0$ and $i$ is a
monomorphism, hence $s\bar{d}_{n+2}d_{n+3}=0$ is also trivial. Now
the proof is finished.
\end{proof}

By this proposition and Remark \ref{coh1} the exact sequences of
functors in (\ref{thetacoh}) define classes in cohomology of
categories
\begin{equation}\label{lasclases}
\set{\theta^n}\in\left\{%
\begin{array}{ll}
    H^1(\chain_2(\ab(T))/\!\simeq,H^4(-,\Gamma_TH_2)), & \hbox{if $n=2$;} \\ & \\
    H^1(\chain_n(\ab(T))/\!\simeq,H^{n+2}(-,H_n\otimes\Z/2)), & \hbox{if $n\geq 3$;} \\
\end{array}\right.
\end{equation}
The suspension of chain complexes
\begin{equation}\label{scc}
\Sigma\colon\chain_n(\ab(T))/\!\simeq\;\To\chain_{n+1}(\ab(T))/\!\simeq
\end{equation}
and the suspension of $T$-complexes induce suspension functors
$$\Sigma\colon\C{H}_n(T)/\!\simeq\;\To\C{H}_{n+1}(T)/\!\simeq,\;\; \Sigma (\cc_*,f_{n+2})=(\Sigma\cc_*,\Sigma
f_{n+2}),$$ such that $\Sigma r=r\Sigma$ and
$\Sigma\lambda=\lambda\Sigma$. It is easy to see by using the
description of $\theta^n$ given above and Theorem \ref{lwm} that
these suspension functors are also compatible with the obstruction
operators, that is the following equality holds in
$H^{n+3}(\Sigma\cc_*,H_{n+1}\Sigma\cc_*'\otimes\Z/2)=H^{n+2}(\cc_*,H_n\cc_*'\otimes\Z/2)$
\begin{equation*}
\theta^{n+1}_{\Sigma(\cc_*,f_{n+2}),\Sigma(\cc_*',g_{n+2})}(\Sigma\xi)=
\left\{%
\begin{array}{ll}
    (\sigma_T)_*\theta^2_{(\cc_*,f_4),(\cc_*',g_4)}(\xi), & \hbox{if $n=2$;} \\ & \\
    \theta^n_{(\cc_*,f_{n+2}),(\cc_*',g_{n+2})}(\xi), & \hbox{if $n\geq 3$.} \\
\end{array}\right.
\end{equation*}
Therefore the cohomology classes defined by the exact sequences in
(\ref{thetacoh}) are compatible with the suspension functor in the
sense that
\begin{equation}\label{compatibles}
(\sigma_T)_*\set{\theta^2}=\Sigma^*\set{\theta^3}\text{ and
}\set{\theta^n}=\Sigma^*\set{\theta^{n+1}}\text{ for }n\geq3.
\end{equation}
Notice that the suspension of chain complexes (\ref{scc}) is an
isomorphism of categories, therefore the induced homomorphisms
$\Sigma^*$ above are abelian group isomorphisms, in particular
$\set{\theta^2}$ determines $\set{\theta^n}$ for all $n\geq 3$ and
all cohomology classes ${\theta^n}$ $(n\geq 3)$ correspond to each
other. From now on we shall simply write $\theta$ for the
obstruction operator $\theta^n$, or $\theta^T$ if we want to specify
the base tree.

Given a proper map $f\colon T\r T'$ between trees the ``change of
tree'' functors $\uf{F}^f\colon\ab(T)\r\ab(T')$ and
$f_*\colon\C{Topp}^T_c/\!\simeq\r\C{Topp}^{T'}_c/\!\simeq$ in
Proposition \ref{induce} and (\ref{cbase}), respectively, give rise
to maps between exact sequences for functors in the sense of
\cite{ah} IV.4.13 $(n\geq 3)$
\begin{equation}\label{lasma}
\begin{array}{c}
\xymatrix{H^3(-,\Gamma_TH_2)\ar[r]^<(.25)+\ar[d]_{\bar{\uf{F}}^f}&\C{H}_2(T)/\!\simeq\ar[r]^<(.2)\lambda\ar[d]^{(\uf{F}^f,f_*)}
&\C{chain}_2(\ab(T))/\!\simeq\ar[r]^<(.25){\theta^T}\ar[d]^{\uf{F}^f}&H^4(-,\Gamma_TH_2)\ar[d]^{\bar{\uf{F}}^f}\\
H^3(-,\Gamma_{T'}H_2)\ar[r]^<(.25)+&\C{H}_2(T')/\!\simeq\ar[r]^<(.2)\lambda
&\C{chain}_2(\ab(T'))/\!\simeq\ar[r]^<(.25){\theta^{T'}}&H^4(-,\Gamma_{T'}H_2)}\\{}\\
\begin{small}
\xymatrix@C=11pt{H^{n+1}(-,H_n\otimes\Z/2)\ar[r]^<(.15)+\ar[d]_{\bar{\uf{F}}^f}&\C{H}_n(T)/\!\simeq\ar[r]^<(.1)\lambda\ar[d]^{(\uf{F}^f,f_*)}
&\C{chain}_n(\ab(T))/\!\simeq\ar[r]^<(.15){\theta^T}\ar[d]^{\uf{F}^f}&H^{n+2}(-,H_n\otimes\Z/2)\ar[d]^{\bar{\uf{F}}^f}\\
H^{n+1}(-,H_n\otimes\Z/2)\ar[r]^<(.15)+&\C{H}_n(T')/\!\simeq\ar[r]^<(.1)\lambda
&\C{chain}_n(\ab(T'))/\!\simeq\ar[r]^<(.15){\theta^{T'}}&H^{n+2}(-,H_n\otimes\Z/2)}
\end{small}
\end{array}
\end{equation}
Here the natural transformations $\bar{\uf{F}}^f\colon
H^m(-,\Gamma_TH_2)\r H^m(\uf{F}^f,\Gamma_{T'}H_2\uf{F}^f)$ between
bimodules over $\C{chain}_2(\ab(T))/\!\simeq$ are induced by taking
$H^m$ on the following composition of cochain homomorphisms
\begin{equation}\label{tecambio}
\xymatrix{\hom_{\ab(T)}(\cc_*,\Gamma_TH_2\cc_*)\ar[r]^<(.18){\uf{F}^f_*}&\hom_{\ab(T')}(\uf{F}^f\cc_*,\uf{F}^f_*\Gamma_TH_2\cc_*)\ar[d]\\
\hom_{\ab(T')}(\uf{F}^f\cc_*,\Gamma_{T'}H_2\uf{F}^f\cc_*)\ar@{=}[r]&
\hom_{\ab(T')}(\uf{F}^f\cc_*,\Gamma_{T'}\uf{F}^f_*H_2\cc_*)}
\end{equation}
Here the first arrow is given by the functor $\uf{F}_*^f$ extending
$\uf{F}^f$ to the category of all $\ab(T)$-modules in the sense of
(\ref{yonedacom}), the second one is given by the natural
transformation $\uf{F}^f_*\Gamma_T\r\Gamma_{T'}\uf{F}^f_*$, see
Proposition \ref{quadpropri2}, and for the equality we use the
facts that $\uf{F}^f_*$ is right-exact and $\cc_*$ is concentrated
in dimensions $\geq 2$. Similarly the natural transformations
$\bar{\uf{F}}^f\colon H^m(-,H_n\otimes\Z/2)\r
H^m(\uf{F}^f,(H_n\uf{F}^f)\otimes\Z/2)$ between bimodules over
$\C{chain}_n(\ab(T))/\!\simeq$ are given by the cochain homomorphism
\begin{equation}\label{tecambio2}
\xymatrix{\hom_{\ab(T)}(\cc_*,H_n\cc_*\otimes\Z/2)\ar[r]^<(.18){\uf{F}^f_*}&\hom_{\ab(T')}(\uf{F}^f\cc_*,\uf{F}^f_*(H_n\cc_*\otimes\Z/2))\ar@{=}[d]\\
\hom_{\ab(T')}(\uf{F}^f\cc_*,H_n(\uf{F}^f\cc_*)\otimes\Z/2)\ar@{=}[r]&
\hom_{\ab(T')}(\uf{F}^f\cc_*,(\uf{F}^f_*H_n\cc_*)\otimes\Z/2)}
\end{equation}
As a consequence of this we have the following result.

\begin{prop}\label{concola}
Given a proper map $f\colon T\r T'$ between trees the cohomology
classes in (\ref{lasclases}) are compatible with the induced ``change of
tree'' functor $\uf{F}^f\colon\ab(T)\r\ab(T')$, i. e. with the
notation above
$$(\uf{F}^f)^*\set{\theta^{T'}}=\bar{\uf{F}}^f_*\set{\theta^T}.$$
\end{prop}

\begin{rem}\label{concola2}
We point out that $\bar{\uf{F}}^f_*$ is an isomorphism provided $f$
induces an injection $\F(f)\colon\F(T)\hookrightarrow\F(T')$. In
fact (\ref{tecambio}) and (\ref{tecambio2}) are cochain isomorphisms
in this case, see Propositions  \ref{fullis} and \ref{quadpropri}.
Therefore the cohomology class $\set{\theta^{T'}}$ determines
$\set{\theta^{T}}$ for all trees $T$ such that $\F(T)$ is
homeomorphic to a subspace of $\F(T')$.
\end{rem}

\chapter{James-Hopf and cohomology invariants in proper homotopy theory}\label{chVI}

In this chapter we introduce new interrelated proper homotopy and cohomology invariants: James-Hopf invariants, interior cup-producs, cup-products for $T$-homotopy systems, $T$-complexes and chain complexes, and reduced or stable versions of them. These invariants are constructed by using the controlled quadratic functors from Chapter \ref{qf}. We show that the (reduced) cup-product invariant for $T$-complexes is related to obstruction theory through a change of coefficient given by  natural transformations between controlled quadratic functors. The cup-product for chain complexes measures the obstruction to the existence of a $T$-complex with a co-H-multiplication and a prescribed proper cellular chain complex, as we show in the last section of this chapter. Such an onstruction is trivial in classical homotopy theory since classical Moore spaces are much better behaved than in proper homotopy theory. In proper homotopy theory the cup-product for chain complexes need not be trivial since we are in projective dimension $2$, and not just $1$ as in the classical case. Indeed the computation of the cup-product for chain complexes is one of the main problems that we tackle in this paper, see Chapter \ref{chVII}. This computation is the basis of the main results of this paper on the proper homotopy classification of locally compact $A^2_n$-polyhedra, see Chapter \ref{chVIII}.

\section{James-Hopf invariants}

\begin{defn}\label{cohspace}
Recall that a \emph{co-H-space} in $\C{Topp}^T_c$ is an object $X$
endowed with a proper map $\mu\colon X\r X\vee X$ under $T$, the
\emph{co-H-multiplication}, such that the composition with the
projections $p_1,p_2\colon X\vee X\r X$ are homotopic to the
identity $p_i\mu\simeq 1_X$ $(i=1,2)$. Spherical objects of
dimension $\geq 1$ are co-H-spaces. More generally, it is well known
that all suspensions are canonically co-H-spaces, see \cite{ah}
II.6.16.
\end{defn}

\begin{defn}\label{JH1}
The \emph{James-Hopf invariant} of a $T$-complex $X$ is the unique
$\ab(T)$-module morphism
$$\gamma_2\colon\Pi_3\Sigma X\To\otimes^2_T\hh_2\Sigma X$$ such that if
$i_1,i_2,\mu\colon\Sigma X\r\Sigma X\vee\Sigma X$ are the
inclusions of the factors of the coproduct and the
co-H-multiplication of the suspension $\Sigma X$, respectively,
then
$$i_3i_{12}\gamma_2=\mu_*-{i_2}_*-{i_1}_*\colon\Pi_3\Sigma X\To\Pi_3(\Sigma X\vee\Sigma X).$$
The morphism $\gamma_2$ exists and is unique by Proposition
\ref{piaditiv}, and it is natural in $X$.
\end{defn}

\begin{prop}\label{conmuta1}
For any $T$-complex $X$ the following diagram is commutative
$$\xymatrix@C=10pt{\Gamma_T\h_2\Sigma X\ar[rr]^{\tau_T}\ar[dr]_{i_3}&&\otimes^2_T\hh_2\Sigma X\\
&\Pi_3\Sigma X\ar[ru]_{\gamma_2}&}$$
\end{prop}

\begin{proof}\renewcommand{\theequation}{\alph{equation}}\setcounter{equation}{0}
On one hand we have that
\begin{eqnarray*}
  i_3i_{12}\gamma_2i_3 &=& (\Pi_3\mu-\Pi_3i_2-\Pi_3i_1)i_3 \\
   &=& i_3(\Gamma_T\h_2\mu-\Gamma_T\h_2i_2-\Gamma_T\h_2i_1) \\
   &=& i_3(\Gamma_T(i_1+i_2)-\Gamma_Ti_2-\Gamma_Ti_1).
\end{eqnarray*}
For the first equality we use the Definition \ref{JH1}, for the
second one we use the naturality of Whitehead's long exact sequence
and Theorem \ref{lwm}, for the third one we use the natural
identifications $\hh_2(\Sigma X\vee\Sigma X)=\hh_2\Sigma
X\oplus\hh_2\Sigma X$, $\hh_2\mu=i_1+i_2$ and $\hh_2i_k=i_k$
$(k=1,2)$.

On the other hand we have that
\begin{eqnarray}
  i_3i_{12}\tau_T &=& i_3i_{12}r_{12}\Gamma_T(i_1+i_2) \\
\nonumber   &=& i_3(\Gamma_T(i_1p_1+i_2p_2)-\Gamma_T(i_1p_1)-\Gamma_T(i_2p_2))\Gamma_T(i_1+i_2) \\
\nonumber   &=& i_3(\Gamma_T(i_1+i_2)-\Gamma_Ti_1-\Gamma_Ti_2).
\end{eqnarray}
Here for the first equality we use the Proposition
\ref{quadpropri2}, for the second one we use the definitions of
$i_{12}$ and $r_{12}$ in Definition \ref{dgqf}, and finally for the
third one we use the properties of the structural morphisms of a
direct sum.

Now that we have proved the equality
$i_3i_{12}\gamma_2i_3=i_3i_{12}\tau_T$ the proposition follows
from the fact that $i_3i_{12}$ is a monomorphism, see Proposition
\ref{piaditiv}.
\end{proof}

\setcounter{equation}{1}


\begin{defn}\label{JHS1}
Given a $T$-complex X and $n\geq 2$ we define the \emph{stable
James-Hopf invariant} as the unique $\ab(T)$-module morphism
$\gamma_2^n$ that fits into the following commutative diagram
$$\xymatrix{\Pi_{n+2}\Sigma^nX\ar[r]^<(.2){\gamma^n_2}&\hat{\otimes}^2_T\hh_{n+1}\Sigma^nX\\
\Pi_3\Sigma
X\ar[r]^<(.3){\gamma_2}\ar@{->>}[u]^{\Sigma^{n-1}}&\otimes^2_T\hh_2\Sigma
X\ar@{->>}[u]_{\bar{\sigma}_T}\\
\otimes^2_T\hh_2\Sigma
X\ar@{=}[r]\ar[u]^{i_3[-,-]_T}&\otimes^2_T\hh_2\Sigma
X\ar[u]_{\tau_T[-,-]_T}}$$ Here we use the suspension isomorphism in
proper homology $\hh_2\Sigma X=\hh_{n+1}\Sigma^nX$ as an
identification. The commutativity of the lower square follows from
Proposition \ref{conmuta1}. Moreover, the morphism $\gamma^n_2$
exists and is unique since the column in the left is exact by
Proposition \ref{nucleosusp} and Freudenthal's suspension theorem,
and $\bar{\sigma}_T\tau_T[-,-]_T=\bar{\tau}_T\sigma_T[-,-]_T=0$, see
(\ref{d2}).
\end{defn}

As a consequence of Proposition \ref{conmuta1}, the properties of
diagram (\ref{d2}) and Theorem \ref{lwm} we have the following
result.

\begin{prop}\label{conmuta2}
The following diagram commutes for any $T$-complex $X$ and
$n\geq 2$
$$\xymatrix@C=10pt{\h_{n+1}\Sigma^n X\otimes\Z_2\ar[rr]^{\bar{\tau}_T}\ar[dr]_{i_{n+2}}&&\hat{\otimes}^2_T\hh_{n+1}\Sigma^n X\\
&\Pi_{n+2}\Sigma^n X\ar[ru]_{\gamma^n_2}&}$$
\end{prop}


\begin{defn}\label{JH2}
Given two $T$-complexes $X^2$ and $Y$ with $\dim X^2\leq 2$ we
define the \emph{James-Hopf invariant} of a map $f\colon\Sigma
X^2\r\Sigma Y$ as the unique element
$$\gamma_2(f)\in H^3(\Sigma X^2,\otimes^2_T\hh_2\Sigma Y)$$ such that
the following equality holds in $[\Sigma X^2,\Sigma Y\vee\Sigma
Y]^T$,
$$\mu f=i_1f+i_2f+j(i_3i_{12})_*\gamma_2(f).$$
Here $i_1,i_2,\mu\colon\Sigma X^2\r\Sigma X^2\vee\Sigma X^2$ are
the inclusions of the factors of the coproduct and the
co-H-multiplication of the suspension $\Sigma X^2$. Moreover,
$\gamma_2(f)$ exists and is unique by Proposition
\ref{portutatis}.
\end{defn}


If $X=S^2_\alpha$ the function $$\gamma_2\colon[\Sigma X^2,\Sigma
Y]\To H^3(\Sigma X^2,\otimes^2_T\hh_2\Sigma Y)$$ coincides with the
$\ab(T)$-module morphism in Definition \ref{JH1} evaluated at the
object in $\ab(T)$ corresponding to the spherical object
$S^2_\alpha$. However in general this function $\gamma_2$ need not
be a homomorphism but a quadratic function in the sense of
Definition \ref{dgqf}. In order to check this fact we define the
following new proper homotopy operation.

\begin{defn}\label{cupext1}
Given three $T$-complexes $X^2$, $Y$ and $Z$ with $\dim X^2\leq 2$
and two maps $f\colon\Sigma X^2\r\Sigma Y$ and $g\colon\Sigma
X^2\r\Sigma Z$ we define their \emph{interior cup-product}
$$f\cup
g\in H^3(\Sigma X^2,\hh_2\Sigma Y\otimes_T\hh_2\Sigma Z)$$ as the
unique element satisfying the following equality in $[\Sigma
X^2,\Sigma Y\vee\Sigma Z]^T$,
$$i_2f+i_1g=i_1g+i_2f+j(i_3i_{12})_*(f\cup g).$$ 
The element $f\cup g$ exists and is unique by Proposition
\ref{portutatis}. Moreover, by Proposition \ref{yanose} the proper
homotopy group $[\Sigma X^2,\Sigma Y\vee\Sigma Z]^T$ has nilpotency
class $2$, so the interior cup-product is bilinear and factors
through the cokernel of $j$ in Proposition \ref{yanose} by a
homomorphism
$$H^2(\Sigma X^2,\h_2\Sigma Y)\otimes H^2(\Sigma X^2,\h_2\Sigma
Z)\st{\cup}\To H^3(\Sigma X^2,\hh_2\Sigma Y\otimes_T\hh_2\Sigma
Z).$$ Here we use the Hurewicz isomorphism as an identification. The
suspension isomorphisms in proper (co)homology allow us to regard
this homomorphism as follows
$$H^1(
X^2,\h_1Y)\otimes H^1(X^2,\h_1Z)\st{\cup}\To H^2( X^2,\hh_1
Y\otimes_T\hh_1Z).$$ Furthermore, by using \cite{comcup} one can
check that this is a generalization of the cup-product in ordinary
homotopy theory. We point out that a complete definition of a
cup-product operation for the proper cohomology in Chapter
\ref{ephai} is not known.
\end{defn}

\begin{prop}\label{cuadra1}
Given two $T$-complexes $X^2$ and $Y$ with $\dim X^2\leq 2$ and
two maps $f,g\colon\Sigma X^2\r\Sigma Y$ the following equality is
satisfied
$$\gamma_2(f+g)=\gamma_2(f)+\gamma_2(g)+f\cup g\in H^3(\Sigma
X^2,\otimes^2_T\hh_2\Sigma Y).$$
\end{prop}

\begin{proof}
The result follows from the equalities in $[\Sigma X^2,\Sigma
Y\vee \Sigma Y]^T$
\begin{eqnarray*}
  i_1(f+g)+i_2(f+g)&&\\+j(i_3i_{12})_*(\gamma_2(f)+\gamma_2(g)+f\cup g) &=& i_1f+i_1g+i_2f+i_2g\\&&+j(i_3i_{12})_*(\gamma_2(f)+\gamma_2(g)+f\cup g) \\
   &=& i_1f+i_2f+i_1g+i_2g\\&&+j(i_3i_{12})_*(\gamma_2(f)+\gamma_2(g)) \\
   &=& \mu f+\mu g \\
   &=& \mu(f+g).
\end{eqnarray*}
Here we use that $j$ in Proposition \ref{yanose} is central and
Definitions \ref{JH2} and \ref{cupext1}.
\end{proof}

This result together with \cite{hth} A.10.2 (f) can be used to
prove that the interior cup-product that we have defined in proper
homotopy theory generalizes the corresponding one in ordinary
homotopy theory, which appears for example in \cite{hth} A.1.18.

\begin{prop}\label{conmuta3}
Given two $T$-complexes $X^2$ and Y with $\dim X^2\leq 2$ the
following diagram commutes
$$\xymatrix{H^3(\Sigma X^2,\Gamma_T\hh_2\Sigma
Y)\ar[r]^<(.18){{i_3}_*}\ar[d]_{{\tau_T}_*}&H^3(\Sigma
X^2,\Pi_3\Sigma Y)\ar@{^{(}->}[d]^j\\H^3(\Sigma
X^2,\otimes^2_T\hh_2\Sigma Y)&[\Sigma X^2,\Sigma
Y]^T\ar[l]^<(.18){\gamma_2}}$$
\end{prop}

\begin{proof}
Given $a\in H^3(\Sigma X^2,\Gamma_T\hh_2\Sigma Y)$ we have the
equalities
\begin{eqnarray*}
  j(i_3i_{12})_*\gamma_2(j{i_3}_*(a)) &=& -i_2j{i_3}_*(a)-i_1j{i_3}_*(a)+\mu j{i_3}_*(a) \\
   &=& j[i_3(-\Gamma_Ti_2-\Gamma_Ti_1+\Gamma_T(i_1+i_2))]_*(a) \\
   &=& j(i_3i_{12}\tau_T)_*(a) \\
   &=& j(i_3i_{12})_*(\tau_T)_*(a).
\end{eqnarray*}
For the first one we use Definition \ref{JH2}; for the second one we
use the naturality of $j$ in Proposition \ref{yanose} and
Whitehead's long exact sequence, as well as Theorem \ref{lwm};
finally for the third one we use equalities (a) in the proof of
Proposition \ref{conmuta1}.

The homomorphism $j$ is injective by Proposition \ref{yanose} and
$(i_3i_{12})_*$ is injective because $i_3i_{12}$ is a split
monomorphism, see Proposition \ref{piaditiv}, therefore the
proposition follows.
\end{proof}


\begin{defn}\label{JHS2}
Given two $T$-complexes $X^2$ and $Y$ with $\dim X^2\leq 2$ we
define the \emph{stable James-Hopf invariant} $\gamma_2^n$ $(n\geq
2)$ as the unique function which fits into the following
commutative diagram
$$\xymatrix{[\Sigma^nX^2,\Sigma^nY]^T\ar[r]^<(.2){\gamma^n_2}&H^{n+2}(\Sigma^nX^2,\hat{\otimes}^2_T\hh_{n+1}\Sigma^nY)\\
[\Sigma X^2,\Sigma
Y]^T\ar[r]^<(.3){\gamma_2}\ar@{->>}[u]^{\Sigma^{n-1}}&H^3(\Sigma
X^2,\otimes^2_T\hh_2\Sigma
Y)\ar@{->>}[u]_{(\bar{\sigma}_T)_*}\\
H^3(\Sigma X^2,\otimes^2_T\hh_2\Sigma
Y)\ar@{=}[r]\ar[u]^{j(i_3[-,-]_T)_*}&H^3(\Sigma
X^2,\otimes^2_T\hh_2\Sigma Y)\ar[u]_{(\tau_T[-,-]_T)_*}}$$ The
commutativity of the lower square follows from Proposition
\ref{conmuta3}. Here we use the suspension isomorphisms in
(co)homology as identifications.
\end{defn}

It is not as clear as in Definition \ref{JHS1} that $\gamma_2^n$
exists and is unique, but we check this in the following lemma.

\begin{lem}
Indeed the function $\gamma_2^n$ is well defined.
\end{lem}

\begin{proof}
The column in the left is exact by Proposition \ref{porcelanosa},
therefore we just need to check that always
$(\bar{\sigma}_T)_*\gamma_2(f+j(i_3[-,-]_T)_*(\alpha))=(\bar{\sigma}_T)_*\gamma_2(f)$.
The properties of the interior cup-product defined in Definition
\ref{cupext1} and Proposition \ref{yanose} show that $f\cup
j(i_3[-,-]_T)_*(\alpha)=0$, hence by Proposition \ref{cuadra1} we
have that
$\gamma_2(f+j(i_3[-,-]_T)_*(\alpha))=\gamma_2(f)+\gamma_2(j(i_3[-,-]_T)_*(\alpha))$,
but using the commutativity of the lower square we get
$(\bar{\sigma}_T)_*\gamma_2(j(i_3[-,-]_T)_*(\alpha))=(\bar{\sigma}_T)_*(\tau_T[-,-]_T)_*(\alpha)=0$
by (\ref{d2}), hence we are done.
\end{proof}

Again $\gamma^n_2$ $(n\geq 2)$ in Definition \ref{JHS2}
generalizes the stable James-Hopf invariant in Definition
\ref{JHS1} when $X^2=S^2_\alpha$ is a $2$-dimensional spherical
object, and this new $\gamma^n_2$ is not in general a homomorphism
but a quadratic map. To check this last statement we define a new
stable proper homotopy operation.

\begin{defn}\label{cupext2}
Given three $T$-complexes $X^2$, $Y$ and $Z$ with $\dim X^2\leq 2$
and maps $f\colon\Sigma^nX^2\r\Sigma^nY$,
$g\colon\Sigma^nX^2\r\Sigma^nZ$ $(n\geq 2)$ then by using the
suspension isomorphisms in proper (co)homology we define their
\emph{interior cup-product}
$$f\cup g\in
H^{n+2}(\Sigma^nX^2,\hh_{n+1}\Sigma^nY\otimes_T\hh_{n+1}\Sigma^nZ)$$
as the interior cup-product $f'\cup g'$ in the sense of Definition
\ref{cupext1} of two maps $f'\colon\Sigma X^2\r\Sigma Y$ and
$g'\colon\Sigma X^2\r\Sigma Z$ with $\Sigma^{n-1}f'=f$ and
$\Sigma^{n-1}g'=g$.

The properties of the interior cup-product in Definition
\ref{cupext1} together with Proposition \ref{porcelanosa} guarantee
that this new interior cup-product is well defined, bilinear, and
factors through the cokernel of $j$ in Proposition \ref{yanose} by a
homomorphism
$$\xymatrix{H^{n+1}(\Sigma^nX^2,\hh_{n+1}\Sigma^nY)\otimes
H^{n+1}(\Sigma^nX^2,\hh_{n+1}\Sigma^nZ)\ar[d]^{\cup}\\
H^{n+2}(\Sigma^nX^2,\hh_{n+1}\Sigma^nY\otimes_T\hh_{n+1}\Sigma^nZ)}.$$
Moreover, by Proposition \ref{cuadra1} and Definition \ref{JHS2}
we get the following result.
\end{defn}

\begin{prop}
Given two $T$-complexes $X^2$ and $Y$ with $\dim X^2\leq 2$ and
two maps $f,g\colon\Sigma^nX^2\r\Sigma^nY$ $(n\geq 2)$ the
following equality is satisfied
$$\gamma^n_2(f+g)=\gamma^n_2(f)+\gamma^n_2(g)+(\bar{\sigma}_T)_*(f\cup g)\in H^{n+2}(\Sigma^nX,\hat{\otimes}^2_T\hh_{n+1}\Sigma^nY).$$
\end{prop}

Finally we will show the behaviour of the James-Hopf invariants in
Definitions \ref{JH2} and \ref{JHS2} with respect to compositions
of maps.

\begin{prop}\label{composicio1}
Given three $T$-complexes $X^2$, $Y^2$ and $Z$ with dimensions $\dim X^2$,
$\dim Y^2\leq 2$ and maps $f\colon\Sigma X^2\r\Sigma Y^2$,
$g\colon\Sigma Y^2\r\Sigma Z$ the following equality holds
$$\gamma_2(gf)=f^*\gamma_2(g)+g_*\gamma_2(f)\in H^3(\Sigma X^2,\otimes^2_T\hh_2\Sigma Z).$$
Here $f^*=H^3(f,\otimes^2_T\hh_2\Sigma Z)$ and $g_*=H^3(\Sigma
X^2,\otimes^2_T\hh_2g)$.
\end{prop}

\begin{proof}\renewcommand{\theequation}{\alph{equation}}\setcounter{equation}{0}
It follows from the equalities
\begin{eqnarray}
  (i_1+i_2)gf &=& (i_1g+i_2g+j(i_3i_{12})_*\gamma_2(g))f \\
   &=& (i_1g+i_2g)f+j(i_3i_{12})_*f^*\gamma_2(g) \\
\nonumber   &=& (g\vee g)(i_1+i_2)f+j(i_3i_{12})_*f^*\gamma_2(g) \\
   &=& (g\vee g)(i_1f+i_2f+j(i_3i_{12})_*\gamma_2(f))+j(i_3i_{12})_*f^*\gamma_2(g) \\
   &=&
   i_1gf+i_2gf+j(i_3i_{12})_*g_*\gamma_2(f)+j(i_3i_{12})_*f^*\gamma_2(g).
\end{eqnarray}
For (a) and (c) we use Definition \ref{JH2}, and for (b) and (d) we
use Proposition \ref{yanose} and the naturality of $j$ and
$i_3i_{12}$ in Propositions \ref{yanose} and \ref{piaditiv},
respectively.
\end{proof}

\setcounter{equation}{1}

The following corollary follows directly from Definitions \ref{JHS2}
and \ref{cupext2}.

\begin{cor}\label{composicio2}
Given three $T$-complexes $X^2$,
$Y^2$ and $Z$ with dimensions $\dim X^2$, $\dim Y^2$ $\leq 2$ and maps
$f\colon\Sigma^n X^2\r\Sigma^n Y^2$, $g\colon\Sigma^n
Y^2\r\Sigma^n Z$ the following equality holds $(n\geq 2)$
$$\gamma^n_2(gf)=f^*\gamma^n_2(g)+g_*\gamma^n_2(f)\in H^{n+2}(\Sigma^n X^2,\hat{\otimes}^2_T\hh_{n+1}\Sigma^n Z).$$
Here $f^*=H^{n+2}(f,\hat{\otimes}^2_T\hh_{n+1}\Sigma^n Z)$ and
$g_*=H^{n+2}(\Sigma^n X^2,\hat{\otimes}^2_T\hh_{n+1}g)$.
\end{cor}

\section{Cup-product invariants}\label{cpi}

\begin{defn}
The \emph{cup-product invariant} of a $1$-reduced homotopy system
$(\cc_*,f_4)$ is the cohomology class $$\cup_{(\cc_*,f_4)}\in
H^4(\cc_*,\otimes^2_TH_2\cc_*)$$ represented by the cocycle
$$\cc_4\st{f_4}\To\Pi_3C_{d_3}\st{\gamma_2}\To\otimes^2_T
H_2\cc_*.$$ This is indeed a cocycle since $f_4$ already is.
Moreover, notice that in order to define the cup-product invariant
we have to choose a suspension structure in the cofiber $C_{d_3}$,
which exists by Freudental's suspension theorem, otherwise the
James-Hopf invariant $\gamma_2$ would not be defined.
\end{defn}

The cup product invariant is a $0$-cocycle in $\C{H}_2(T)/\!\simeq$.

\begin{prop}\label{cupnat}
Given a morphism $(\xi,\eta)\colon (\cc_*,f_4)\r (\cc_*',g_4)$ in
$\C{H}_2(T)/\!\simeq$ the following equality holds
$$\xi_*\cup_{(\cc_*,f_4)}=\xi^*\cup_{(\cc_*',g_4)}\in H^4(\cc_*,\otimes^2_TH_2\cc_*').$$
Here $\xi^*=H^4(\xi,\otimes^2_TH_2\cc_*')$ and
$\xi_*=H^4(\cc_*,\otimes^2_TH_2\xi)$.
\end{prop}

\begin{proof}
Recall from Definition \ref{htpysys} that, by using the
isomorphisms in Proposition \ref{fund} as identifications, we have
the equality $\eta f_4=g_4\xi_4\in
[S^3_{\alpha_4},C_{d_{n+1}'}]^T$ where $S^3_{\alpha_4}$
corresponds to $\cc_4$, hence $\gamma_2(\eta
f_4)=\gamma_2(g_4\xi_4)$. Moreover, by Proposition
\ref{composicio1} we have that
$$\gamma_2(\eta f_4)=f_4^*\gamma_2(\eta)+\eta_*\gamma_2(f_4)\st{\text{(a)}}=f_4^*\gamma_2(\eta)+(\otimes^2_TH_2\xi)\gamma_2f_4,$$
$$\gamma_2(g_4\xi_4)=\xi_4^*\gamma_2(g_4)+{g_4}_*\gamma_2(\xi_4)\st{\text{(b)}}=\gamma_2g_4\xi_4.$$
For (a) we use that $\eta$ induces $\xi$ in dimensions $2$ and
$3$, and for (b) we use that $\gamma_2(\xi_4)=0$ since the target
of $\xi_4$ is a $3$-dimensional spherical object, whose
$2$-dimensional proper homology module is trivial.

Moreover, $f_4^*\gamma_2(\eta)$ is a coboundary since by
Definition \ref{htpysys} the morphism induced in cohomology
$H^*(-,\otimes^2_TH_2\cc_*')$ by the composite
$$S^3_{\alpha_4}\st{f_4}\To C_{d_3}\st{q}\To S^3_{\alpha_3}$$ coincides
with
$$(\otimes^2_TH_2\cc_*')(d_4)\colon(\otimes^2_TH_2\cc_*')(\cc_3)\st{q^*}\twoheadrightarrow H^3(C_{d_3},\otimes^2_TH_2\cc_*')\st{f_4^*}\To
(\otimes^2_TH_2\cc_*')(\cc_4),$$ and $q^*$ is surjective because
$q$ is part of the long cofiber sequence
$$S^2_{\alpha_3}\st{d_3}\r S^2_{\alpha_2}\rightarrowtail C_{d_3}\st{q}\r S^3_{\alpha_3}\st{d_3}\r S^3_{\alpha_2}\r\cdots.$$
Therefore the cocycle $(\otimes^2_TH_2\xi)\gamma_2f_4$ represents
the same class as $\gamma_2g_4\xi_4$ in cohomology, and this is
what we exactly wanted to prove.
\end{proof}

\begin{cor}
The cup-product invariant of a $1$-reduced homotopy system
$(\cc_*,f_4)$ does not depend on the suspension structure chosen
in $C_{d_3}$ for its definition.
\end{cor}

This can be regarded as a special case of Proposition \ref{cupnat}
for the identity morphism in $(\cc_*,f_4)$.

\begin{defn}\label{redcup}
The \emph{reduced cup-product invariant} of an $(n-1)$-reduced
homotopy system $(\cc_*,f_{n+2})$ $(n\geq 3)$ is the cohomology
class
$$\hat{\cup}_{(\cc_*,f_{n+2})}\in H^{n+2}(\cc_*,\hat{\otimes}^2_TH_n\cc_*)$$
represented by the cocycle
$$\cc_{n+2}\st{f_{n+2}}\To\Pi_{n+1}C_{d_{n+1}}\st{\gamma_2^{n-1}}\To\hat{\otimes}^2_T
H_n\cc_*.$$ This is a cocycle because $f_{n+2}$ is a cocycle. Notice
that in order to define the reduced cup-product we have to choose an
$(n-1)$-fold suspension structure in the cofiber $C_{d_{n+1}}$, this
is possible by Freudental's suspension theorem, otherwise the
reduced James-Hopf invariant $\gamma_2^{n-1}$ would not be defined.

The reduced cup-product invariant of a $1$-reduced homotopy system
$(\cc_*,f_4)$ is the following change of coefficients of the
cup-product invariant
$$\hat{\cup}_{(\cc_*,f_4)}=(\bar{\sigma}_T)_*\cup_{(\cc_*,f_4)}\in
H^4(\cc_*,\hat{\otimes}^2_TH_2\cc_*).$$
\end{defn}

One can check that the reduced cup-product is a $0$-cocycle in
$\C{H}_n(T)/\!\simeq$, for $n=2$ it follows directly from
Proposition \ref{cupnat} and for $n\geq 3$ one argues as in the
proof of that proposition but using Proposition \ref{composicio2}
instead of Proposition \ref{composicio1}.

\begin{prop}\label{cuprednat}
Given a morphism $(\xi,\eta)\colon (\cc_*,f_{n+2})\r
(\cc_*',g_{n+2})$ in the category $\C{H}_n(T)/\!\simeq$ $(n\geq 2)$ the
following equality holds
$$\xi_*\hat{\cup}_{(\cc_*,f_{n+2})}=\xi^*\hat{\cup}_{(\cc_*',g_{n+2})}\in H^{n+2}(\cc_*,\hat{\otimes}^2_TH_n\cc_*').$$
Here $\xi^*=H^{n+2}(\xi,\hat{\otimes}^2_TH_n\cc_*')$ and
$\xi_*=H^{n+2}(\cc_*,\hat{\otimes}^2_TH_n\xi)$.
\end{prop}

\begin{cor}
The reduced cup-product invariant of an $(n-1)$-reduced homotopy
system $(\cc_*,f_{n+2})$ $(n\geq 2)$ does not depend on the
$(n-1)$-fold suspension structure chosen in $C_{d_{n+1}}$ for its
definition.
\end{cor}

By Definition \ref{JHS1} the reduced cup-product invariant is
preserved by suspensions in the sense of the statement of the
following proposition.

\begin{prop}\label{redcupst}
The following equality holds for any $(n-1)$-reduced homotopy
system $(\cc_*,f_{n+2})$ $(n\geq 2)$
$$\hat{\cup}_{(\cc_*,f_{n+2})}=\hat{\cup}_{\Sigma (\cc_*,f_{n+2})}\in
H^{n+2}(\cc_*,\hat{\otimes}^2_TH_n\cc_*)=H^{n+3}(\Sigma\cc_*,\hat{\otimes}^2_TH_{n+1}\Sigma\cc_*).$$
\end{prop}

We can now can define cup-product invariants for $T$-complexes by
using the functors $r$ in (\ref{funtoR}).

\begin{defn}
The \emph{cup-product invariant} of a $1$-reduced $T$-complex $X$
is
$$\cup_X=\cup_{rX}\in H^4(X,\otimes^2_T\hh_2X).$$
Moreover, the \emph{reduced cup-product invariant} of an
$(n-1)$-reduced $T$-complex $X$  $(n\geq 2)$ is
$$\hat{\cup}_X=\hat{\cup}_{rX}\in H^{n+2}(X,\hat{\otimes}^2_T\hh_nX).$$
\end{defn}

These invariants are $0$-cocycles in $\C{CW}^T_n/\!\simeq$ by
Propositions \ref{cupnat} and \ref{cuprednat}, as we state in the
following proposition.

\begin{prop}
Given a proper map $f\colon X\r Y$ between $1$-reduced
$T$-complexes the following equality holds
$$f_*\cup_X=f^*\cup_Y\in H^4(X,\otimes^2_T\hh_2Y).$$
Here $f^*=H^4(f,\otimes^2_T\hh_2Y)$ and
$f_*=H^4(X,\otimes^2_T\hh_2f)$.
\end{prop}

\begin{prop}\label{kaka}
Given a proper map $f\colon X\r Y$ between $(n-1)$-reduced
$T$-complexes $(n\geq 2)$ the following equality holds
$$f_*\hat{\cup}_X=f^*\hat{\cup}_Y\in H^{n+2}(X,\hat{\otimes}^2_T\hh_nY).$$
Here $f^*=H^{n+2}(f,\hat{\otimes}^2_T\hh_nY)$ and
$f_*=H^{n+2}(X,\hat{\otimes}^2_T\hh_nf)$.
\end{prop}

\section{Cup-product and obstruction theory}

The following theorem shows the significance of the cup-product invariant in obstruction theory. In the statement we use the natural transformations $\tau_T$ and $\bar{\tau}_T$ in (\ref{d2}).

\begin{thm}\label{cupobs}
Given two $(n-1)$-reduced homotopy systems $(\cc_*,f_{n+2})$ and
$(\cc_*',g_{n+2})$ and a chain homotopy class
$\xi\colon\cc_*\r\cc_*'$ the following equalities hold: if $n=2$
$$(\tau_T)_*\theta_{(\cc_*,f_4),(\cc_*',g_4)}(\xi)=\xi_*\cup_{(\cc_*,f_4)}-\xi^*\cup_{(\cc_*',g_4)}\in H^4(\cc_*,\otimes^2_TH_2\cc_*'),$$
and if $n\geq 3$
$$(\bar{\tau}_T)_*\theta_{(\cc_*,f_{n+2}),(\cc_*',g_{n+2})}(\xi)=\xi_*\hat{\cup}_{(\cc_*,f_{n+2})}-
\xi^*\hat{\cup}_{(\cc_*',g_{n+2})}\in
H^{n+2}(\cc_*,\hat{\otimes}^2_TH_n\cc_*').$$
\end{thm}

Here $\theta$ is the obstruction operator in one of the exact
sequences for functors appearing in (\ref{thetacoh}).

\begin{proof}[Proof of Theorem \ref{cupobs}]
Let us prove the case $n=2$. With the same notation as in Section
\ref{hs} recall that the cohomology class
$\theta_{(\cc_*,f_4),(\cc'_*,g_4)}(\xi)$ is represented by a cocycle
$\beta$ fitting into the following commutative diagram
$$\xymatrix{\cc_4\ar[rr]^{-g_4\xi_4+\eta f_4}\ar[rd]_\beta&&\Pi_3C_{d_3'}\\
&\Gamma_TH_2\cc_*'\ar@{^{(}->}[ru]_{i_3}&}$$ where $\eta\colon
C_{d_3}\r C_{d_3'}$ is a map inducing $\xi$ in dimensions $2$ and
$3$. By Proposition \ref{conmuta1} the change of coefficients
$(\tau_T)_*\theta_{(\cc_*,f_4),(\cc'_*,g_4)}(\xi)$ is represented by
the cocycle
\begin{eqnarray*}
  \gamma_2(-g_4\xi_4+\eta f_4) &=&   -\gamma_2(g_4\xi_4)+\gamma_2(\eta f_4) \\
   &=&
   -\xi_4^*\gamma_2(g_4)+f_4^*\gamma_2(\eta)+\xi_*\gamma_2(f_4)\\
   &=&
   -\gamma_2g_4\xi_4+f_4^*\gamma_2(\eta)+\xi_*\gamma_2(f_4).
\end{eqnarray*}
For the second equality we use Proposition \ref{composicio1} and the
fact that $\gamma_2(\xi_4)=0$, see the proof of Proposition
\ref{cupnat}. Moreover, one can check, also as in the proof of
Proposition \ref{cupnat}, that $f_4^*\gamma_2(\eta)$ is a coboundary
and the case $n=2$ follows since $\gamma_2g_4\xi_4$ and
$\xi_*\gamma_2(f_4)$ are cocycles representing
$\xi^*\cup_{(\cc_*',g_4)}$ and $\xi_*\cup_{(\cc_*,f_4)}$,
respectively.

The case $n\geq 3$ can be obtained analogously by using Propositions
\ref{conmuta2} and \ref{composicio2} instead of Propositions
\ref{conmuta1} and \ref{composicio1}.
\end{proof}

The following corollary is a vanishing result on the change of
coefficients of the cohomology classes $\set{\theta}$ in
(\ref{lasclases}). It is a consequence of Theorem \ref{cupobs} and
Propositions \ref{cupnat} and \ref{cuprednat}, see Remark
\ref{coh1}.

\begin{cor}\label{seanula}
For $n=2$
$$0=(\tau_T)_*\set{\theta}\in H^1(\chain_2(\ab(T))/\!\simeq\;,H^4(-,\otimes^2_TH_2)),$$
and for any $n\geq 3$
$$0=(\bar{\tau}_T)_*\set{\theta}\in H^1(\chain_n(\ab(T))/\!\simeq\;,H^{n+2}(-,\hat{\otimes}^2_TH_n)).$$
\end{cor}

In the second corollary of Theorem \ref{cupobs} we use the notation
introduced in Definition \ref{esf} (5).

\begin{cor}\label{cupaccion}
Let $(\cc_*,f_{n+2})$ be an $(n-1)$-reduced homotopy system, for
$n=2$ if $\alpha\in H^4(\cc_*,\Gamma_TH_2\cc_*)$ then
$$\cup_{(\cc_*,f_4)+\alpha}=\cup_{(\cc_*,f_4)}+(\tau_T)_*(\alpha),$$
and for any $n\geq 3$ if $\alpha\in
H^{n+2}(\cc_*,H_n\cc_*\otimes\Z/2)$ then
$$\hat{\cup}_{(\cc_*,f_{n+2})+\alpha}=\hat{\cup}_{(\cc_*,f_{n+2})}+(\bar{\tau}_T)_*(\alpha).$$
\end{cor}

\section{The chain cup-product and cohomology of
categories}\label{chaincup}

This section contains the definition and elementary properties of the chain cup-produc invariant. We carry out a thorough study of this invariant in the next chapter. This invariant plays a crucial role in the main results obtained in this paper concerning the proper homotopy classification problem for locally compact $A^2_n$-polyhedra.

\begin{defn}\label{cupcc}
The \emph{chain cup-product invariant} of a bounded chain complex
$\cc_*$ in $\ab(T)$ concentrated in degrees $\geq n$ $(n\geq 2)$ is
defined as
$$\bar{\cup}_{\cc_*}=(\bar{q}_T)_*\hat{\cup}_{(\cc_*,f_{n+2})}\in
H^{n+2}(\cc_*,\wedge^2_TH_n\cc_*),$$ where $(\cc_*,f_{n+2})$ is an
$(n-1)$-reduced homotopy system with $\lambda
(\cc_*,f_{n+2})=\cc_*$. Such a homotopy system exists by Proposition
\ref{so}.
\end{defn}

\begin{prop}\label{clasecupcc}
The chain cup-product is a $0$-cocycle of
$\chain_n(\ab(T))/\!\simeq$ with coefficients in
$H^{n+2}(-,H_n\otimes\Z_2)$ $(n\geq 2)$, i. e. given a chain
homotopy class $\xi\colon\cc_*\r\cc_*'$ between bounded chain
complexes in $\ab(T)$ concentrated in degrees $\geq n$ the following
equality holds
$$\xi_*\bar{\cup}_{\cc_*}=\xi^*\bar{\cup}_{\cc_*'}\in H^{n+2}(\cc_*,\wedge^2_TH_n\cc_*').$$
\end{prop}

\begin{proof}
Let $(\cc_*,f_{n+2})$ and $(\cc_*',g_{n+2})$ be $(n-1)$-reduced
homotopy systems. By Proposition \ref{cupobs} and (\ref{d2}) we have
for all $n\geq 3$
\begin{eqnarray*}
  0 &=& (\bar{q}_T)_*(\bar{\tau}_T)_*\theta_{(\cc_*,f_{n+2}),(\cc_*',g_{n+2})}(\xi) \\
   &=& (\bar{q}_T)_*(\xi_*\hat{\cup}_{(\cc_*,f_{n+2})}-\xi^*\hat{\cup}_{(\cc_*',g_{n+2})}) \\
   &=& \xi_*\bar{\cup}_{\cc_*}-\xi^*\bar{\cup}_{\cc_*'},
\end{eqnarray*}
and for $n=2$
\begin{eqnarray*}
  0 &=& (\bar{q}_T)_*(\bar{\tau}_T)_*(\sigma_T)_*\theta_{(\cc_*,f_4),(\cc_*',g_4)}(\xi) \\
   &=& (\bar{q}_T)_*(\bar{\sigma}_T)_*(\tau_T)_*\theta_{(\cc_*,f_4),(\cc_*',g_4)}(\xi)  \\
   &=& (\bar{q}_T)_*(\bar{\sigma}_T)_*(\xi_*\cup_{(\cc_*,f_4)}-\xi^*\cup_{(\cc_*',g_4)}) \\
   &=& (\bar{q}_T)_*(\xi_*\hat{\cup}_{(\cc_*,f_4)}-\xi^*\hat{\cup}_{(\cc_*',g_4)}) \\
   &=& \xi_*\bar{\cup}_{\cc_*}-\xi^*\bar{\cup}_{\cc_*'}.
\end{eqnarray*}
\end{proof}

We will denote $$\bar{\cup}_n\in
H^0(\C{chain}_n(\ab(T))/\!\simeq\;,H^{n+2}(-,\wedge^2_TH_n))$$ to
the element in cohomology defined by the chain cup-product $(n\geq
2)$.

\begin{cor}
The chain cup-product invariant of a bounded chain complex $\cc_*$
in $\ab(T)$ does not depend on the $(n-1)$-reduced homotopy system
$(\cc_*,f_{n+2})$ chosen for its definition.
\end{cor}

This can be regarded as a special case of Proposition
\ref{clasecupcc} for the identity morphism in $\cc_*$.

As a consequence of the stability of the reduced cup-product, see
Proposition \ref{redcupst}, we have that the chain cup-product
invariant is also stable.

\begin{prop}\label{cupccestable}
Given a bounded chain complex $\cc_*$ in $\ab(T)$ concentrated in
degrees $\geq n$ $(n\geq 2)$ we have that
$$\bar{\cup}_{\cc_*}=\bar{\cup}_{\Sigma\cc_*}\in
H^{n+2}(\cc_*,\wedge^2_TH_n\cc_*)=H^{n+3}(\Sigma\cc_*,\wedge^2_TH_{n+1}\Sigma\cc_*),$$
in particular
$$\bar{\cup}_n=\Sigma^*\bar{\cup}_{n+1}.$$
\end{prop}

From now one we shall simply write $\bar{\cup}$ for the cohomology
element $\bar{\cup}_n$, or $\bar{\cup}^T$ if we want to specify the
tree.

Given a proper map $f\colon T\r T'$ between trees the ``change of
tree'' functor $\uf{F}^f\colon\ab(T)\r\ab(T')$ in Proposition
\ref{induce} induces an obvious functor
$$\uf{F}^f\colon\C{chain}_n(\ab(T))/\!\simeq\,\To\C{chain}_n(\ab(T'))/\!\simeq,$$
already considered in (\ref{lasma}), and a natural transformation
between bimodules over $\C{chain}_n(\ab(T))/\!\simeq$
$$\bar{\uf{F}}^f\colon H^{n+2}(-,\wedge^2_TH_n)\To H^{n+2}(\uf{F}^f,\wedge^2_TH_n\uf{F}^f)$$
by taking $H^{n+2}$ in the following composition of cochain
homomorphisms
\begin{equation}\label{tecambio3}
\xymatrix{\hom_{\ab(T)}(\cc_*,\wedge^2_TH_n\cc_*)\ar[r]^<(.18){\uf{F}^f_*}&\hom_{\ab(T')}(\uf{F}^f\cc_*,\uf{F}^f_*\wedge^2_TH_n\cc_*)\ar[d]\\
\hom_{\ab(T')}(\uf{F}^f\cc_*,\wedge^2_{T'}H_n\uf{F}^f\cc_*)\ar@{=}[r]&
\hom_{\ab(T')}(\uf{F}^f\cc_*,\wedge^2_{T'}\uf{F}^f_*H_n\cc_*)}
\end{equation}
defined as (\ref{tecambio}) and (\ref{tecambio2}).

\begin{prop}\label{compata}
Given a proper map $f\colon T\r T'$ between trees the chain
cup-product in cohomology of categories is compatible with the
``change of tree'' functor $\uf{F}^f\colon\ab(T)\r\ab(T')$, i. e.
with the notation above
$$(\uf{F}^f)^*\bar{\cup}^T=\bar{\uf{F}}^f_*\bar{\cup}^T.$$
\end{prop}

\begin{proof}
By the stability of the chain cup-product, see Proposition
\ref{cupccestable}, it is enough to make the proof of $n=2$. Recall
that we also have a topological ``change of tree'' functor
$f_*\colon\C{Topp}^T_c/\!\simeq\,\r\C{Topp}^{T'}_c/\!\simeq$, see
(\ref{cbase}). If $Z$ is any $T$ complex this functor induces an
$\ab(T)$-module morphism
$$\Pi_3Z\To(\uf{F}^f)^*\Pi_3f_*Z,$$
given by
$$(\Pi_3Z)(\Z\grupo{A}_\alpha)=[S^3_\alpha,Z]^T\st{f_*}\To[f_*S^3_\alpha,f_*Z]^{T'}=((\uf{F}^f)^*\Pi_3Z)(\Z\grupo{A}_\alpha).$$
This morphism has an adjoint which is an $\ab(T')$-module morphism
$$\uf{F}^f_*\Pi_3Z\To\Pi_3f_*Z.$$ The functor $f_*$ commutes up to natural equivalence with the suspension,
and it is easy to check that for $Z=\Sigma Y$ this adjoint
$\ab(T')$-module morphism above is compatible with the James-Hopf
invariant in Definition \ref{JH1}, i. e. there is a commutative
diagram
$$\xymatrix{\uf{F}^f_*\Pi_3\Sigma Y\ar[rr]^{\uf{F}^f_*\gamma_2}\ar[d]&&\uf{F}^f_*\otimes^2_T\hh_2\Sigma Y\ar[d]\\
\Pi_3\Sigma f_*Y\ar[r]^{\gamma_2}&\otimes^2_{T'}\hh_2\Sigma f_*
Y\ar@{=}[r]&\otimes^2_{T'}\uf{F}^f_*\hh_2\Sigma Y}$$ The vertical
arrow in the right is given by Proposition \ref{quadpropri2} and
$\uf{F}^f_*\hh_2\Sigma Y=\hh_2\Sigma f_* Y$ because $\cc_*\Sigma f_*
Y=\uf{F}^f\cc_*\Sigma Y$ is concentrated in dimensions $\geq 2$ and
$\uf{F}^f_*$ is right-exact.

This implies that if $(\cc_*,f_4)$ is a $1$-reduced $T$-homotopy
system then $(\uf{F}^f\cc_*,f_*f_4)$ is also a $1$-reduced
$T'$-homotopy system and
$$\bar{\uf{F}}^f\cup_{(\cc_*,f_4)}=\cup_{(\uf{F}^f\cc_*,f_*f_4)}\in H^4(\uf{F}^f\cc_*,\otimes^2_TH_2\uf{F}^f\cc_*),$$
in particular
$$\bar{\uf{F}}^f\bar{\cup}_{\cc_*}=\bar{\cup}_{\uf{F}^f\cc_*},$$
as we wanted to show.
\end{proof}

\begin{rem}\label{compata2}
Remark \ref{concola2} also applies in this case, i. e. if
$\F(f)\colon\F(T)\hookrightarrow\F(T')$ is injective then
$\bar{\uf{F}}^f$ is a natura isomorphism because (\ref{tecambio3})
is a cochain isomorphism. In particular if $\bar{\cup}^T\neq 0$ for
a tree $T$ with a finite number of ends then $\bar{\cup}^{T'}\neq 0$
provided $T'$ has more ends that $T$, compare Corollary \ref{nohay}.
\end{rem}

\section{Cup-product and co-H-multiplications}

In this section we characterize under some assumptions the existence of co-H-multiplications by using the chain cup-product. As a consequence of this and the computations in the next chapter we obtain in Corollary \ref{kio} examples of degree $2$ proper Moore spaces which do not admit co-H-multiplications.

\begin{prop}\label{cohcup}
A $1$-connected $T$-reduced $X$ with $\dim X\leq 4$ possesses a
co-H-multiplication if and only if $0=\cup_X\in
H^4(X,\otimes^2_T\hh_2X)$.
\end{prop}

\begin{proof}
Let $f_n\colon S^{n-1}_{\alpha_n}\r X^{n-1}$ be the attaching map of
$n$-cells on $X$ $(n\geq 2)$, and $\mu\colon X^3\r X^3\vee X^3$ the
co-H-multiplication associated to the suspension structure chosen
over $X^3$ for the definition of the James-Hopf invariant
$\gamma_2\colon\Pi_3X^3\r\otimes^2_T\hh_2X$. Suppose that
$\cup_X=0$, then there exists $\alpha\in(\otimes^2_T\hh_2X)(\cc_3X)$
such that $\gamma_2(f_4)=\alpha d_4$. If we write $\tilde{\alpha}$
for the image of $\alpha$ under the natural projection
$$(\otimes^2_T\hh_2X)(\cc_3X)\twoheadrightarrow
H^3(X^3,\otimes^2_T\hh_2X),$$ the equality $\gamma_2(f_4)=\alpha
d_4$ is equivalent to say that $\gamma_2(f_4)=f_4^*(\tilde{\alpha})$
since $\alpha d_4=f_4^*\tilde{\alpha}$, compare with the proof of
Proposition \ref{cupnat}. By Proposition \ref{portutatis}
$-j(i_3i_{12})_*(\tilde{\alpha})+\mu$ is also a co-H-multiplication
in $X^3$. In addition if $\mu'$ is the co-H-multiplication of
$S^3_{\alpha_4}$, which is unique up to homotopy, then
\begin{eqnarray*}
(-j(i_3i_{12})_*(\tilde{\alpha})+\mu)f_4 & = &
-j(i_3i_{12})_*f_4^*(\tilde{\alpha})+\mu f_4 \\ &=&
-j(i_3i_{12})_*\gamma_2(f_4)+\mu f_4\\
   &=& i_1f_4+i_2f_4-\mu f_4+\mu f_4 \\
   &=& i_1f_4+i_2f_4 \\
   &=& (f_4\vee f_4)\mu'.
\end{eqnarray*}
Here we use in the first equality the naturality of the central
extension in Proposition \ref{portutatis}, and in the third one
Definition \ref{JH1}. The previous equalities show that the
following diagram commutes up to homotopy,
$$\xymatrix@C=70pt{X^3\ar[r]^<(.35){-j(i_3i_{12})_*(\tilde{\alpha})+\mu} & X^3\vee X^3\\
S^3_{\alpha_4}\ar[u]^{f_4}\ar[r]_{\mu'}&S^3_{\alpha_4}\vee
S^3_{\alpha_4}\ar[u]_{f_4\vee f_4}}$$ This diagram induces a map
(a principal map in the sense of \cite{ah} V.2) between the
cofibers of the vertical arrows which is a co-H-multiplication in
$X$.

Suppose now that $X$ has a co-H-multiplication $\bar{\mu}\colon X\r
X\vee X$. In the following chains of equalities we will use the
inclusions $i_1,i_2\colon X\r X\vee X$ and the projections $p_1,
p_2\colon X\vee X\r X$ of the factors of the coproduct; the
inclusions $i_1,i_2\colon\hh_2X\r\hh_2X\oplus \hh_2X$ and the
projections $p_1,p_2\colon\hh_2X\oplus \hh_2X\r\hh_2X$ of the
factors of the direct sum; the natural identifications $\hh_2(X\vee
X)=\hh_2X\oplus \hh_2X$, $\hh_2i_1=i_1$, $\hh_2i_2=i_2$,
$\hh_2\bar{\mu}=i_1+i_2$, $\hh_2p_1=p_1$ and $\hh_2p_1=p_2$; and
Proposition \ref{cupnat}. We have that
\begin{eqnarray*}
\cup_{X\vee X}&=&p_1^*{i_1}_*\cup_X+p_2^*{i_2}_*\cup_X \\
   &=& p_1^*(\otimes^2_Ti_1)_*\cup_X +
   p_2^*(\otimes^2_Ti_2)_*\cup_X,
\end{eqnarray*}
Moreover, by using these equalities
\begin{eqnarray*}
  (\otimes^2_T(i_1+i_2))_*\cup_X &=& \mu_*\cup_X \\
   &=& \mu^*\cup_{X\vee X} \\
   &=& (\otimes^2_Ti_1)_*\cup_X +
   (\otimes^2_Ti_2)_*\cup_X,
\end{eqnarray*}
i. e.
$$0=(i_1\otimes_Ti_2)_*\cup_X+(i_2\otimes_Ti_1)_*\cup_X\in
H^4(X,\otimes^2_T(\hh_2X\oplus \hh_2X)),$$ therefore
\begin{eqnarray*}
  0 &=& (p_1\otimes_T
p_2)_*((i_1\otimes_Ti_2)_*\cup_X+(i_2\otimes_Ti_1)_*\cup_X) \\
   &=& (p_1i_1\otimes_Tp_2i_2)_*\cup_X+(p_1i_2\otimes_Tp_2i_1)_*\cup_X \\
   &=& \cup_X.
\end{eqnarray*}
Now the proof is finished.
\end{proof}

\begin{prop}\label{cohcupcc}
Let $\cc_*$ be a chain complex in $\ab(T)$ concentrated in degrees
$2$, $3$ and $4$. There exists a $1$-reduced co-H-space $X$ with
$\dim X\leq 4$ and a proper cellular chain complex $\cc_*X=\cc_*$ if
and only if the chain cup-product of $\cc_*$ vanishes
$0=\bar{\cup}_{\cc_*}\in H^4(\cc_*,\wedge^2_TH_2\cc_*)$.
\end{prop}

\begin{proof}
Since $\cc_*$ is concentrated in dimensions $\leq 4$ the following
sequence is exact, see (\ref{d2}),
$$H^4(\cc_*,\Gamma_TH_2\cc_*)\st{(\tau_T)_*}\To H^4(\cc_*,\otimes^2_TH_2\cc_*)\st{(q_T)_*}\twoheadrightarrow
H^4(\cc_*,\wedge^2_TH_2\cc_*).$$ Let $(\cc_*,f_4)$ be a
$1$-reduced homotopy system. By Definition \ref{cupcc}
$$(q_T)_*\cup_{(\cc_*,f_4)}=(\bar{q}_T)_*(\bar{\sigma}_T)_*\cup_{(\cc_*,f_4)}=(\bar{q}_T)_*\hat{\cup}_{(\cc_*,f_4)}=\bar{\cup}_{\cc_*},$$
therefore $\bar{\cup}_{\cc_*}=0$ if and only if there exists
$\alpha\in H^4(\cc_*,\Gamma_TH_2\cc_*)$ such that
$(\tau_T)_*(\alpha)=\cup_{(\cc_*,f_4)}$. If there is such an
$\alpha$, by using Corollary \ref{cupaccion} we see that
$\cup_{(\cc_*,f_4)-\alpha}=0$, so if $X$ is a $T$-complex with
$rX=(\cc_*,f_4)-\alpha=(\cc_*,f_4')$, for example $X=C_{f_4'}$ the
cofiber of the map $S^3_{\alpha_4}\r C_{d_2}$ determined by $f_4'$,
then by Proposition \ref{cohcup} this $X$ is a co-H-space, and
$\cc_*X=\lambda rX=\lambda(\cc_*,f_4)=\cc_*$.

On the other hand if $X$ is a co-H-space with $\cc_*X=\cc_*$ then
$\cup_X=0$ by Proposition \ref{cohcup}, hence
$$\bar{\cup}_{\cc_*}=(\bar{q}_T)_*\hat{\cup}_{rX}=(\bar{q}_T)_*(\bar{\sigma}_T)_*\cup_X=0.$$
\end{proof}


\chapter{The computation of the chain cup-product}\label{chVII}

In the first section of this chapter we give a purely algebraic formula for the chain cup-product invariant which was already defined in a homotopical way in Chapter \ref{chVI}. In the second section we compute the chain cup-product mod $2$ by using the algebraic description in the first section, the localization theorem in cohomology of categories (\cite{fcc}), and the algebrain computations in Chapter \ref{chIX}, which are based on the representation theory in controlled algebra developed in \cite{rca}. The computations in this chapter are crucial for the main homotopical results of this paper.

\section{A purely algebraic description of the chain cup-product}

In the statement of the following theorem we use the linear
extension of categories in Proposition \ref{leab}.

\begin{thm}\label{calculocupcc}
Let $\cc_*$ be a chain complex in $\ab(T)$ concentrated in degrees
$\geq n$,
$$\grupo{A_{n+2}}_{\alpha_{n+2}}^\ni\st{\partial_{n+2}}\To\grupo{A_{n+1}}_{\alpha_{n+1}}^\ni\st{\partial_{n+1}}\To\grupo{A_n}_{\alpha_n}^\ni$$
a sequence of morphisms in $\nil(T)$ whose abelianization is
$$\cc_{n+2}\st{d_{n+2}}\To\cc_{n+1}\st{d_{n+1}}\To\cc_n,$$
$\vartheta\colon\cc_{n+2}\r\wedge^2_T\cc_n$ the unique morphism
such that $\partial_{n+1}\partial_{n+2}=0+\vartheta$, and
$\tilde{p}\colon\cc_n\twoheadrightarrow H_n\cc_*$ the natural
projection. Then the composite
$$\cc_{n+2}\st{\vartheta}\To\wedge^2_T\cc_n\st{\wedge^2_T\tilde{p}}\twoheadrightarrow\wedge^2_TH_n\cc_*$$
is a cocycle representing chain cup-product
$$\bar{\cup}_{\cc_*}\in H^{n+2}(\cc_*,\wedge^2_TH_n\cc_*).$$
\end{thm}

In the proof of this theorem we use an extended technique in
proper homotopy theory which however we only use here. It consists
of lifting natural properties in ordinary homotopy theory to
proper homotopy theory by using the naturality with respect to the
inclusions of a basis of neighbourhoods of the space of
Freudenthal ends into the Freudenthal compactification. We also
use without explicit mention the fact that the ordinary homotopy
theory of pointed spaces coincides with the ordinary homotopy
theory of spaces with a base tree, since trees are all
contractible from the ordinary homotopy viewpoint.

\begin{proof}[Proof of Theorem \ref{calculocupcc}]

By the stability of the chain cup-product, see Proposition
\ref{cupccestable}, it is enough to consider the case $n=2$. Suppose
first that the result is already known for a collection of
length-two sequences in $\nil(T)$ whose abelianizations are
resolutions of all f. p. $\ab(T)$-modules. Let
\begin{equation*}\tag{a}
\grupo{B_{4}}^\ni_{\beta_{4}}\st{\partial_{4}'}\To\grupo{B_{3}}^\ni_{\beta_{3}}\st{\partial_{3}'}\To
\grupo{B_{2}}^\ni_{\beta_{2}}
\end{equation*} be such a sequence
whose abelianization is a resolution of $H_2\cc_*$. Then there is a
diagram in $\nil(T)$
$$\xymatrix{\grupo{A_{4}}_{\alpha_{4}}^\ni\ar[r]^{\partial_{4}}\ar[d]_{f_{4}}&
\grupo{A_{3}}_{\alpha_{3}}^\ni\ar[r]^{\partial_{3}}\ar[d]_{f_{3}}&\grupo{A_2}_{\alpha_2}^\ni\ar[d]_{f_2}\\
\grupo{B_{4}}^\ni_{\beta_{4}}\ar[r]_{\partial_{4}'}&\grupo{B_{3}}^\ni_{\beta_{3}}\ar[r]_{\partial_{3}'}&
\grupo{B_{2}}^\ni_{\beta_{2}}}$$ whose abelianization if a chain
morphism inducing the identity on $H_2$. There are well-defined
controlled homomorphisms
$\vartheta_2\colon\Z\grupo{A_{4}}_{\alpha_4}\r\wedge^2_T\Z\grupo{B_3}_{\beta_3}$,
$\vartheta_1\colon\Z\grupo{A_{3}}_{\alpha_3}\r\wedge^2_T\Z\grupo{B_2}_{\beta_2}$
and
$\vartheta'\colon\Z\grupo{B_4}_{\beta_4}\r\wedge^2_T\Z\grupo{B_2}_{\beta_2}$
such that
\begin{eqnarray*}
\partial'_4f_4+\vartheta_2&=&f_3\partial_4,\\
\partial'_3f_3+\vartheta_1&=&f_2\partial_3,\\
\partial_3'\partial_4'&=&0+\vartheta'.
\end{eqnarray*}
Moreover,
\begin{eqnarray*}
0+(f_n^\abb)_*\vartheta&=&f_2\partial_3\partial_4\\
&=&\partial_3'f_3\partial_4+d_4^*\vartheta_1\\
&=&\partial_3'\partial_4'f_4+d_4^*\vartheta_1+(\partial_3')^\abb_*\vartheta_2\\
&=&0+(f_4^\abb)^*\vartheta'+d_4^*\vartheta_1+(\partial_3')^\abb_*\vartheta_2
\end{eqnarray*}
so
$$(\wedge^2_Tf^\abb_2)\vartheta=\vartheta'f_4^\abb+\vartheta_1d_4+(\wedge^2_T(\partial_3')^\abb)\vartheta_2.$$
If $\tilde{p}'\colon\Z\grupo{B_2}_{\beta_2}\twoheadrightarrow
H_2\cc_*$ is the natural projection then
$\tilde{p}=\tilde{p}'f_2^\abb$, therefore
$$(\wedge^2_T\tilde{p})\vartheta=(\wedge^2_T\tilde{p}')\vartheta'f_4^\abb+(\wedge^2_T\tilde{p}')\vartheta_1d_4,$$
i. e. $(\wedge^2_T\tilde{p})\vartheta$ represents the same
cohomology class as $(\wedge^2_T\tilde{p}')\vartheta'f_4^\abb$ in the cohomology group
$H^4(\cc_*,\wedge^2_TH_2\cc_*)$, but by hypothesis
$(\wedge^2_T\tilde{p}')\vartheta'$ represents the chain cup of the
abelianization of (a), hence by Proposition \ref{kaka}
$(\wedge^2_T\tilde{p}')\vartheta'f_4^\abb$ represents the chain
cup-product of $\cc_*$, and the theorem would follow.

Now we are going to prove the the theorem for a particular
\begin{equation*}\tag{b}
\grupo{A_{4}}_{\alpha_{4}}^\ni\st{\partial_{4}}\To\grupo{A_{3}}_{\alpha_{3}}^\ni\st{\partial_{3}}\To\grupo{A_2}_{\alpha_2}^\ni
\end{equation*}
whose abelianization $\cc_*$ is a resolution of an arbitrary f. p.
$\ab(T)$-module $\um{M}=H_2\cc_*$.

The construction in \cite{projdim} of the kernel of a morphism
between free $T$-controlled $\Z$-modules shows that we can choose
the resolution $\cc_*$ in such a way that there are countable bases
of connected neighborhoods $\set{T_{v^i_n}\sqcup
T_{v^i_n}^\F}_{n\geq 0}$ as in Remark \ref{nice2} $(i=,2,3,4)$ of
the points of $\F(T)$ in the Freudenthal compactification $\hat{T}$
such that for all $n\geq 0$ $T_{v_n^4}\subset T_{v_n^3}\subset
T_{v_n^2}$ and the restriction of $\cc_*$
\begin{equation*}\tag{c}
\Z\grupo{\alpha^{-1}_4(T_{v_n^4})}\st{d_4^n}\hookrightarrow\Z\grupo{\alpha_3^{-1}(T_{v_n^3})}\st{d_3^n}\To
\Z\grupo{\alpha_2^{-1}(T_{v_n^2})}
\end{equation*}
is defined an it is an exact sequence of abelian groups.

The resolution $\cc_*$ is the abelianization of a sequence in
$\gr(T)$
$$\grupo{A_{4}}_{\alpha_{4}}\st{\bar{\partial}_{4}}\To\grupo{A_{3}}_{\alpha_{3}}\st{\bar{\partial}_{3}}\To\grupo{A_2}_{\alpha_2}$$
such that the restrictions
\begin{equation*}\tag{d}
\grupo{\alpha^{-1}_4(T_{v_n^4})}\st{\bar{\partial}_4^n}\hookrightarrow
\grupo{\alpha_3^{-1}(T_{v_n^3})}\st{\bar{\partial}_3^n}\To
\grupo{\alpha_2^{-1}(T_{v_n^2})}
\end{equation*}
are also defined for all $n\geq 0$. We define (b) as
$\partial_i=\partial_i^\ni$ $(i=3,4)$.

In order to compute the chain cup-product $\bar{\cup}_{\cc_*}$ we
are going to choose a convenient $1$-reduced $T$-homotopy system
$(\cc_*,f_4)$, see Definition \ref{cupcc}. Consider the $T$-complex
$Y=C_{\bar{\partial}_3}$ which is the cofiber of the map
$S^1_{\alpha_3}\r S^1_{\alpha_2}$ induced by $\bar{\partial}_3$, see
Proposition \ref{fund2}. Notice that $\Sigma Y= C_{d_3}$. Let
$\alpha_i^n\colon\alpha_i^{-1}(T_{v^i_n})\r T_{v^2_n}$ be the
restriction of $\alpha_i$ $(i=2,3,4)$, and
$Y_n=C_{\bar{\partial}^n_3}$ the cofiber of the map
$S^1_{\alpha^n_3}\r S^1_{\alpha^n_2}$ induced by the controlled
homomorphism $\bar{\partial}^n_3$ in (d).
\begin{equation*}\tag{e}
\begin{array}{l}
\text{The collection $\set{\Sigma Y_n\sqcup T^\F_{v^2_n}}_{n\geq0}$
is a basis of connected neighborhoods of}\\\text{$\F(T)=\F(\Sigma
Y)$ in the Freudenthal compactification $\widehat{\Sigma Y}$,}
\end{array}
\end{equation*}
i. e. any point of $\F(T)$ has an ``arbitrary small'' neighborhood
in $\widehat{\Sigma Y}$ belonging to this collection, despite
$\Sigma Y_n\sqcup T^\F_{v^2_n}$ need not be a neighborhood of all
points in $T^\F_{v^2_n}$, but only of those $T^\F_{v^3_n}\subset
T^\F_{v^2_n}$.

By the exactness of (c) the ordinary homology group $H_3$ of $\Sigma
Y_n$ is
$$H_3\Sigma Y_n=\Z\grupo{\alpha_4^{-1}(T_{v_n^4})},\;\;\;n\geq 0.$$
As we remarked in \cite{cosqc} 6.2 in the proof of \cite{cosqc} 6.1
we construct sections $\zeta_{Y_n}$ of the ordinary Hurewicz
homomorphism $\pi_3\Sigma Y_n\twoheadrightarrow H_3\Sigma Y_n$ which
are compatible with the inclusion of subcomplexes, therefore by (e)
they determine a proper map $f_4\colon S^3_{\alpha_4}\r\Sigma Y$
which restricts to maps $S^3_{\alpha^n_4}\r\Sigma Y_n$ inducing
$\zeta_{Y_n}$ on $\pi_3$ $(n\geq 0)$. Now it is immediate to notice
that $(\cc_*,f_4)$ is a well-defined $1$-reduced $T$-homotopy
system.

The homomorphism $\eta$ defined at the beginning of the proof of
\cite{cosqc} 6.1 determines a controlled homomorphism
$\eta\colon\wedge^2_T\Z\grupo{A_2}_{\alpha_2}\r\otimes^2_T\Z\grupo{A_2}_{\alpha_2}$.
It satisfies $q_T\eta=1$. If $i_1$ and $i_2$ are the inclusions of
the factors of the coproduct
$\grupo{A_3}_{\alpha_3}^\ni\vee\grupo{A_3}_{\alpha_3}^\ni$ and
$i_{12}$ is the inclusion of the quadratic crossed effect of
$\wedge^2_T$ then there is a unique controlled homomorphism
$\nabla\partial_4\colon\Z\grupo{A_4}_{\alpha_4}\r\otimes^2_T\Z\grupo{A_3}_{\alpha_3}$
with
$$(i_2+i_1)\partial_4=(i_2\partial_4+i_1\partial_4)+i_{12}(\nabla\partial_4).$$
By Proposition \ref{piaditiv} the controlled homomorphism
$$\eta(\vartheta-(\wedge^2_Td_2)q_T(\nabla\partial_4))\colon\Z\grupo{A_4}_{\alpha_4}\r\otimes^2_T\Z\grupo{A_2}_{\alpha_2}$$
can be identified with the proper homotopy class of a map
$\kappa\colon S^3_{\alpha_4}\r S^2_{\alpha_2}\vee S^2_{\alpha_2}$
under $T$ which becomes null-homotopic in we project onto any of the
factors of the target.

Let $i_1, i_2, \mu\colon\Sigma Y\r\Sigma Y\vee\Sigma Y$ be now the
inclusion of the factors and the co-H-multiplication of the
suspension of $Y$. The homotopy class of
$$S^3_{\alpha_4}\st{\kappa}\To S^2_{\alpha_2}\vee S^2_{\alpha_2}\subset\Sigma Y\vee\Sigma Y$$
represents $$\mu f_4-i_2f_4-i_1f_4\in[S^3_{\alpha_4},\Sigma
Y\vee\Sigma Y]^T.$$ This follows from (d) and the fact that the
corresponding formula in ordinary homotopy theory holds for all
restrictions
$$S^3_{\alpha_4^n}\st{\kappa}\To S^2_{\alpha_2^n}\vee S^2_{\alpha_2^n}\subset\Sigma Y_n\vee\Sigma Y_n.$$
In particular by Definition \ref{JH1}
$$\gamma_2(f_4)=(\otimes^2_T\tilde{p})\eta(\vartheta-(\wedge^2_Td_2)q_T(\nabla d_4)),$$
and by Definitions \ref{redcup} and \ref{cupcc} the chain
cup-product $\bar{\cup}_{\cc_*}$ in indeed represented by
\begin{eqnarray*}
\bar{q}_T\bar{\sigma}_T\gamma_2(f_4)&=&q_T\gamma_2(f_4)\\
&=&q_T(\otimes^2_T\tilde{p})\eta(\vartheta-(\wedge^2_Td_2)q_T(\nabla
d_4))\\
&=&(\wedge^2_T\tilde{p})q_T\eta(\vartheta-(\wedge^2_Td_2)q_T(\nabla
d_4))\\
&=&(\wedge^2_T\tilde{p})(\vartheta-(\wedge^2_Td_2)q_T(\nabla d_4))\\
&=&(\otimes^2_T\tilde{p})\vartheta.
\end{eqnarray*}
\end{proof}

In the statement of the following proposition we consider the
$\C{vect}(T_3)$-module $\uf{M}\ul{V}^{(3,5)}$ in Theorem
\ref{clasclas}, see also Proposition \ref{losn}.

\begin{prop}\label{comput1}
If $\cc_*$ is a finite-type resolution of the f. p.
$\C{vect}(T_3)$-module $\uf{M}\ul{V}^{(3,5)}$ as an
$\ab(T_3)$-module then $$0\neq\bar{\cup}_{\cc_*}\in
H^{n+2}(\cc_*,\wedge^2_{T_3}H_n\cc_*)=\Z/2.$$
\end{prop}

\begin{proof}
By Proposition \ref{ZaZ2} and Lemma \ref{03} we have that
\begin{eqnarray*}
H^{n+2}(\cc_*,\wedge^2_{T_3}H_n\cc_*)&=&\ext^2_{\ab(T_3)}(\uf{M}\ul{V}^{(3,5)},\wedge^2_{T_3}\uf{M}\ul{V}^{(3,5)})\\
&=&\ext^1_{\C{vect}(T_3)}(\uf{M}\ul{V}^{(3,5)},\wedge^2_{T_3}\uf{M}\ul{V}^{(3,5)})\\
&=&\Z/2.
\end{eqnarray*}

Let us check now that the cup-product invariant of $\cc_*$ is not
zero. For this we will make a convenient choice of $\cc_*$. In the
proof of \cite{rca} 9.1 we construct a finite presentation of
$\uf{M}\ul{V}^{(3,5)}$ as a $\C{vect}(T_3)$-module
$$\uf{F}_2\grupo{D}_\delta\st{\rho}\hookrightarrow\uf{F}_2\grupo{C}_\gamma\st{\tilde{p}}\twoheadrightarrow\uf{M}\ul{V}^{(3,5)}.$$
Here $D=\set{{}_mw^i\,;\,1\leq i\leq 3,m\geq1}$,
$C=D\sqcup\set{w_1,w_2}$, $\delta({}_mw^i)=\gamma({}_mw^i)=v_m^i$,
$\gamma(w_j)=v_0$ and
\begin{equation*}
\rho({}_mw^i)=\left\{%
\begin{array}{ll}
    {}_mw^i-{}_{m-1}w^i, & \hbox{for all $i=1,2,3$ provided $m>1$;} \\
    {}_1w^1-w_1, & \hbox{if $i=1$ and $m=1$;} \\
    {}_1w^2-w_2-w_1, & \hbox{if $i=2$ and $m=1$;} \\
    {}_1w^3-w_2, & \hbox{if $i=3$ and $m=1$.} \\
\end{array}%
\right.
\end{equation*}
Moreover the projection
$\tilde{p}\colon\uf{F}_2\grupo{C}_\gamma\twoheadrightarrow
\uf{M}\ul{V}^{(3,5)}$ is determined by the vector space
homomorphism
$\tilde{p}_0\colon\uf{F}_2\grupo{C}\twoheadrightarrow\uf{F}_2\grupo{x,
y}$ with
\begin{equation*}
\begin{array}{c}\tilde{p}_0({}_mw^i)=\left\{%
\begin{array}{ll}
    x, & \hbox{$i=1$, $m\geq 1$;} \\
    x+y, & \hbox{$i=2$, $m\geq 1$;} \\
    y, & \hbox{$i=3$, $m\geq 1$;}
\end{array}\right.\\{}\\
\tilde{p}_0(w_1)=x;\\{}\\ \tilde{p}_0(w_2)=y.
\end{array}%
\end{equation*}

If
$\bar{\varphi}\colon\Z\grupo{D}_\delta\hookrightarrow\Z\grupo{C}_\gamma$
is the controlled homomorphism defined by the same formulas as
$\rho$ then
$$\Z\grupo{D}_\delta\st{(2,-\bar{\varphi})}\hookrightarrow\Z\grupo{D}_\delta\oplus\Z\grupo{C}_\gamma\st{(\bar{\varphi},2)}\To\Z\grupo{C}_\gamma$$
is a resolution of $\uf{M}\ul{V}^{(3,5)}$ as an $\ab(T_3)$-module.
We take $\cc_*$ as the translation of this resolution to degree
$n$.

Let
$\varphi\colon\grupo{D}_\delta^{nil}\hookrightarrow\grupo{C}_\gamma^{nil}$
be the controlled homomorphism defined as $\rho$ and
$\bar{\varphi}$. The sequence
\begin{equation*}
\grupo{D}^\ni_\delta\st{\partial_{n+2}}\To\grupo{D}_\delta^\ni\vee\grupo{C}_\gamma^\ni\st{\partial_{n+1}}\To\grupo{C}_\gamma^\ni
\end{equation*}
with $\partial_{n+2}(d)=i_1(d)+i_1(d)-i_2\varphi(d)$ and
$\partial_{n+1}(d,c)=(\varphi(d),c+c)$ $(c\in C, d\in D)$ satisfies
the hypothesis of Theorem \ref{calculocupcc}.

Given  $1\leq i\leq 3$ and $m>1$ we have that
\begin{eqnarray*}
  \partial_{n+1}\partial_{n+2}({}_mw^i) &=& \partial_{n+1}(i_1({}_mw^i)+i_1({}_mw^i)-i_2({}_mw^i-{}_{m-1}w^i)) \\
  &=&
  \partial_{n+1}(i_1({}_mw^i)+i_1({}_mw^i)+i_2({}_{m-1}w^i)-i_2({}_mw^i))\\
   &=& {}_mw^i-{}_{m-1}w^i+{}_mw^i-{}_{m-1}w^i\\&&+{}_{m-1}w^i+{}_{m-1}w^i-{}_mw^i-{}_mw^i \\
   &=&  {}_mw^i+[{}_{m-1}w^i,-{}_mw^i]-{}_mw^i\\
   &=& [{}_mw^i,{}_{m-1}w^i].
\end{eqnarray*}
Moreover
\begin{eqnarray*}
    \partial_{n+1}\partial_{n+2}({}_1w^1) &=& \partial_{n+1}(i_1({}_1w^1)+i_1({}_1w^1)-i_2({}_1w^1-w_1))  \\
  &=&    \partial_{n+1}(i_1({}_1w^1)+i_1({}_1w^1)+i_2(w_1)-i_2({}_1w^1))\\
   &=& {}_1w^1-w_1+{}_1w^1-w_1\\&&+w_1+w_1-{}_1w^1-{}_1w^1 \\
   &=& {}_1w^1+[w_1,-{}_1w^1]-{}_1w^1 \\
   &=& [{}_1w^1,w_1],
\end{eqnarray*}

\begin{eqnarray*}
    \partial_{n+1}\partial_{n+2}({}_1w^2) &=& \partial_{n+1}(i_1({}_1w^2)+i_1({}_1w^2)-i_2({}_1w^2-w_2-w_1))  \\
&=&          \partial_{n+1}(i_1({}_1w^2)+i_1({}_1w^2)+i_2(w_1)+i_2(w_2)-i_2({}_1w^2))\\
   &=& {}_1w^2-w_2-w_1+{}_1w^2-w_2-w_1\\&&+w_1+w_1+w_2+w_2-{}_1w^2-{}_1w^2 \\
   &=& {}_1w^2+[w_1+w_2,w_2-{}_1w^2]-{}_1w^2 \\
   &=& [{}_1w^2,w_1]+[{}_1w^2,w_2]+[w_1,w_2],
\end{eqnarray*}

\begin{eqnarray*}
    \partial_{n+1}\partial_{n+2}({}_1w^3) &=& \partial_{n+1}(i_1({}_1w^3)+i_1({}_1w^3)-i_2({}_1w^3-w_2))  \\
  &=&        \partial_{n+1}(i_1({}_1w^3)+i_1({}_1w^3)+i_2(w_2)-i_2({}_1w^3))\\
   &=& {}_1w^3-w_2+{}_1w^3-w_2\\&&+w_2+w_2-{}_1w^3-{}_1w^3 \\
   &=& {}_1w^3+[w_2,-{}_1w^3]-{}_1w^3 \\
   &=& [{}_1w^3,w_2],
\end{eqnarray*} therefore the controlled homomorphism
$\vartheta\colon\Z\grupo{D}_\delta\r\wedge^2_{T_3}\Z\grupo{C}_\gamma$
satisfying the equation $\partial_{n+1}\partial_{n+2}=0+\vartheta$ is given by
$$\vartheta({}_mw^i)=\left\{%
\begin{array}{ll}
    {}_mw^i\wedge{}_{m-1}w^i, & \hbox{for all $i=1,2,3$ if $m>1$;} \\
    {}_1w^1\wedge w_1, & \hbox{if $i=1$ and $m=1$;} \\
    {}_1w^2\wedge w_1+{}_1w^2\wedge w_2+w_1\wedge w_2, & \hbox{if $i=2$ and $m=1$;} \\
    {}_1w^3\wedge w_2, & \hbox{if $i=3$ and $m=1$.} \\
\end{array}%
\right.    $$

For any $\ab(T_3)$-module $\um{M}$ we write
$\hat{p}\colon\um{M}\twoheadrightarrow\um{M}\otimes\Z/2$ for the
natural projection. By Theorem \ref{calculocupcc}
$\bar{\cup}_{\cc_*}$ is represented by the cocycle
$$\Z\grupo{D}_\delta\st{\vartheta}\To\Z\grupo{C}_\gamma\st{\wedge^2_{T_3} {(\tilde{p}\hat{p})}}\To\wedge^2_{T_3}\uf{M}\ul{V}^{(3,5)}.$$
By Proposition \ref{calculazo} and \cite{rca} 7.7
$\wedge^2_{T_3}\uf{M}\ul{V}^{(3,5)}=\uf{F}_2\grupo{x\wedge y}_\phi$
is a free $T_3$-controlled $\uf{F}_2$-module with only one generator
$x\wedge y$. Under this identification the cocycle $(\wedge^2_{T_3}
{(\tilde{p}\hat{p})})\vartheta\colon\Z\grupo{D}_\delta\r\uf{F}_2\grupo{x\wedge
 y}_\phi$ corresponds to the controlled homomorphism
$(\wedge^2\tilde{p}_0)\vartheta\colon\uf{F}_2\grupo{D}_\delta\r\uf{F}_2\grupo{x\wedge
y}_\phi$. 

For any $m>1$ and $i=1,2,3$
\begin{eqnarray*}
 (\wedge^2\tilde{p}_0)\vartheta({}_mw^i) &=& (\wedge^2\tilde{p}_0)({}_mw^i\wedge{}_{m-1}w^i) \\
   &=& \left\{%
\begin{array}{ll}
    x\wedge x, & \hbox{if $i=1$;} \\
    (x+y)\wedge(x+y), & \hbox{if $i=2$;} \\
    y\wedge y, & \hbox{if $i=3$;} \\
\end{array}%
\right\}     \\
   &=& 0; \\
  (\wedge^2\tilde{p}_0)\vartheta({}_1w^1)  &=& (\wedge^2\tilde{p}_0)({}_1w^1\wedge w_1) \\
   &=& x\wedge x \\
   &=& 0; \\
   (\wedge^2\tilde{p}_0)\vartheta({}_1w^2)   &=& (\wedge^2\tilde{p}_0)({}_1w^2\wedge w_1+{}_1w^2\wedge w_2+w_1\wedge w_2) \\
   &=& (x+y)\wedge x+(x+y)\wedge y +x\wedge y \\
   &=& x\wedge y\\
     (\wedge^2\tilde{p}_0)\vartheta({}_1w^3)  &=& (\wedge^2\tilde{p}_0)({}_1w^3\wedge w_2) \\
   &=& y\wedge y\\
   &=& 0.
\end{eqnarray*}

If $(\wedge^2_{T_3} {(\tilde{p}p)})\vartheta$ represented the
trivial cohomology class then there would be a controlled
homomorphism
$$(\xi_1,\xi_2)\colon\uf{F}_2\grupo{D}_\delta\oplus\uf{F}_2\grupo{C}_\gamma\To\uf{F}_2\grupo{x\wedge
y}_\phi$$ such that
\begin{eqnarray*}
(\wedge^2\tilde{p}_0)\vartheta&=&(\xi_1,\xi_2)((2,-\bar{\varphi})\otimes\Z/2)\\
&=&(\xi_1,\xi_2)(0,\rho)\\
&=&(0,\xi_2\rho).
\end{eqnarray*}
In particular for all $m\geq 1$ and $i=1,2,3$ the following
equalities should hold
\begin{eqnarray*}
  0 &=& (\wedge^2\tilde{p}_0)\vartheta({}_{m+1}w^i) \\
   &=& \xi_2\rho({}_{m+1}w^i) \\
   &=& \xi_2({}_{m+1}w^2-{}_mw^i),
\end{eqnarray*}
i. e.
$$\xi_2({}_mw^i)=\xi_2({}_{m+1}w^i),\;\; m\geq 1, i=1,2,3.$$
Since $\xi_2$ is supposed to be controlled and
$\gamma({}_mw^i)=v^i_m$ then $\xi_2({}_mw^i)$ should vanish for $m$
big enough and therefore for all $m$ by using the previous equality
$$\xi_2({}_mw^i)=0,\;\; m\geq 1, i=1,2,3;$$
hence
\begin{eqnarray*}
  0 &=& (\wedge^2\tilde{p}_0)\vartheta({}_1w^1) \\
   &=& \xi_2\rho({}_1w^1) \\
   &=& \xi_2({}_1w^1-w_1)\\
   &=& \xi_2(w_1),
\end{eqnarray*}
and
\begin{eqnarray*}
  0 &=& (\wedge^2\tilde{p}_0)\vartheta({}_1w^3) \\
   &=& \xi_2\rho({}_1w^3) \\
   &=& \xi_2({}_1w^3-w_2)\\
   &=& \xi_2(w_2),
\end{eqnarray*}
so
\begin{eqnarray*}
  x\wedge y &=& (\wedge^2\tilde{p}_0)\vartheta({}_1w^2) \\
   &=& \xi_2\rho({}_1w^2) \\
   &=& \xi_2({}_1w^2-w_2-w_1)\\
   &=&\xi_2(w_1)+\xi_2(w_2)\\
   &=& 0,
\end{eqnarray*}
and we reach a contradiction derived from supposing that the
cohomology class $\bar{\cup}_{\cc_*}$ represented by the cocycle
$(\wedge^2_{T_3} {(\tilde{p}p)})\vartheta$ was trivial.
\end{proof}

The following corollary follows from this last proposition and
Remark \ref{compata2}.

\begin{cor}\label{nohay}
If $T$ has more than $2$ ends then the chain cup-product in
cohomology of categories does not vanish
$$0\neq\bar{\cup}\in H^0(\C{chain}_n(\ab(T))/\!\simeq\;,H^{n+2}(-,\wedge^2_TH_n)).$$
\end{cor}

We do not know wether the chain cup-product vanishes for trees
with $1$ or $2$ ends, however in the following section we
completely compute its mod $2$ version for trees with less than
four ends.

\section[The computation of the mod $2$ chain cup-product\dots]{The computation of the mod $2$ chain cup-product element in cohomology of
categories}\label{ccm2}

Recall from Section \ref{chaincup} that the chain cup-product can
be regarded as an element
$$\bar{\cup}\in
H^0(\C{chain}_n(\ab(T))/\!\simeq\;,H^{n+2}(-,\wedge^2_TH_n)).$$
This element determines obstructions to the existence of
co-H-spaces with a given proper cellular chain complex, see
Proposition \ref{cohcupcc}.

The best we can do to compute this cohomology group is to apply
the localization theorem in \cite{fcc}. By Theorem \ref{kernel}
the homology functor
$$H_n\colon\C{chain}_n(\ab(T))/\!\simeq\;\To\C{fp}(\ab(T))$$ has a
right adjoint determined by the choice of a finite-type resolution
in degree $n$ for any f. p. $\ab(T)$-module, therefore by
\cite{fcc} 6.5
$$H^0(\C{chain}_n(\ab(T))/\!\simeq\;,H^{n+2}(-,\wedge^2_TH_n))\simeq
H^0(\C{fp}(\ab(T)),\ext^2_{\ab(T)}(-,\wedge^2_T)).$$ The image of
$\bar{\cup}$ under this isomorphism is the universal obstruction
to the existence of co-H-structures in degree $2$ proper Moore
spaces.

As we mentioned in Section \ref{cohcat} the cohomology of
$\C{fp}(\ab(T))$ is strongly related to its representation theory,
i. e. the representation theory of $\Z(\F(T))$. As a consequence of
the results in \cite{rca} 4 we can not expect to obtain satisfactory
classification theorems for f. p. $\ab(T)$-modules because such a
result would imply a classification of all countable abelian groups,
therefore we are unable to compute the previous cohomology group of
$\C{fp}(\ab(T))$. For this reason we will consider the image of
$\bar{\cup}$ under the change of coefficients
$$\xymatrix{H^0(\C{chain}_n(\ab(T))/\!\simeq\;,H^{n+2}(-,\wedge^2_TH_n))\ar[d]^{\hat{p}_*}\\
H^0(\C{chain}_n(\ab(T))/\!\simeq\;,H^{n+2}(-,\wedge^2_TH_n\otimes\Z_2))}$$
induced by the natural projection
$\hat{p}\colon\um{M}\twoheadrightarrow\um{M}\otimes\Z_2$. The
element $\hat{p}_*\bar{\cup}$ will be termed \emph{chain cup-product
mod $2$}. It is the universal obstruction to the existence of
co-H-structures in degree $2$ proper Moore spaces with exponent $2$
in proper homology.

Once again the localization theorem in \cite{fcc} yields an
isomorphism
\begin{small}
\begin{equation*}
H^0(\C{chain}_n(\ab(T))/\!\simeq\;,H^{n+2}(-,\wedge^2_TH_n\otimes\Z/2))\simeq
H^0(\C{fp}(\ab(T)),\ext^2_{\ab(T)}(-,\wedge^2_T\otimes\Z/2)).
\end{equation*}
\end{small}
Moreover, the adjoint ``change of coefficients'' functors in
(\ref{amod2}) together with Proposition \ref{ZaZ2} and the
localization theorem in \cite{fcc} yield another isomorphism
\begin{equation*}
H^0(\C{fp}(\ab(T)),\ext^2_{\ab(T)}(-,\wedge^2_T\otimes\Z/2))\simeq
H^0(\C{fp}(\C{vect}(T)),\ext^1_{\C{vect}(T)}(-,\wedge^2_T)),
\end{equation*}
and the representation theory of $\C{fp}(\C{vect}(T))$ is not as
complicated as that of $\C{fp}(\ab(T))$. More precisely, as we
proved in \cite{rca} $\C{fp}(\C{vect}(T))$ has finite
representation type if (and only if) $T$ has less than four ends
and we have explicit classification theorems for f. p.
$\C{vect}(T)$-modules in these cases, see Section \ref{rrt},
therefore we can expect to be able to compute the previous
cohomology group of $\C{fp}(\C{vect}(T))$ for trees with less than
four ends. Indeed we will do it in the following theorem with the
help of the computations carried out in Sections \ref{homext} and
\ref{quadcomp}.

\begin{thm}\label{vayacalculillo}
If $T$ has $1$ or $2$ ends then
$$H^0(\C{chain}_n(\ab(T))/\!\simeq\;,H^{n+2}(-,\wedge^2_TH_n\otimes\Z_2))=0,$$
so $\hat{p}_*\bar{\cup}=0$. Moreover, if $T$ has $3$ ends
$$0\neq\hat{p}_*\bar{\cup}\in H^0(\C{chain}_n(\ab(T))/\!\simeq\;,H^{n+2}(-,\wedge^2_TH_n\otimes\Z_2))=\Z_2.$$
\end{thm}

\begin{proof}
By Propositions \ref{conZ2} and \ref{comput1} $\bar{\cup}\neq 0$
for $T$ a tree with $3$ ends, therefore in this case the
cohomology group in the statement is non-trivial and all we have
to do now is to compute it for trees with less than four ends. As
we have seen above it is isomorphic to
$$H^0(\C{fp}(\C{vect}(T)),\ext^1_{\C{vect}(T)}(-,\wedge^2_T))$$ by
the localization theorem in cohomology of categories, see
\cite{fcc}.

Recall from Section \ref{cohcat} that any $0$-cocycle $c$ is
additive, i. e. it is completely determined by its values
$$c_\um{M}\in\ext^1_{\C{vect}(T)}(\um{M},\wedge^2_T\um{M})$$ where
$\um{M}$ runs along all elementary f. p. $\C{vect}(T)$-modules,
see Section \ref{rrt}. This fact and Lemmas \ref{01}, \ref{02} and
\ref{03} prove that the cohomology group
$$H^0(\C{fp}(\C{vect}(T)),\ext^1_{\C{vect}(T)}(-,\wedge^2_T))$$ is
trivial if $T$ has $1$ or $2$ ends, while it can be trivial or
$\Z/2$ if $T$ has $3$ ends, but as we have noticed before it is
non-trivial in this case, hence we are done.
\end{proof}

The isomorphism observed above and given by the localization
theorem for cohomology of categories (\cite{fcc} 6.5)
\begin{small}
\begin{equation*}
H^0(\C{fp}(\C{vect}(T)),\ext^1_{\C{vect}(T)}(-,\wedge^2_T))\st{\simeq}\To
H^0(\C{chain}_n(\ab(T))/\!\simeq\;,H^{n+2}(-,\wedge^2_TH_n\otimes\Z_2))
\end{equation*}
\end{small}
is $\kappa_*(H_n\otimes\Z/2)^*$, where $H_n\otimes\Z/2$ is the
homology functor
$$H_n\otimes\Z/2\colon\C{chain}_n(\ab(T))/\!\simeq\;\To\C{fp}(\C{vect}(T))$$
and $\kappa$ is the following composite of natural transformations
between bimodules over $\C{chain}_n(\ab(T))$
\begin{equation*}
\xymatrix@R=15pt{\ext^1_{\C{vect}(T)}(H_n\otimes\Z/2,\wedge^2_TH_n\otimes\Z/2)\ar[d]^\simeq\\
\ext^2_{\ab(T)}(H_n\otimes\Z/2,\wedge^2_TH_n\otimes\Z/2)\ar[d]^{\hat{p}^*}\\
\ext^2_{\ab(T)}(H_n,\wedge^2_TH_n\otimes\Z/2)\ar[d]\\
H^{n+2}(-,\wedge^2_TH_n\otimes\Z_2)}
\end{equation*}
which is induced by Proposition \ref{ZaZ2}, the natural projection
$\hat{p}\colon\um{M}\twoheadrightarrow\um{M}\otimes\Z_2$, and the
universal coefficients spectral sequence for the computation of the
cohomology of a chain complex $\cc_*$ in $\ab(T)$
\begin{equation}\label{ucss}
E^{p,q}_2=\ext^p_{\ab(T)}(H_q\cc_*,\um{M})\Rightarrow
H^{p+q}(\cc_*,\um{M}).
\end{equation}
Therefore Proposition
\ref{comput1}, Lemma \ref{03} and the additivity of $0$-cocycles in
cohomology of categories, see Section \ref{cohcat}, prove more than
Theorem \ref{vayacalculillo} for a tree with three ends.

\begin{thm}\label{calculote}
Given a tree $T$ with $3$ ends and a bounded chain complex $\cc_*$
in $\ab(T)$ concentrated in dimensions $\geq n$, if
$H_n\cc_*\otimes\Z/2$ contains $\uf{M}\ul{V}^{(3,5)}$ as a direct
summand the cup-product mod $2$ $\hat{p}_*\bar{\cup}_{\cc_*}$ is
the image of the generator
\begin{equation*}
\xymatrix@R=15pt{
\Z/2=\ext^1_{\C{vect}(T)}(\uf{M}\ul{V}^{(3,5)},\wedge^2_{T}\uf{M}\ul{V}^{(3,5)})\ar@{^{(}->}[d]\\
\ext^1_{\C{vect}(T)}(H_n\cc_*\otimes\Z/2,\wedge^2_TH_n\cc_*\otimes\Z/2)\ar[d]^{\kappa}\\
H^{n+2}(\cc_*,\wedge^2_TH_n\cc_*\otimes\Z/2)}
\end{equation*}
Otherwise $\hat{p}_*\bar{\cup}_{\cc_*}=0$.
\end{thm}

\chapter{Pontrjagin-Steenrod invariants and proper homotopy types}\label{chVIII}

This is the central chapter of this paper. Here we obtain our main proper homotopy classification results for $A^2_n$-polyhedra when the involved representation theory is of finite type. For this we answer the questions on proper Pontrjagin-Steenrod invariants raised in the introduction. In order to obtain these results we use the computations carried out in Chapters \ref{chVII} and \ref{chIX}.

\section{On the existence of Pontrjagin-Steenrod invariants in proper homotopy
theory}\label{epsi}

The Pontrjagin-Steenrod invariant in ordinary homotopy theory is one
of the tools used by J. H. C. Whitehead to classify homotopy types
of $A^2_n$-polyhedra in terms of algebraic data, see \cite{htskp},
\cite{sc4p} and \cite{aces}. The alternative tool was the lower part
of Whitehead's exact sequence for the Hurewicz homomorphism. In fact
the lower part of this sequence can be derived from the
Pontrjagin-Steenrod invariant and the universal coefficients exact
sequence. As we mentioned in Section \ref{ctpht} there is a proper
analogue of Whitehead's exact sequence, however it can not be used
to classify proper homotopy types of $A^2_n$-polyhedra because there
are proper Moore spaces with the same homology but a different
proper homotopy type, see \cite{tsukuba95} or Remarks \ref{tata} and
\ref{tata2} below. For this reason it becomes more interesting to
study the existence of Pontrjagin-Steenrod invariants in proper
homotopy theory.

A Pontrjagin-Steenrod invariant in proper homotopy theory should be
a cohomology invariant ($0$-cocycle) in the category of
$(n-1)$-reduced $T$-homotopy systems
\begin{equation*}
\wp_n{(\cc_*,f_{n+2})}\in\left\{%
\begin{array}{ll}
    H^4(\cc_*,\Gamma_TH_2\cc_*), & \hbox{$n=2$ (Pontrjagin invariant);} \\
    {} \\
    H^{n+2}(\cc_*,H_n\cc_*\otimes\Z/2), & \hbox{$n>2$ (Steenrod invariant).} \\
\end{array}%
\right.
\end{equation*}
This cohomology invariant should determine the obstruction operator
$\theta$ in the exact sequence of functors (\ref{thetacoh}), i. e.
given a chain morphism $\xi\colon\cc_*\r\cc_*'$ the following
equality should hold
$$\theta_{(\cc_*,f_{n+2}),(\cc_*',g_{n+2})}(\xi)=\xi_*\wp_n{(\cc_*,f_{n+2})}-\xi^*\wp_n{(\cc_*',g_{n+2})}.$$
Therefore Remark \ref{coh1} shows that the existence of this
invariant is equivalent to the vanishing of the characteristic class
in cohomology of categories determined by the obstruction operator,
see (\ref{lasclases}),
\begin{equation*}
\begin{array}{cl}
\set{\theta}\in H^1(\C{chain}_2(\ab(T))/\!\simeq\;,H^4(-,\Gamma_TH_2)), & \hbox{if $n=2$;} \\
    {} \\
\set{\theta}\in H^1(\C{chain}_n(\ab(T))/\!\simeq\;,H^{n+2}(-,H_n\otimes\Z/2)), & \hbox{if $n>2$.} \\
\end{array}%
\end{equation*}

If the proper Pontrjagin-Steenrod invariant were completely
analogous to the ordinary one then it should be compatible with the
(reduced) cup-product invariant defined in Section \ref{cpi}, i. e.
\begin{equation*}
\begin{array}{cl}
(\tau_T)_*\wp_2{(\cc_*,f_{4})}=\cup_{(\cc_*,f_{4})}, &
\text{if } n=2;\\{}\\
(\bar{\tau}_T)_*\wp_n{(\cc_*,f_{n+2})}=\hat{\cup}_{(\cc_*,f_{n+2})},
& \text{if } n>2;
\end{array}
\end{equation*}
see also (\ref{d2}). The existence of a proper Pontrjagin-Steenrod
invariant compatible with the (reduced) cup-product would
immediately imply the vanishing of the chain cup-product, see
(\ref{d2}) and Definition \ref{cupcc}, therefore we can already
discard this possibility for trees with more than two ends, see
Corollary \ref{nohay}. However for our purposes, the classification
of proper homotopy types in terms of algebraic data, this
compatibility is not necessary, therefore we will now concentrate in
the study of $\set{\theta}$.

\section{Steenrod invariants for spaces with less than $4$ ends}

Suppose that $T$ is a tree such that the natural transformation in
(\ref{d2})
\begin{equation}\label{doesplit}
\bar{\tau}_T\colon\um{M}\To\hat{\otimes}^2_T\um{M}
\end{equation} is a split monomorphism for any f. p. $\C{vect}(T)$-module
$\um{M}$. In this case we would have a short exact sequence of
natural transformations between bimodules over
$\C{chain}_n(\ab(T))/\!\simeq$
\begin{equation}\label{shocorto}
H^{n+2}(-,H_n\otimes\Z_2)\st{(\bar{\tau}_T)_*}\hookrightarrow
H^{n+2}(-,\hat{\otimes}^2_TH_n)\st{(\bar{q}_T)_*}\twoheadrightarrow
H^{n+2}(-,\wedge^2_TH_n),
\end{equation}
which would give rise to a Bockstein long exact sequence in
cohomology of categories
\begin{equation}\label{bocks}
\xymatrix{\vdots\ar[d]\\
H^k(\C{chain}_n(\ab(T))/\!\simeq\;,H^{n+2}(-,H_n\otimes\Z_2))\ar[d]^{(\bar{\tau}_T)_*}\\
H^k(\C{chain}_n(\ab(T))/\!\simeq\;,H^{n+2}(-,\hat{\otimes}^2_TH_n\otimes\Z_2))\ar[d]^{(\bar{q}_T)_*}\\
H^k(\C{chain}_n(\ab(T))/\!\simeq\;,H^{n+2}(-,\wedge^2_TH_n\otimes\Z_2))\ar[d]^{\beta}\\
H^{k+1}(\C{chain}_n(\ab(T))/\!\simeq\;,H^{n+2}(-,H_n\otimes\Z_2))\ar[d]\\\vdots}
\end{equation}
Otherwise the sequence of bimodules would only be exact in the
middle and the Bockstein exact sequence would not even be defined.

\begin{thm}\label{image}
In the previous conditions the characteristic class $\set{\theta}$
represented by the obstruction operator $\theta$ in (\ref{thetacoh})
$(n> 2)$ is the image of the chain cup-product mod $2$
$\hat{p}_*\bar{\cup}$ in Section \ref{ccm2} under the Bockstein
homomorphism
$$\set{\theta}=\beta\hat{p}_*\bar{\cup}\in H^1(\C{chain}_n(\ab(T))/\!\simeq\;,H^{n+2}(-,H_n\otimes\Z_2)).$$
In particular this is satisfied for $T$ a tree with less than four
ends.
\end{thm}

The first part of the statement is an immediate consequence of
Definition \ref{cupcc} and Theorem \ref{cupobs}. The second part
follows from Proposition \ref{splitmod2}.

The weaker fact that $\set{\theta}$ is in the image of the
Bockstein homomorphism can also be deduced from Corollary
\ref{seanula}.

Moreover, by Remark \ref{herculillo} Theorem \ref{image} can not
hold for trees with more than six ends because, as one can readily
see, the splitability of (\ref{doesplit}) is equivalent to the
exactness of (\ref{shocorto}). For trees with $4$, $5$ or $6$ ends
nothing is known.

The following corollary follows from Theorems \ref{image} and
\ref{vayacalculillo}.

\begin{cor}
If $T$ has $1$ or $2$ ends the characteristic class $\set{\theta}$
represented by the obstruction operator $\theta$ in
(\ref{thetacoh}) $(n> 2)$ is trivial
$$0=\set{\theta}\in H^1(\C{chain}_n(\ab(T))/\!\simeq\;,H^{n+2}(-,H_n\otimes\Z_2)),$$
therefore Steenrod invariants in the sense of Section \ref{epsi}
do exist for trees with $1$ or $2$ ends.
\end{cor}

An explicit definition of these Steenrod invariants will be given
later in this section. Before we consider the case of trees with
$3$ ends, which is a bit more complicated because the cup-product
mod $2$ does not vanish in this case, see Theorem
\ref{vayacalculillo}.

\begin{thm}\label{3fin}
If $T$ has $3$ ends then
$$\xymatrix{H^0(\C{chain}_n(\ab(T))/\!\simeq\;,H^{n+2}(-,\hat{\otimes}^2_TH_n\otimes\Z_2))\ar[d]^{(\bar{q}_T)_*}_\simeq\\
H^0(\C{chain}_n(\ab(T))/\!\simeq\;,H^{n+2}(-,\wedge^2_TH_n\otimes\Z_2))=\Z/2}$$
is an isomorphism. In particular the characteristic class
$\set{\theta}$ represented by the obstruction operator $\theta$ in
(\ref{thetacoh}) $(n> 2)$ is also trivial
$$0=\set{\theta}\in H^1(\C{chain}_n(\ab(T))/\!\simeq\;,H^{n+2}(-,H_n\otimes\Z_2)),$$
therefore Steenrod invariants in the sense of Section \ref{epsi}
does exist for trees with $3$ ends as well.
\end{thm}

\begin{proof}
We will use the following Bockstein long exact sequence
\begin{equation*}
\xymatrix{\vdots\ar[d]\\
H^k(\C{fp}(\C{vect}(T)),\ext^1_{\C{vect}(T)})\ar[d]^{(\bar{\tau}_T)_*}\\
H^k(\C{fp}(\C{vect}(T)),\ext^1_{\C{vect}(T)}(-,\hat{\otimes}^2_T))\ar[d]^{(\bar{q}_T)_*}\\
H^k(\C{fp}(\C{vect}(T)),\ext^1_{\C{vect}(T)}(-,\wedge^2_T))\ar[d]^{\beta}\\
H^{k+1}(\C{fp}(\C{vect}(T)),\ext^1_{\C{vect}(T)})\ar[d]\\\vdots}
\end{equation*}
which is isomorphic to (\ref{bocks}) by the localization theorem
for cohomology of categories in \cite{fcc}.

We showed in Theorem \ref{vayacalculillo} that
$$H^0(\C{fp}(\C{vect}(T)),\ext^1_{\C{vect}(T)}(-,\wedge^2_T))\simeq\Z/2$$
generated by the $0$-cocycle $c$ corresponding to the chain
cup-product mod $2$. This $0$-cocycle $c$ is additive, hence it is
determined by its values over elementary f. p.
$\C{vect}(T)$-modules, see Theorem \ref{clasclas},
$$c_\um{M}\in\ext^1_{\C{vect}(T)}(\um{M},\wedge^2_T\um{M}).$$ In
fact $c_\um{M}=0$ unless $\um{M}=\uf{M}\ul{V}^{(3,5)}$, and in
this case
$$0\neq c_{\uf{M}\ul{V}^{(3,5)}}\in\ext^1_{\C{vect}(T)}(\uf{M}\ul{V}^{(3,5)},\wedge^2_T\uf{M}\ul{V}^{(3,5)})=\Z/2,$$
see Lemma \ref{03}.

Let
$$\varepsilon\colon\wedge^2_T\uf{M}\ul{V}^{(3,5)}\twoheadrightarrow\hat{\otimes}^2_T\uf{M}\ul{V}^{(3,5)}$$
be a splitting of $\bar{q}_T$. Such a splitting exists by
Proposition \ref{splitmod2}. We define the $0$-cochain $d$ of
$\C{fp}(\C{vect}(T))$ with coefficients in
$\ext^1_{\C{vect}(T)}(-,\hat{\otimes}^2_T)$ as
\begin{equation*}
\begin{array}{ll}
d_\um{M}=0,&\hbox{if $\um{M}$ is an elementary
$\C{vect}(T)$-module $\um{M}\neq\uf{M}\ul{V}^{(3,5)}$;}\\
d_{\uf{M}\ul{V}^{(3,5)}}=\varepsilon_*(c_{\uf{M}\ul{V}^{(3,5)}}).
\end{array}
\end{equation*}
If we manage to prove that $d$ is a cocycle then immediately
$(\bar{q}_T)_*d=c$ and the surjectivity will follow.

The elementary $\C{vect}(T)$-module $\uf{M}\ul{V}^{(3,5)}$ has
only the identity and the trivial endomorphisms, see Proposition
\ref{subhom} (4), therefore all we have to do to prove that $d$ is
a cocycle is to check that given an an elementary
$\C{vect}(T)$-module $\um{M}\neq\uf{M}\ul{V}^{(3,5)}$ and
arbitrary morphisms $f\colon\um{M}\r\uf{M}\ul{V}^{(3,5)}$ and
$g\colon\uf{M}\ul{V}^{(3,5)}\r\um{M}$ the following equalities
hold
\begin{equation*}
\begin{array}{c}
0=f^*\varepsilon_*c_{\uf{M}\ul{V}^{(3,5)}}
\in\ext^1_{\C{vect}(T)}(\um{M},\hat{\otimes}^2_{T}\uf{M}\ul{V}^{(3,5)}),
\\{}\\
0=g_*\varepsilon_*c_{\uf{M}\ul{V}^{(3,5)}}
\in\ext^1_{\C{vect}(T)}(\uf{M}\ul{V}^{(3,5)},\hat{\otimes}^2_{T}\um{M}).
\end{array}
\end{equation*}
The first one follows easily from the equality
$f^*\varepsilon_*=\varepsilon_*f^*$ and the fact that $c$ is a
cocycle. In order to prove the second one it is enough to check that
for any elementary $\C{vect}(T)$-module
$\um{M}\neq\uf{M}\ul{V}^{(3,5)}$ either
\begin{equation*}
\hom_{\C{vect}(T)}(\uf{M}\ul{V}^{(3,5)},\um{M})=0
\end{equation*}
or
$$\ext^1_{\C{vect}(T)}(\uf{M}\ul{V}^{(3,5)},\hat{\otimes}^2_{T}\um{M})=0.$$
This is proved in Lemma \ref{04} below.

In order to finish the proof we just have to check that
$$H^0(\C{fp}(\C{vect}(T)),\ext^1_{\C{vect}(T)})=0.$$
This follows from the additivity of $0$-cocycles and the fact that
elementary f. p. $\C{vect}(T)$-modules have no self-extensions,
see Proposition \ref{noself}.
\end{proof}

In the proof of Theorem \ref{3fin} it is also implicitly contained
the proof of the following one, compare with Theorem
\ref{calculote}. In the statement we use the natural
transformation $\hat{\kappa}$ which is the composite
\begin{equation*}
\xymatrix@R=15pt{\ext^1_{\C{vect}(T)}(H_n\otimes\Z/2,\hat{\otimes}^2_TH_n\otimes\Z/2)\ar[d]^\simeq\\
\ext^2_{\ab(T)}(H_n\otimes\Z/2,\hat{\otimes}^2_TH_n\otimes\Z/2)\ar[d]^{\hat{p}^*}\\
\ext^2_{\ab(T)}(H_n,\hat{\otimes}^2_TH_n\otimes\Z/2)\ar[d]\\
H^{n+2}(-,\hat{\otimes}^2_TH_n\otimes\Z_2)}
\end{equation*}
given by Proposition \ref{ZaZ2}, the natural projection
$\hat{p}\colon\um{M}\twoheadrightarrow\um{M}\otimes\Z_2$, and the
universal coefficients spectral sequence (\ref{ucss}).

\begin{thm}\label{calculote2}
Given a tree $T$ with $3$ ends and a bounded chain complex $\cc_*$
in $\ab(T)$ concentrated in dimensions $\geq n$, if
$H_n\cc_*\otimes\Z/2$ contains $\uf{M}\ul{V}^{(3,5)}$ as a direct
summand $((\bar{q}_T)_*^{-1}\hat{p}_*\bar{\cup})_{\cc_*}$ then is
the image of the generator
\begin{equation*}
\xymatrix@R=15pt{
\Z/2=\ext^1_{\C{vect}(T)}(\uf{M}\ul{V}^{(3,5)},\hat{\otimes}^2_{T}\uf{M}\ul{V}^{(3,5)})\ar@{^{(}->}[d]\\
\ext^1_{\C{vect}(T)}(H_n\cc_*\otimes\Z/2,\hat{\otimes}^2_TH_n\cc_*\otimes\Z/2)\ar[d]^{\hat{\kappa}}\\
H^{n+2}(\cc_*,\hat{\otimes}^2_TH_n\cc_*\otimes\Z/2)}
\end{equation*}
Otherwise $((\bar{q}_T)_*^{-1}\hat{p}_*\bar{\cup})_{\cc_*}=0$.
\end{thm}

Now that we have proved the existence of Steenrod invariants for
trees $T$ with less than $4$ ends we proceed to define them
explicitly.

\begin{defn}
Let $(\cc_*,f_{n+2})$ be an $(n-1)$-reduced $T$-homotopy system
$(n>2)$. If $T$ has $1$ or $2$ ends its \emph{Steenrod invariant}
$$\wp_n(\cc_*,f_{n+2})\in H^{n+2}(\cc_*,H_n\cc_*\otimes\Z/2)$$ is the unique
element such that
$$(\bar{\tau}_T)_*\wp_n(\cc_*,f_{n+2})=\hat{p}_*\hat{\cup}_{(\cc_*,f_{n+2})}\in H^{n+2}(\cc_*,\hat{\otimes}^2_TH_n\cc_*\otimes\Z/2),$$
and if $T$ has $3$ ends it is the unique element satisfying
$$(\bar{\tau}_T)_*\wp_n(\cc_*,f_{n+2})=\hat{p}_*\hat{\cup}_{(\cc_*,f_{n+2})}-((\bar{q}_T)^{-1}_*\hat{p}_*\bar{\cup})_{\cc_*}\in H^{n+2}(\cc_*,\hat{\otimes}^2_TH_n\cc_*\otimes\Z/2).$$
\end{defn}

\begin{thm}\label{steen}
If $T$ has less than $4$ ends the Steenrod invariant of an
$(n-1)$-reduced $T$-homotopy system $(n>2)$ is well defined and
determines the obstruction operator $\theta$ in (\ref{thetacoh})
in the following sense, given a chain homomorphism
$\xi\colon\cc_*\r\cc_*'$ then
$$\theta_{(\cc_*,f_{n+2}),(\cc_*',g_{n+2})}(\xi)=\xi_*\wp_n{(\cc_*,f_{n+2})}-\xi^*\wp_n{(\cc_*',g_{n+2})}.$$
\end{thm}

\begin{proof}
The existence is derived from the fact that (\ref{shocorto}) is
short exact in this case, see Proposition \ref{splitmod2}, and the
equalities
\begin{equation*}
(\bar{q}_T)_*\hat{p}_*\hat{\cup}_{(\cc_*,f_{n+2})}=\hat{p}_*\bar{\cup}_{\cc_*}=0
\end{equation*}
if $T$ has $1$ or $2$ ends, see Definition \ref{cupcc} and Theorem
\ref{vayacalculillo}; or
\begin{equation*}
(\bar{q}_T)_*(\hat{p}_*\hat{\cup}_{(\cc_*,f_{n+2})}-((\bar{q}_T)^{-1}_*\hat{p}_*\bar{\cup})_{\cc_*})=\hat{p}_*\bar{\cup}_{\cc_*}-\hat{p}_*\bar{\cup}_{\cc_*}=0
\end{equation*}
if $T$ has $3$, see Definition \ref{cupcc} and Theorem \ref{3fin}.

The second part of the statement follows from Theorem \ref{cupobs}
and the facts that $(\bar{q}_T)^{-1}\hat{p}_*\bar{\cup}$ is a
$0$-cocycle and $(\bar{\tau}_T)_*$ is a monomorphism in this case,
see Proposition \ref{splitmod2}.
\end{proof}

Steenrod invariants for $T$-complexes are defined by using the
functor $r$ in Section \ref{hs}.

\begin{defn}
Let $T$ be a tree with less than $4$ ends, the \emph{Steenrod
invariant} of an $(n-1)$-reduced $T$-complex $X$ $(n>2)$ is
defined as
$$\wp_n(X)=\wp_n(rX)\in H^{n+2}(X,\hh_nX\otimes\Z/2).$$
\end{defn}

\section{The proper homotopy classification of $A^2_n$-polyhedra with less than $4$ ends $(n>2)$}\label{kia}

Given a fixed tree $T$ an \emph{$A^2_n$-polyhedron} is an
$(n-1)$-reduced $T$-complex $X$ of dimension $\leq n+2$, the
proper homotopy category of $A^2_n$-polyhedra was already
considered in Section \ref{hs} and is denoted by $\C{A}^2_n(T)$.
This category is equivalent to the proper homotopy category of
$(n-1)$-connected $T$-complexes of dimension $\leq n+2$, see
Proposition \ref{conexred}.


\begin{defn}
The objects of the category $\C{P}^2_n(T)$ are pairs $(\cc_*,\wp)$
where $\cc_*$ is a chain complex in $\ab(T)$ concentrated in
dimensions $n$, $n+1$ and $n+2$, together with an element
$$\wp\in H^{n+2}(\cc_*,H_n\cc_*\otimes\Z/2).$$
Morphisms $\xi\colon(\cc_*,\wp)\r(\cc_*',\wp')$ are chain homotopy
classes $\xi\colon\cc_*\r\cc_*'$ such that
$$\xi_*\wp=\xi^*\wp'\in H^{n+2}(\cc_*,H_n\cc_*'\otimes\Z/2).$$
\end{defn}

\begin{thm}\label{clasi1}
If $T$ has less than $4$ ends and $n>2$ the functor
$$\C{A}^2_n(T)\To\C{P}^2_n(T)\colon X\mapsto (\cc_*X,\wp_n(X))$$
induces a bijection between the sets of isomorphism classes of
objects.
\end{thm}

\begin{proof}
Recall from Section \ref{hs} that $\C{H}^2_n(T)$ denotes the
category of $(n-1)$-reduced $T$-homotopy systems $(\cc_*,f_{n+2})$
such that the chain complex $\cc_*$ is concentrated in dimensions
$n$, $n+1$ and $n+2$. By proposition \ref{isoob} this theorem is
equivalent to prove that the functor
$$\C{H}^2_n(T)/\!\simeq\;\To\C{P}^2_n(T)\colon (\cc_*,f_{n+2})\mapsto (\cc_*,\wp_n(\cc_*,f_{n+2}))$$
induces a bijection between the sets of isomorphism classes of
objects in both categories.

By Proposition \ref{so} given $(\cc_*,\wp)$ in $\C{P}^2_n(T)$
there exists an $(n-1)$-reduced $T$-homotopy system
$(\cc_*,f_{n+2})$. With the notation in Section \ref{cohcat} it is
easy to see by applying Proposition \ref{steen} that the new
homotopy system $(\cc_*,f_{n+1})+(\wp-\wp_n(\cc_*,f_{n+1}))$ has
Steenrod invariant $\wp$.

Suppose now that we have two homotopy systems $(\cc_*,f_{n+2})$,
$(\cc_*',g_{n+2})$ and an isomorphism in $\C{P}^2_n(T)$
$$\xi\colon(\cc_*,\wp_n(\cc_*,f_{n+2}))\To (\cc_*',\wp_n(\cc_*',g_{n+2})).$$
By Proposition \ref{steen} there exists a morphism of homotopy
systems $$(\xi,\eta)\colon(\cc_*,f_{n+2})\To(\cc_*',g_{n+2}).$$ This
morphism is an isomorphism in the homotopy category because the
functor $(\cc_*,f_{n+2})\mapsto\cc_*$ fits into an exact sequence in
the sense of Section \ref{cohcat}.
\end{proof}

\begin{cor}\label{clasi2}
If $T$ has less than $4$ ends, given a chain complex $\cc_*$ in
$\ab(T)$ concentrated in dimensions $n$, $n+1$ and $n+2$ for some
$n>2$ the set of proper homotopy classes of $A^2_n$-polyhedra $X$
with proper cellular chain complex $\cc_*X$ homotopy equivalent to
$\cc_*$ is in bijective correspondence with the orbit set
$$H^{n+2}(\cc_*,H_n\cc_*\otimes\Z/2)/\operatorname{Aut}(\cc_*)$$
of the right action of the group $\operatorname{Aut}(\cc_*)$ of
self-homotopy equivalences of $\cc_*$ on the cohomology group
$H^{n+2}(\cc_*,H_n\cc_*\otimes\Z/2)$ given by
$$a^\xi=\xi^*\xi_*^{-1}(a).$$
\end{cor}

\begin{proof}
One can readily check that the function sending $X$ to the class
of $\wp_n(X)$ in the orbit set induces the desired bijection.
\end{proof}

\section{Moore spaces in proper homotopy theory}\label{mspht}

\begin{defn}
Given a fixed tree $T$ a Moore space of degree $n$ is a
$1$-connected $T$-complex $X$ whose unique non-trivial proper
homology $\ab(T)$-module is $\hh_nX$.
\end{defn}

By Theorem \ref{kernel} and Proposition \ref{so} there are Moore
spaces in degree $n\geq 2$ for any finitely presented
$\ab(T)$-module. Moreover, the proper homotopy category
$\C{M}_n(T)$ of Moore spaces of degree $n$ can be regarded as a
full subcategory of $\C{A}^2_n(T)$. In order to specify the
results in the previous section to this subcategory we give the
following definition.

\begin{defn}
Objects in the category $\C{PM}^2_n(T)$ are pairs $(\um{M},\wp)$
given by a finitely presented $\ab(T)$-module $\um{M}$ together
with an element
$$\wp\in \ext^2_{\ab(T)}(\um{M},\um{M}\otimes\Z/2).$$
Morphisms $\xi\colon(\um{M},\wp)\r(\um{M}',\wp')$ are
$\ab(T)$-module morphisms $\varphi\colon\um{M}\r\um{M}'$ such that
$$\varphi_*\wp=\varphi^*\wp'\in \ext^2_{\ab(T)}(\um{M},\um{M}'\otimes\Z/2).$$
\end{defn}

\begin{thm}\label{clasimoore1}
If $T$ has less than $4$ ends and $n>2$ the functor
$$\C{M}^2_n(T)\To\C{PM}^2_n(T)\colon X\mapsto (\hh_nX,\wp_n(X))$$
induces a bijection between the sets of isomorphism classes of
objects.
\end{thm}

\begin{cor}\label{clasimoore2}
If $T$ has less than $4$ ends, given a f. p. $\ab(T)$-module
$\um{M}$ the set of proper homotopy classes of Moore spaces $X$ of
degree $n$ with proper cellular homology $\hh_nX$ isomorphic to
$\um{M}$ is in bijective correspondence with the orbit set
$$\ext^2_{\ab(T)}(\um{M},\um{M}\otimes\Z/2)/\operatorname{Aut}(\um{M})$$
of the right action given by
$$a^\varphi=\varphi^*\varphi_*^{-1}(a).$$
\end{cor}

\begin{rem}\label{tata}
The easiest explicit computation one can do with the previous
corollary is the following. Consider the $\C{vect}(T_1)$-module
$\um{A}\oplus\um{C}$, see Theorem \ref{clasclas}. By Proposition
\ref{ZaZ2} and \cite{rca} 7.17 and A.3
$$\ext^2_{\ab(T)}(\um{A}\oplus\um{C},(\um{A}\oplus\um{C})\otimes\Z/2)=\Z/2,$$
in particular the action of
$\operatorname{Aut}(\um{A}\oplus\um{C})$ on this group is trivial
and there are exactly $2$ proper Moore space of degree $n>2$ with
proper homology $\um{A}\oplus\um{C}$.

The first example of non homotopy equivalent Moore spaces with the
same proper homology $\ab(T_1)$-module $\um{S}$ was discovered in
\cite{tsukuba95}. Although they used the pro-categorical approach to
proper homology one can check by using the Brown-Grossman functor,
see \cite{iht} V.3.10, and the projective resolution for $\um{S}$ in
\cite{tsukuba95} Appendix A that $\um{S}$ in our algebraic context
is $\um{S}=\um{R}\oplus\um{B}$. By Proposition \ref{ZaZ2},
\cite{rca} 7.12 and A.3
$$\dim_{\uf{F}_2}\ext^2_{\ab(T)}(\um{R}\oplus\um{B},(\um{R}\oplus\um{B})\otimes\Z/2)=2^{\aleph_0},$$
\end{rem}

An immediate consequence of Corollary \ref{clasimoore2} is the
following.

\begin{cor}
If $T$ has less than $4$ ends and $n>2$ there is a unique Moore
space of degree $n$ with homology $\um{M}$ if and only if
$\ext^2_{\ab(T)}(\um{M},\um{M}\otimes\Z/2)=0$.
\end{cor}

In the next corollary we completely determine the uniqueness of
proper Moore spaces of degree $n>2$ with only one end and exponent
$2$ in proper homology. It follows from the previous corollary,
Proposition \ref{ZaZ2} and the homological computations in the
appendix of \cite{rca}.

\begin{cor}
If $T$ has only $1$ end and $\um{M}$ is a f. p. $\C{vect}(T)$-module
then there exists a unique Moore space $X$ of degree $n>2$ with
homology $\um{M}$ if and only if this $\C{vect}(T)$-module satisfies
one of the following conditions:
\begin{enumerate}
\item it does not contain $\um{A}$ as a direct summand,

\item it does not contain either $\um{R}$ or $\um{C}$ as a direct summand,

\item it does not contain either $\um{B}$ or $\um{C}$ as a direct
summand.
\end{enumerate}
\end{cor}

Now by Theorem \ref{clasclas} we get the following result.

\begin{cor}
It $T$ has only one end there are infinite f. p.
$\C{vect}(T)$-modules $\um{M}$ such that there exist different
Moore spaces in degree $n>2$ with homology $\um{M}$.
\end{cor}

\begin{rem}\label{tata2}
The lower part of the Whitehead long exact sequence
$$\hh_{n+2}X\st{b_{n+2}}\r\Gamma_{n+1}X\st{i_{n+1}}\r\Pi_{n+1}X\st{h_{n+1}}\twoheadrightarrow\hh_{n+1}X,\;\;\hh_nX,$$
does not classify proper homotopy types of $A^2_n$-polyhedra. By
Theorem \ref{lwm} the sequence is the same for any two Moore spaces
of degree $n$ with the same homology. However, as we can see for
example in the previous corollary, the homology does not determine
the proper homotopy type.
\end{rem}

The following result on the existence of co-H-structures on degree
$2$ Moore spaces is a consequence of Theorem \ref{calculote} and
Proposition \ref{cohcupcc}.

\begin{cor}\label{kio}
Let $\um{M}$ be a f. p. $\C{vect}(T)$-module. If $T$ has $1$ or
$2$ ends all Moore spaces in degree $2$ with homology $\um{M}$ are
co-H-spaces. If $T$ has three ends such a Moore space is a
co-H-space if and only if $\um{M}$ does not contain
$\uf{M}\ul{V}^{(3,5)}$ as a direct summand.
\end{cor}

\begin{rem}\label{kio2}
The module $\uf{M}\ul{V}^{(3,5)}$ has the curious property that
there exists a unique Moore space of degree $n>2$ with that
homology, but all possible Moore spaces with
homology $\uf{M}\ul{V}^{(3,5)}$ in degree $2$ are not co-H-spaces.
\end{rem}

\chapter{Computations in controlled algebra}\label{chIX}

In this section we perform technical algebraic computations leading to the main homotopical results of this paper already presented in the previous chapter. For this we use the representation theory of the algebras $k(\F(T))$ for $k$ a field considered in \cite{rca}. 

\section{Review of the representation theory of the algebras
$k(\F(T))$}\label{rrt}

Let $k$ be an arbitrary field. Recall from Section \ref{freecontrol}
that $k(\F(T))$ is a $k$-algebra Morita equivalent to the small
additive category $\C{M}_k(T)$.

In \cite{rca} 1.1 we determined the representation type of the
$k$-algebras $k(\F(T))$ for $k$ an arbitrary field in terms of the
number of Freudenthal ends of $T$.

\begin{thm}\label{rep}
The representation type of $k(\F(T))$ is
\medskip
\begin{center}
\begin{tabular}{|c|c|}
  \hline
  $\card \F(T)$ & type \\
  \hline
  $<4$ & finite \\
  $=4$ & tame \\
  $>4$ & wild \\
  \hline
\end{tabular}
\end{center}
\end{thm}

In the finite and tame cases we give in \cite{rca} explicit
classification theorems for f. p. $\C{M}_k(T)$-modules. In order
to state them we fix an explicit tree $T_n$ with a finite number
$n$ of Freudenthal ends, see Section 3.1 in \cite{rca}. The vertex
set of $T_n$ is $$T_n^0=\set{v_0}\cup\set{v_m^1,\dots
v_m^n}_{m\geq 1}$$ and there are edges joining $v_0$ with $v_1^i$
and $v_m^i$ with $v^i_{m+1}$ $(1\leq i\leq n, m\geq 1)$. Notice
that we can identify $T_1$ with the half-line $\Real_+$ and
$T^0_1$ with the non-negative integers $\N_0$. Consider the
$\N_0\times\N_0$ matrices $\mat{A}$ and $\mat{B}$ with entries in
$k$ given by
\begin{itemize}
\item $\mat{a}_{i+1,i}=1$ $(i\in\NN)$ and $\mat{a}_{ij}=0$ in
other cases,

\item $\mat{b}_{\frac{n(n+1)}{2}+i,\frac{(n-1)n}{2}+i}=1$ for any
$n>0$ and $0\leq i<n$, and $\mat{b}_{ij}=0$ otherwise.
\end{itemize}
Moreover, let $\mat{I}$ be the identity matrix. The f. p.
$\C{M}_k(\Real_+)$-modules $\um{A},\um{B},\um{C},\um{B}_\infty$ and
$\um{C}_\infty$ are defined as the cokernels of the endomorphisms of
the free $T$-controlled $k$-module $k\grupo{\N_0}_\delta$, where
$\delta\colon\N_0\hookrightarrow \Real_+$ is the inclusion, given by
the matrices $\mat{A}$, $\mat{I}-\mat{A}$, $\mat{I}-\mat{A}^t$,
$\mat{I}-\mat{B}$ and $\mat{I}-\mat{B}^t$, respectively. Here
$(-)^t$ denotes transposition. Furthermore, we call $\um{R}$ to the
f. g. free $\C{M}_k(\Real_+)$-module $k\grupo{\N_0}_\delta$.

There are $n$ proper homotopy classes of maps $f_i\colon\Real_+\r
T_n$ $(1\leq i\leq n)$ corresponding to the inclusions of the $n$
different Freudenthal ends of $T_n$. We can take $f_i$ as the
simplicial map with $f_i(0)=v_0$ and $f_i(j)=w_n^i$ $(j\geq 1)$. As
we know from Section \ref{freecontrol} these maps induce ``change of
tree'' functors
\begin{equation*}
\uf{F}^i=\uf{F}^{f_i}\colon\C{M}_k(\Real_+)\To\C{M}_k(T_n).
\end{equation*}

For trees with more than one Freudenthal end the classification
theorems depend on the classification of what we call
finite-dimensional rigid $n$-subspaces, see Section 8 of
\cite{rca}, which are special representations of the $n$-subspace
quiver $Q_n$
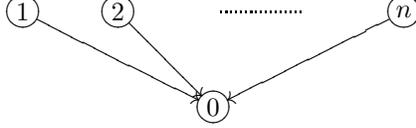
$$\xygraph{*+=[o]+[F]{0} (:@{<-}[ull] *+=[o]+[F]{1}, :@{<-}[ul]
*+=[o]+[F]{2}, :@{}[u] {} :@{.}[r] {}, :@{<-}[urr] *+=[o]+[F]{n})}$$
Recall that a representation of $Q_n$ is a diagram of $k$-vector
spaces indexed by $Q_n$, i. e. $n+1$ vector spaces
$V_0,V_1,\dots,V_n$ together with homomorphisms $V_i\r V_0$ $(0\leq
i\leq n)$. These representations form an abelian category which
coincides with the category of modules over the path algebra $kQ_n$.
An $n$-subspace is a representation such that all the morphisms
$V_i\r V_0$ are inclusions of subspaces. There is an exact, full and
faithful functor between the categories of finite-dimensional
$n$-subspaces and f. p. $\C{M}_k(T_n)$-modules, 
see Section 9 of \cite{rca},
$$\uf{M}\colon\C{sub}_n^{\mathrm{fin}}\To\C{fp}(\C{M}_k(T_n)).$$

The following classification theorem is proved in \cite{rca} 10.5
for arbitrary fields. In the statement we use the following
terminology: a solution to the \emph{decomposition problem} in a
small additive category $\mathbf{A}$ consists of a set of objects
(which we call \emph{elementary objects}) and of a set of
isomorphisms (\emph{elementary isomorphisms}) between finite direct
sums of elementary objects. These sets must satisfy that any object
in $\mathbf{A}$ is isomorphic to a finite direct sum of elementary
ones, and any isomorphism relation between two such direct sums can
be derived from the elementary isomorphisms.

\begin{thm}\label{clasclas}
Let $$\set{\ul{V}^{(n,j)}}_{j\in J_n}$$ be the set of indecomposable
rigid $n$-subspaces $(n\geq 1)$, see Proposition \ref{losn}. There
is a solution to the decomposition problem in the category of f. p.
$\C{M}_k(T_n)$-modules given by the following $1+5n+\card J_n$
elementary modules $(1\leq i\leq n,j\in J_n)$
$$\uf{F}_*^1\um{A},\;\uf{F}^i_*\um{R},\;\uf{F}^i_*\um{B},\;\uf{F}^i_*\um{B}_\infty,
\;\uf{F}^i_*\um{C},\;\uf{F}^i_*\um{C}_\infty,\;\uf{M}\ul{V}^{(n,j)},\;\;\;$$
and $6n$ elementary isomorphisms $(1\leq i\leq n)$
$$\uf{F}_*^1\um{A}\oplus\uf{F}^i_*\um{R}\simeq\uf{F}^i_*\um{R},\;
\uf{F}^i_*\um{R}\oplus\uf{F}^i_*\um{R}\simeq\uf{F}^i_*\um{R},\;
\uf{F}_*^i\um{B}\oplus\uf{F}^i_*\um{B}_\infty\simeq\uf{F}^i_*\um{B}_\infty,\;$$
$$\uf{F}_*^i\um{B}_\infty\oplus\uf{F}^i_*\um{B}_\infty\simeq\uf{F}^i_*\um{B}_\infty,\;
\uf{F}_*^i\um{C}\oplus\uf{F}^i_*\um{C}_\infty\simeq\uf{F}^i_*\um{C}_\infty,\;
\uf{F}_*^i\um{C}_\infty\oplus\uf{F}^i_*\um{C}_\infty\simeq\uf{F}^i_*\um{C}_\infty.$$
\end{thm}

For $n=1, 2$ and $3$ indecomposable rigid $n$-subspaces are listed
in the following proposition, see \cite{rca} 8.7.

\begin{prop}\label{losn}
The following are complete lists of (representatives of the
isomorphism classes of) indecomposable rigid $n$-subspaces for
$n<4$
\begin{itemize}
\item $n=1$, none, \item $n=2$, $\ul{V}^{(2,1)}=\left(k\rightarrow
k\leftarrow k\right)$, \item $n=3$,
$$\ul{V}^{(3,1)}=\left(\begin{array}{ccccc}
   &  & k &  &  \\
   &  & \downarrow &  &  \\
  k & \rightarrow & k & \leftarrow & 0
\end{array}\right),\;\;\;
\ul{V}^{(3,2)}=\left(\begin{array}{ccccc}
   &  & 0 &  &  \\
   &  & \downarrow &  &  \\
  k & \rightarrow & k & \leftarrow & k
\end{array}\right),$$
$$\ul{V}^{(3,3)}=\left(\begin{array}{ccccc}
   &  & k &  &  \\
   &  & \downarrow &  &  \\
  0 & \rightarrow & k & \leftarrow & k
\end{array}\right),\;\;\;
\ul{V}^{(3,4)}=\left(\begin{array}{ccccc}
   &  & k &  &  \\
   &  & \downarrow &  &  \\
  k & \rightarrow & k & \leftarrow & k
\end{array}\right),$$
$$\ul{V}^{(3,5)}=\left(\begin{array}{ccccc}
   &  & k\grupo{x+y} &  &  \\
   &  & \downarrow &  &  \\
  k\grupo{x} & \rightarrow & k\grupo{x,y} & \leftarrow & k\grupo{y}
\end{array}\right).$$
\end{itemize}
\end{prop}

For $n=4$ the list of indecomposable rigid $n$-subspaces is infinite
and very complicated to describe. It can be extracted from
\cite{nazarova}, see the remark in \cite{rca} 8.8.

\section{Some ``change of tree'' computations}

\begin{prop}\label{cosadea}
For all $n\geq 1$ and $1\leq i\leq n$ we have that
$\uf{F}^1_*\um{A}=\uf{F}^i_*\um{A}$.
\end{prop}

This proposition follows either from \cite{rca} 7.7 or from
\cite{rca} 10.3.

\begin{lem}\label{0casito}
Let $k\grupo{A}_\alpha$ be  a $T_1$-controlled $k$-module and let
$k\grupo{B}_\beta$ be a $T_n$-controlled $k$-module such that there
exists a neighbourhood $U$ of the $i^{\text{th}}$ Freudenthal end of
$T_n$ in $\hat{T}_n$ with $\beta^{-1}(U)=\emptyset$. Then there is a
natural abelian group isomorphism
$$\hom_{\C{M}_k(T_n)}(\uf{F}^i_*k\grupo{A}_\alpha,k\grupo{B}_\beta)\simeq\bigoplus_{A}k\grupo{B}
\subset\prod_{A}k\grupo{B}=\hom_{k}(k\grupo{A},k\grupo{B}).$$
\end{lem}

\begin{proof}
It is enough to check that any controlled homomorphism
$\varphi\colon\uf{F}^i_*k\grupo{A}_\alpha=k\grupo{A}_{f_i\alpha}\r
k\grupo{B}_\beta$ vanishes in almost all $A$. Since $\varphi$ is
controlled there exists another neighbourhood $V\subset U$ of the
$i^{\mathrm{th}}$ Freudenthal end of $T_n$ in $\hat{T}_n$ such
that
$\varphi(\alpha^{-1}(f_i^{-1}(U)))=\varphi((f_i\alpha)^{-1}(U))\subset
k\grupo{\beta^{-1}(U)}=0$. The closure of the set
$T_1-f_i^{-1}(U)\subset T_1$ is compact because
$(f_i\alpha)^{-1}(U)$ is a neighbourhood of the unique Freudenthal
end of $T_1$ hence $\alpha^{-1}(T_1-f_i^{-1}(U))$ is finite and
the lemma follows.
\end{proof}

\begin{cor}\label{0casito2}
Given two $T_1$-controlled $k$-modules $k\grupo{A}_\alpha$ and
$k\grupo{B}_\beta$, if $i\neq j$ there is a natural isomorphism
$$\hom_{\C{M}_k(T_1)}(k\grupo{A}_\alpha,(\uf{F}^i)^*\uf{F}^j_*k\grupo{B}_\beta)\simeq\bigoplus_{A}k\grupo{B}
\subset\prod_{A}k\grupo{B}=\hom_{k}(k\grupo{A},k\grupo{B}).$$
\end{cor}

This corollary follows from (\ref{yonedacom}), Lemma \ref{0casito}
and the fact that $(\uf{F}^i)^*$ is right adjoint to $\uf{F}^i_*$.

\begin{cor}\label{0casito3}
If $\um{M}$ a $\C{M}_k(T_1)$-module f. p. by the controlled
homomorphism $\varphi\colon k\grupo{B}_\beta\r k\grupo{C}_\gamma$
and $i\neq j$ then there is a natural isomorphism
\begin{small}
$$\hom_{\C{M}_k(T_1)}(k\grupo{A}_\alpha,(\uf{F}^i)^*\uf{F}^j_*\um{M})\simeq\bigoplus_{A}
\frac{k\grupo{C}}{\varphi(k\grupo{B})}
\subset\prod_{A}\frac{k\grupo{C}}{\varphi(k\grupo{B})}=
\hom_{k}\left(k\grupo{A},\frac{k\grupo{C}}{\varphi(k\grupo{B})}\right).$$
\end{small}
\end{cor}

This follows from Corollary \ref{0casito2} and the right-exactness
of $(\uf{F}^i)^*$ and $\uf{F}^j_*$.

\begin{lem}\label{0casito4}
There is an isomorphism natural in the $T_1$-controlled $k$-module
$k\grupo{A}_\alpha$
$$\hom_{\C{M}_k(T_1)}(k\grupo{A}_\alpha,\bigoplus\limits_{B}\um{A})\simeq\bigoplus_{A}
k\grupo{B} \subset\prod_{A}k\grupo{B}=
\hom_{k}(k\grupo{A},k\grupo{B}).$$
\end{lem}

\begin{proof}
Since
$$\hom_{\C{M}_k(T_1)}(k\grupo{A}_\alpha,\bigoplus\limits_{B}\um{A})=\bigoplus\limits_{B}\hom_{\C{M}_k(T_1)}(k\grupo{A}_\alpha,\um{A})$$
it is enough to prove the lemma for $B$ a singleton. By \cite{rca}
7.7 we have to prove that if $k\grupo{e}_\phi$ is a
$T_1$-controlled $k$-module with one generator $e$ and
$\phi(e)=v_0$ there is an isomorphism natural in
$k\grupo{A}_\alpha$
$$\hom_{\C{M}_k(T_1)}(k\grupo{A}_\alpha,k\grupo{e}_\phi)\subset\bigoplus_{A}
k\grupo{e} \subset\prod_{A}k\grupo{e}=
\hom_{k}(k\grupo{A},k\grupo{e}),$$ i. e. any controlled homomorphism
$k\grupo{A}_\alpha\r k\grupo{e}_\phi$ vanishes in almost all $A$.
This can be easily checked by using the definition of controlled
homomorphism as in the proof of Lemma \ref{0casito}. In fact it can
also be derived from Proposition \ref{cosadea} and Corollary
\ref{0casito2} by using the fact that $\uf{F}^i$ is fully faithful,
see Proposition \ref{fullis}, and hence
$(\uf{F}^i)^*\uf{F}^i_*\simeq 1$, see Section \ref{modules}.
\end{proof}

\begin{prop}\label{cambios}
Given $i\neq j$ there are isomorphisms
\begin{enumerate}
\item $(\uf{F}^i)^*\uf{F}^j_*\um{B}\simeq\um{A}$,
\item $(\uf{F}^i)^*\uf{F}^j_*\um{B}_\infty\simeq
\bigoplus\limits_0^\infty\um{A}$,
\item $(\uf{F}^i)^*\uf{F}^j_*\um{C}=0$,
\item $(\uf{F}^i)^*\uf{F}^j_*\um{C}_\infty=0$.
\end{enumerate}
\end{prop}

\begin{proof}
By Corollary \ref{0casito3} and Lemma \ref{0casito4} it is enough
to notice that 
\begin{eqnarray*}
\dim k\grupo{\N_0}/(\mat{I}-\mat{A})k\grupo{\N_0}&=&1, \\
\dim k\grupo{\N_0}/(\mat{I}-\mat{B})k\grupo{\N_0}&=&\aleph_0,\\
\dim k\grupo{\N_0}/(\mat{I}-\mat{A}^t)k\grupo{\N_0}&=&0,\\
\dim k\grupo{\N_0}/(\mat{I}-\mat{B}^t)k\grupo{\N_0}&=&0.
\end{eqnarray*}
See the proof of \cite{rca} 7.3
\end{proof}

\begin{prop}\label{submod}
Let $\ul{W}^i$ $(0\leq i\leq n)$ be the $n$-subspace with $\dim
W_0^i=\dim W_i^i=1$ (if $i>0$) and $\dim W_j^i=0$ otherwise, then
$(1\leq i\leq n)$
\begin{enumerate}
\item $\uf{M}\ul{W}^0=\uf{F}^i_*\um{A}$,

\item $\uf{M}\ul{W}^i=\uf{F}^i_*\um{B}$.
\end{enumerate}
\end{prop}

\begin{proof}
By using the commutativity of (\ref{yonedacom}) and the
right-exactness of $\uf{F}^i_*$ one obtains a finite presentation
of $\uf{F}^i_*\um{A}$ and $\uf{F}^i_*\um{B}$ as an
$\C{M}_k(T_n)$-module from \cite{rca} 7.7 and the very definition
$\um{B}$. It is easy to se that this presentation coincides with
the presentation of $\uf{M}\ul{W}^0$ or $\uf{M}\ul{W}^i$,
respectively, given in the proof of \cite{rca} 9.1.
\end{proof}

\section{Computation of some $\hom$ and $\ext$
groups}\label{homext}

In this section we calculate some $\hom$ and $\ext$ groups of
$\ab(T)$-modules and $\C{vect}(T)$-modules which play a role in
computations of cohomology groups of categories carried out in
previous sections. The main tools are contained in the appendix of
\cite{rca}.

\begin{prop}\label{mashom}
The following equalities hold
\begin{enumerate}
\item $\dim\hom_{\C{vect}(T_1)}(\um{B},\um{A})=0$,

\item $\dim\hom_{\C{vect}(T_1)}(\um{A},\um{B})=1$,

\item $\dim\hom_{\C{vect}(T_n)}(\uf{F}^i_*\um{B},\uf{F}^j_*\um{B}_\infty)=0$ for $i\neq
j$,

\item $\dim\hom_{\C{vect}(T_1)}(\um{C},\um{B}_\infty)=0$,

\item
$\dim\ext^1_{\C{vect}(T_n)}(\uf{F}^j_*\um{B},\uf{F}^i_*\um{B}_\infty)=0$.
\end{enumerate}
\end{prop}

\begin{proof}
The equalities (1) and (2) follow easily from Proposition
\ref{submod}. Moreover
\begin{eqnarray*}
\hom_{\C{vect}(T_n)}(\uf{F}^i_*\um{B},\uf{F}^j_*\um{B}_\infty)&=&\hom_{\C{vect}(T_1)}(\um{B},(\uf{F}^i)^*\uf{F}^j_*\um{B}_\infty)\\
&=&\hom_{\C{vect}(T_1)}(\um{B},\bigoplus_0^\infty\um{A})\\
&=&\bigoplus_0^\infty\hom_{\C{vect}(T_1)}(\um{B},\um{A})\\&=&0.
\end{eqnarray*}
Here we use that $(\uf{F}^i)^*$ is right-adjoint to $\uf{F}^i_*$,
Proposition \ref{cambios} (2), and the fact that all finitely
presented $\C{vect}(T_1)$-modules are small as a consequence of
\cite{rca} 7.15.

By using \cite{rca} 7.11 one readily sees that
$\hom_{\C{vect}(T_1)}(\um{M},\um{B}_\infty)=0$ whenever $\um{M}$
is the cokernel of a controlled homomorphism between free
$T_1$-controlled $\uf{F}_2$-modules $\varphi\colon
\uf{F}_2\grupo{A}_\alpha\r \uf{F}_2\grupo{B}_\beta$ whose
underlying vector space homomorphism is surjective $\varphi\colon
\uf{F}_2\grupo{A}\twoheadrightarrow \uf{F}_2\grupo{B}$, and this
happens with $\um{C}$, see the proof of \cite{rca} 7.3.

Finally we have that
\begin{eqnarray*}
\ext^1_{\C{vect}(T_n)}(\uf{F}^j_*\um{B},\uf{F}^i_*\um{B}_\infty)&=&\ext^1_{\C{vect}(T_1)}(\um{B},(\uf{F}^j)^*\uf{F}^i_*\um{B}_\infty)\\
(\hbox{if $i=j$}) &=&
    \ext^1_{\C{vect}(T_1)}(\um{B},\um{B}_\infty)=0,  \\
  (\hbox{if $i\neq j$}) &=&
  \ext^1_{\C{vect}(T_1)}(\um{B},\bigoplus_0^\infty\um{A})=\bigoplus_0^\infty\ext^1_{\C{vect}(T_1)}(\um{B},\um{A})=0.
\end{eqnarray*}
Here we use Propositions \ref{esexacto}, \ref{cambios} (2),
\cite{rca} 7.12 and A.3 (5), and the fact that f. p.
$\C{vect}(T_1)$-modules are small, as noticed above.
\end{proof}

\begin{prop}\label{noself}
If $T$ has less than $4$ ends the elementary finitely presented
$\C{vect}(T)$-modules in Theorem \ref{clasclas} have no
self-extensions.
\end{prop}

\begin{proof}
For $T$ with only one end this follows from \cite{rca} 7.12, 7.17
and the fact that $\um{A}$ and $\um{R}$ are projective, see
\cite{rca} 7.7. For trees with $2$ or $3$ ends it follows in
addition from Proposition \ref{esexacto}, \cite{rca} A.4 and the
well-known fact from representation theory that indecomposable
representations over Dynkin diagrams have no self-extensions.
\end{proof}

We are interested in computing $\ext^2$ groups of
$\C{vect}(T)$-modules regarded as $\ab(T)$-modules. For this we
use the change-of-rings spectral sequence associated to the
natural projection $\Z(\F(T))\twoheadrightarrow\uf{F}_2(\F(T))$,
compare \cite{ugss} 10.2 (c), or equivalently the
change-of-ringoids spectral sequence associated to the functor
$-\otimes\Z/2\colon\ab(T)\rightarrow\C{vect}(T)$, see
(\ref{expli}).

\begin{prop}\label{ZaZ2}
For any pair of $\C{vect}(T)$-modules $\um{M}$, $\um{N}$, the first
one finitely presented, there is a natural isomorphism
$$\ext^2_{\ab(T)}(\um{M},\um{N})\simeq\ext^1_{\C{vect}(T)}(\um{M},\um{N}).$$
Moreover, finitely presented $\C{vect}(T)$-modules have projective
dimension $\leq 1$.
\end{prop}

\begin{proof}
For the proof of this proposition we prefer to work with modules
over rings. To simplify we write $\F=\F(T)$. As a $\Z(\F)$-module
$\uf{F}_2(\F)$ is the cokernel of the multiplication by $2$
homomorphism $2\colon\Z(\F)\hookrightarrow\Z(\F)$ therefore the
$E_2$-term of the change-of-rings spectral sequence
$$E^{p,q}_2=\ext^q_{\uf{F}_2(\F)}(\tor^{\Z(\F)}_p(\um{M},\uf{F}_2(\F)),\um{N})\Rightarrow\ext^{p+q}_{\Z(\F)}(\um{M},\um{N})$$
is concentrated in $0\leq p\leq 1$. Moreover,
$$\tor^{\Z(\F)}_p(\um{M},\uf{F}_2(\F))=\um{M},\;\;p=0,1.$$

F. p. $\uf{F}_2(\F)$-modules are also f. p. as $\Z(\F)$-modules so
by Corollary \ref{pd2} they have projective dimension $\leq 2$ and
hence the $E_\infty$-term is concentrated in $0\leq p+q\leq2$, in
particular for any $q\geq 2$
$$\ext^q_{\uf{F}_2(\F)}(\um{M},\um{N})=E^{1,q}_2=E^{1,q}_\infty\hookrightarrow\ext^{1+q}_{\Z(\F)}(\um{M},\um{N})=0,$$
therefore the $\uf{F}_2(\F)$-module $\um{M}$ has projective
dimension $\leq 1$ and the $E_2$-term is concentrated in $0\leq
p,q\leq 1$, so
$$\ext^1_{\uf{F}_2(\F)}(\um{M},\um{N})=E^{1,1}_2=E^{1,1}_\infty=\ext^2_{\Z(\F)}(\um{M},\um{N}).$$
\end{proof}

By \cite{rca} A.4 the computation of $\ext^1$ groups of
$\C{vect}(T_n)$-modules coming from finite-dimensional
$n$-subspaces can be reduced to the computation of the
corresponding $\ext^1$ groups in the category of representations
of $Q_n$. In this category we have a powerful tool to carry out
this kind of computations, the bilinear form $\grupo{-,-}_{Q_n}$.
This is a homomorphism
$$\grupo{-,-}_{Q_n}\colon(\oplus_0^n\Z)\otimes_\Z(\oplus_0^n\Z)\To\Z,$$
which satisfies the following formula
\begin{equation}\label{extform}
\grupo{\ul{\dim}\;\ul{V},\ul{\dim}\;\ul{W}}_{Q_n}=\dim\hom_{kQ_n}(\ul{V},\ul{W})-\dim\ext^1_{kQ_n}(\ul{V},\ul{W})
\end{equation}
for any pair of finite-dimensional $n$-subspaces $\ul{V}$ and
$\ul{W}$. Here $$\ul{\dim}\;\ul{V}=(\dim V_0,\dots,\dim
V_n)\in\oplus_0^n\Z$$ is the \emph{dimension vector} of the
$n$-subspace $\ul{V}$. The bilinear form is defined by the
following formula
$$\grupo{\ul{a},\ul{b}}_{Q_n}=\sum_{i=0}^na_ib_i-\sum_{i=1}^na_ib_0,$$
where $\ul{a}=(a_0,\dots, a_n)$ and $\ul{b}=(b_0,\dots, b_n)$ are
elements in $\oplus_0^n\Z$.

The following result can be easily checked by hand.

\begin{prop}\label{subhom}
Consider the $3$-subspaces defined in Propositions \ref{losn} and
\ref{submod}. We have that
\begin{enumerate}
\item $\dim\hom_{kQ_3}(\ul{V}^{(3,5)},\ul{W}^i)=0$ $(0\leq i\leq
3)$,

\item $\dim\hom_{kQ_3}(\ul{V}^{(3,5)},\ul{V}^{(3,i)})=1$ $(1\leq i\leq
3)$,

\item $\dim\hom_{kQ_3}(\ul{V}^{(3,5)},\ul{V}^{(3,4)})=2$,

\item $\dim\hom_{kQ_3}(\ul{V}^{(3,5)},\ul{V}^{(3,5)})=1$.
\end{enumerate}
\end{prop}

The following proposition follows immediately from the previous
one by using formula (\ref{extform}).

\begin{prop}\label{subext}
For the $3$-subspaces defined in Propositions \ref{losn} and
\ref{submod} we have that
\begin{enumerate}
\item $\dim\ext^1_{kQ_3}(\ul{V}^{(3,5)},\ul{W}^0)=1$,
\item $\dim\ext^1_{kQ_3}(\ul{V}^{(3,5)},\ul{W}^i)=0$ $(1\leq i\leq 3)$,
\item $\dim\ext^1_{kQ_3}(\ul{V}^{(3,5)},\ul{V}^{(3,i)})=0$ $(0\leq i\leq
4)$.
\end{enumerate}
\end{prop}

Finally we include some computations of $\hom$ and $\ext$ groups
related to some other calculations with controlled quadratic
functors in the following section.

\begin{lem}\label{01}
$\ext^1_{\C{vect}(T_1)}(\um{M},\wedge^2_{T_1}\um{M})=0$ for any
elementary f. p. $\C{vect}(T_1)$-module $\um{M}$ in the sense of
Theorem \ref{clasclas}.
\end{lem}

\begin{proof}
This is true for $\um{A}$ and $\um{R}$ because they are projective.
For $\um{B}$ and $\um{B}_\infty$ it follows from Proposition
\ref{calculazo} and \cite{rca} 7.12. Since $\wedge^2_{T_1}$ is
right-exact, by Proposition \ref{calculazo} (1) and \cite{rca} 7.10
(1) we have that $\wedge^2_{T_1}\um{C}=0$, hence the case
$\um{M}=\um{C}$ follows. 
By \cite{rca} 7.10 (2) and 7.16 there is an epimorphism
$$\ext^1_{\C{vect}(T_1)}(\um{C}_\infty,\wedge^2_{T_1}\um{B}_\infty)\twoheadrightarrow\ext^1_{\C{vect}(T_1)}(\um{C}_\infty,\wedge^2_{T_1}\um{C}_\infty),$$
hence the case of $\um{C}_\infty$ follows from Proposition
\ref{calculazo} (2) and \cite{rca} 7.12.
\end{proof}

\begin{lem}\label{02}
The equality
$\ext^1_{\C{vect}(T_2)}(\um{M},\wedge^2_{T_2}\um{M})=0$ also holds
for any elementary f. p. $\C{vect}(T_2)$-module $\um{M}$ in the
sense of Theorem \ref{clasclas}.
\end{lem}

\begin{proof}
By Proposition \ref{esexacto} and Lemma \ref{01} it is enough to
check it for the elementary module $\uf{M}\ul{V}^{(2,1)}$, see
Theorem \ref{clasclas} and Proposition \ref{losn}, and it is easy
to prove by using Proposition \ref{diagiso1} that
$\wedge^2_{T_2}\uf{M}\ul{V}^{(2,1)}=0$.
\end{proof}

\begin{lem}\label{03}
For an elementary f. p. $\C{vect}(T_3)$-module $\um{M}$ in the
sense of Theorem \ref{clasclas} the equality
$\ext^1_{\C{vect}(T_3)}(\um{M},\wedge^2_{T_3}\um{M})=0$ is
satisfied unless $\um{M}=\uf{M}\ul{V}^{(3,5)}$. In this case
$$\ext^1_{\C{vect}(T_3)}(\uf{M}\ul{V}^{(3,5)},\wedge^2_{T_3}\uf{M}\ul{V}^{(3,5)})=\Z/2.$$
\end{lem}

\begin{proof}
By Proposition \ref{esexacto} and Lemma \ref{01} in order to prove
the first part of the statement it is enough to check the equality
$\ext^1_{\C{vect}(T_3)}(\um{M},\wedge^2_{T_3}\um{M})=0$ for the
elementary modules $\uf{M}\ul{V}^{(3,i)}$ for $1\leq i\leq 4$, see
Theorem \ref{clasclas} and Proposition \ref{losn}, and by using
Proposition \ref{diagiso1} one readily checks that
$\wedge^2_{T_3}\uf{M}\ul{V}^{(3,i)}=0$ in these cases.

The equality
$$\ext^1_{\C{vect}(T_3)}(\uf{M}\ul{V}^{(3,5)},\wedge^2_{T_3}\uf{M}\ul{V}^{(3,5)})=\Z/2$$
follows easily from Propositions \ref{diagiso1}, \ref{subext} and
\cite{rca} A.4.
\end{proof}

\begin{lem}\label{04}
For any elementary f. p. $\C{vect}(T_3)$-module
$\um{M}\neq\uf{F}^1_*\um{A}$, $\uf{F}^i_*\um{R}$, $\uf{M}\ul{V}^{(3,5)}$
$(1\leq i\leq 3)$ in the sense of Theorem \ref{clasclas}
$$\ext^1_{\C{vect}(T_3)}(\uf{M}\ul{V}^{(3,5)},\hat{\otimes}^2_{T_3}\um{M})=0.$$
Moreover if $\um{M}=\uf{F}^1_*\um{A}$ or $\uf{F}^i_*\um{R}$ then
$$\hom_{\C{vect}(T_3)}(\uf{M}\ul{V}^{(3,5)},\um{M})=0.$$
\end{lem}

\begin{proof}
If
$\um{M}\neq\uf{F}^i_*\um{B}_\infty,\uf{F}^i_*\um{C},\uf{F}^i_*\um{C}_\infty$
$(1\leq i\leq 3)$ the natural transformation
$\bar{\tau}_{T_3}\colon\um{M}\r\hat{\otimes}^2_{T_3}\um{M}$ is an
isomorphism. This can be easily checked by using for example
Proposition \ref{submod} and \ref{diagiso1}. Hence the lemma for
these modules follows from Proposition \ref{submod}, \ref{subhom}
and \ref{subext}. For $\um{M}=\uf{F}^i_*\um{C}$ $(1\leq i\leq 3)$
since $\hat{\otimes}^2_{T_3}$ is right-exact and all
$\C{vect}(T_3)$-modules have projective dimension $1$ by Proposition
\ref{ZaZ2}, then the epimorphism in \cite{rca} 7.10 (1) induces a
surjection
$$\ext^1_{\C{vect}(T_3)}(\uf{M}\ul{V}^{(3,5)},\hat{\otimes}^2_{T_3}\uf{F}^i_*\um{B})\twoheadrightarrow
\ext^1_{\C{vect}(T_3)}(\uf{M}\ul{V}^{(3,5)},\hat{\otimes}^2_{T_3}\uf{F}^i_*\um{C}),$$
hence this case follows from the case $\um{M}=\uf{F}^i_*\um{B}$
already checked.

The cases $\um{M}=\uf{F}^i_*\um{B}_\infty$ $(1\leq i\leq 3)$ are a
bit more complicated. Let $k=\uf{F}_2$ be the field with $2$
elements and let $$\varphi\colon (0\r k)\To (k\r k)$$ be the unique
non-trivial morphism of $1$-subspaces and
$$\uf{F}^i\colon\C{sub}^\mathrm{fin}_1\r\C{sub}^\mathrm{fin}_3,\;\;1\leq i\leq 3,$$
the functors considered in \cite{rca} 8.4, i. e.
$\uf{F}^i\ul{V}=\ul{W}$ with $W_0=V_0$, $W_i=V_1$ and $W_j=0$
otherwise. One can easily check that the cokernel of the
monomorphism
$$\phi=\left(\begin{array}{c}\uf{F}^1\varphi\\\uf{F}^2\varphi\\\uf{F}^3\varphi\end{array}\right)\colon \uf{F}^1(0\r k)=\uf{F}^2(0\r k)=\uf{F}^3(0\r k)\hookrightarrow
\begin{array}{c}\uf{F}^1(k\r k)\\\oplus \\\uf{F}^2(k\r k)\\\oplus\\\uf{F}^3(k\r k)\end{array}$$
is $\ul{V}^{(3,5)}$. By applying the functor $\uf{M}$ we obtain a
momomorphism of $\C{vect}(T_3)$-modules
$$\uf{M}\phi=\left(\begin{array}{c}\uf{M}\uf{F}^1\varphi\\\uf{M}\uf{F}^2\varphi\\\uf{M}\uf{F}^3\varphi\end{array}\right)\colon
\uf{F}^1_*\um{A}=\uf{F}^2_*\um{A}=\uf{F}^3_*\um{A}\hookrightarrow
\begin{array}{c}\uf{F}^1_*\um{B}\\\oplus
\\\uf{F}^2_*\um{B}\\\oplus\\\uf{F}^3_*\um{B}\end{array},$$ see Proposition \ref{submod}, whose
cokernel is $\uf{M}\ul{V}^{(3,5)}$. We are now going to prove that
$$\hom_{\C{vect}(T_3)}(\uf{M}\phi,\uf{F}^i_*\um{B}_\infty)$$ is an
isomorphism. By Propositions \ref{cambios} (1) and \ref{mashom} (2)
and the fact that $\uf{F}^i_*$ is a fully faithful left-adjoint of
$(\uf{F}^i)^*$ it is enough to see that
\begin{equation*}\tag{a}
\hom_{\C{vect}(T_1)}(\uf{M}\varphi,\um{B}_\infty)
\end{equation*}
is an isomorphism. Moreover, since $\varphi\neq 0$ by Proposition
\ref{mashom} (2) $\uf{M}\varphi\colon\um{A}\hookrightarrow\um{B}$
coincides with the monomorphism in \cite{rca} 7.10 (1), hence (a) is
an isomorphism by Proposition \ref{mashom} (4) and \cite{rca} 7.12.
Now the case $\um{M}=\uf{F}^i_*\um{B}_\infty$ follows from
Propositions \ref{mashom} (5) and \ref{calculazo} (7) by applying
the long exact sequence for the derived functors of
$\hom_{\C{vect}(T_3)}(-,\uf{F}^i_*\um{B}_\infty)$ to the short exact
sequence
$$\uf{F}^1_*\um{A}\st{\uf{M}\phi}\hookrightarrow\uf{F}^1_*\um{B}\oplus\uf{F}^2_*\um{B}\oplus\uf{F}^3_*\um{B}
\twoheadrightarrow\uf{M}\ul{V}^{(3,5)}.$$ Moreover, this also
proves that
$$\hom_{\C{vect}(T_3)}(\uf{M}\ul{V}^{(3,5)},\uf{F}^i_*\um{B}_\infty)=0,$$
therefore the case $\um{M}=\uf{F}^i_*\um{R}$ follows from
Proposition \ref{esexacto} and the existence of a monomorphism
$\um{R}\hookrightarrow\um{B}_\infty$, see \cite{rca} 7.10 (2).

For $\um{M}=\uf{F}^i_*\um{C}_\infty$ $(1\leq i\leq 3)$ by
Proposition \ref{ZaZ2} and \cite{rca} 7.10 (2) there is an
epimorphism
$$\ext^1_{\C{vect}(T_3)}(\uf{M}\ul{V}^{(3,5)},\hat{\otimes}^2_{T_3}\uf{F}^i_*\um{B}_\infty)\twoheadrightarrow
\ext^1_{\C{vect}(T_3)}(\uf{M}\ul{V}^{(3,5)},\hat{\otimes}^2_{T_3}\uf{F}^i_*\um{C}_\infty),$$
so it follows from the case $\um{M}=\uf{F}^i_*\um{B}_\infty$ already
proved.
\end{proof}

\section{Some computations with the controlled quadratic
functors}\label{quadcomp}

In this section the ground field will be $k=\uf{F}_2$. One of the
main results of this section is the following theorem.

\begin{thm}\label{split}
If $T$ has less than four ends the morphism
$\bar{\tau}_T\colon\um{M}\otimes\Z/2\r\hat{\otimes}^2_T\um{M}$ is
a split monomorphism for any f. p. $\ab(T)$-module $\um{M}$.
\end{thm}

\begin{proof}
We will see that it is enough to prove the theorem for
$\C{vect}(T)$-modules, hence it will follow from Proposition
\ref{splitmod2} below.

Since $-\otimes\Z/2$ is additive and right-exact $\um{M}\otimes\Z/2$
is f. p. as a $\C{vect}(T)$-module. The natural projection
$\hat{p}\colon\um{M}\twoheadrightarrow\um{M}\otimes\Z/2$ induces a
commutative diagram
$$\xymatrix{\um{M}\otimes\Z/2\ar[d]_{\hat{p}\otimes\Z/2}^\simeq\ar[r]^{\bar{\tau}_T}&\hat{\otimes}^2\um{M}\ar[d]^{\hat{\otimes}^2_T\hat{p}}\\
(\um{M}\otimes\Z/2)\otimes\Z/2\ar[r]^{\bar{\tau}_T}&\hat{\otimes}^2_T(\um{M}\otimes\Z/2)}$$
If $s$ is a left-inverse of the lower horizontal row then
$(\hat{p}\otimes\Z/2)^{-1}s(\hat{\otimes}^2_T\hat{p})$ is a
left-inverse of the upper one.
\end{proof}

All endofunctors of the category of finite-dimensional vector spaces
extend to $n$-subspaces in the obvious way. The same happens with
natural transformations, in particular diagram (\ref{d1}) tensored
by $\Z/2$ gives rise to a natural diagram with exact rows and column
in $\C{sub}_n^\mathrm{fin}$,
\begin{equation*}
\xymatrix{\otimes^2\ul{V}\ar[d]_{[-,-]\otimes\Z/2}&&\\(\Gamma
\ul{V})\otimes\Z/2\ar[r]^<(.35){\tau\otimes\Z/2}\ar@{->>}[d]_{\sigma\otimes\Z/2}\ar@{}[dr]|{\mathrm{push}}&\otimes^2\ul{V}\ar@{->>}[d]^{\bar{\sigma}}\ar@{->>}[r]^q&\wedge^2\ul{V}\ar@{=}[d]\\
\ul{V}\ar@{^{(}->}[r]^{\bar{\tau}}&\hat{\otimes}^2\ul{V}\ar@{->>}[r]^{\bar{q}}&\wedge^2\ul{V}}
\end{equation*}
The category $\C{sub}_n^\mathrm{fin}$ is not abelian, but it is
fully included into the abelian category of representations of
$Q_n$. It is in this last category where the previous sequences of
$n$-subspaces are exact.

As it happened with (\ref{d1}) and (\ref{d2}) this diagram is
completely determined by the functor $(\Gamma-)\otimes\Z/2$.

\begin{prop}\label{diagiso1}
The following diagram of $\C{vect}(T_n)$-modules
$$\xymatrix@C=40pt{\uf{M}\otimes^2\ul{V}\ar[d]_{\uf{M}([-,-]\otimes\Z/2)}&&\\\uf{M}((\Gamma
\ul{V})\otimes\Z/2)\ar@{}[dr]|{\mathrm{push}}\ar@{->>}[d]_{\uf{M}(\sigma\otimes\Z/2)}\ar[r]^<(.3){\uf{M}(\tau\otimes\Z/2)}&
\uf{M}\otimes^2\ul{V}\ar@{->>}[r]^{\uf{M}q}\ar@{->>}[d]_{\uf{M}\bar{\sigma}}&\uf{M}\wedge^2\ul{V}\ar@{=}[d]\\
\uf{M}\ul{V}\ar@{^{(}->}[r]^<(.3){\uf{M}\bar{\tau}}&\uf{M}\hat{\otimes}^2\ul{V}\ar@{->>}[r]^{\uf{M}\bar{q}}&\uf{M}\wedge^2\ul{V}}$$
is naturally isomorphic to
$$\xymatrix@C=40pt{\otimes^2_{T_n}\uf{M}\ul{V}\ar[d]_{
[-,-]_{T_n}\otimes\Z/2}&&\\(\Gamma_{T_n}\uf{M}\ul{V})\otimes\Z/2\ar@{}[dr]|{\mathrm{push}}\ar@{->>}[d]_{\sigma_{T_n}\otimes\Z/2}\ar[r]^<(.3){\tau_{T_n}\otimes\Z/2}&
\otimes^2_{T_n}\uf{M}\ul{V}\ar@{->>}[r]^{q_{T_n}}\ar@{->>}[d]_{\bar{\sigma}_{T_n}}&\wedge^2_{T_n}\uf{M}\ul{V}\ar@{=}[d]\\
\uf{M}\ul{V}\ar[r]^<(.3){\bar{\tau}_{T_n}}&\hat{\otimes}^2_{T_n}\uf{M}\ul{V}\ar@{->>}[r]^{\bar{q}_{T_n}}&\wedge^2_{T_n}\uf{M}\ul{V}}$$
\end{prop}

\begin{proof}
Since the first diagram is determined by
$\uf{M}((\Gamma-)\otimes\Z/2)$ and the second one by
$(\Gamma_{T_n}-)\otimes\Z/2$ it will be enough to construct a
natural isomorphism
$$(\Gamma_{T_n}\uf{M}\ul{V})\otimes\Z/2\simeq\uf{M}((\Gamma
\ul{V})\otimes\Z/2).$$

In \cite{rca} 9.4 we constructed a finite presentation
$$\uf{F}_2\grupo{D}_\delta\st{\rho}\To\uf{F}_2\grupo{C}_\gamma\st{p}\twoheadrightarrow\uf{M}\ul{V}.$$
Here $p$ is defined by a certain vector space homomorphism
$p_0\colon\uf{F}_2\grupo{C}\twoheadrightarrow V_0$, see the
definition of $\uf{M}$ in \cite{rca} (9.a).
By right-exactness we have
an exact sequence
\begin{small}
$$(\Gamma_{T_n}\uf{F}_2\grupo{D}_\delta)\otimes\Z/2\oplus\uf{F}_2\grupo{D}_\delta\otimes_{T_n}
\uf{F}_2\grupo{C}_\gamma\st{\xi}\To (\Gamma_{T_n}\uf{F}_2\grupo{
C}_\gamma)\otimes\Z/2\st{(\Gamma_{T_n}p)\otimes\Z/2}\twoheadrightarrow(\Gamma_{T_n}\uf{M}\ul{V})\otimes\Z/2$$
\end{small}
with
$\xi=((\Gamma_{T_n}\rho)\otimes\Z/2,([-,-]_{T_n}\otimes\Z/2)(\rho\otimes_{T_n}1))$.

We leave to the reader to check that the vector space homomorphism
$(\Gamma
p_0)\otimes\Z/2\colon(\Gamma\uf{F}_2\grupo{C})\otimes\Z/2\r(\Gamma
V_0)\otimes\Z/2$ determines a morphism $p'\colon
(\Gamma_{T_n}\uf{F}_2\grupo{
C}_\gamma)\otimes\Z/2\r\uf{M}((\Gamma\ul{V})\otimes\Z/2)$ such
that the sequence
$$(\Gamma_{T_n}\uf{F}_2\grupo{D}_\delta)\otimes\Z/2\oplus\uf{F}_2\grupo{D}_\delta\otimes_{T_n}
\uf{F}_2\grupo{C}_\gamma\st{\xi}\To (\Gamma_{T_n}\uf{F}_2\grupo{
C}_\gamma)\otimes\Z/2\st{p'}\twoheadrightarrow\uf{M}((\Gamma\ul{V})\otimes\Z/2)$$
is also exact, hence we get the desired isomorphism.

\end{proof}

Recall from \cite{rca} Section 4 that there is an exact, full and
faithful functor from the category of vector spaces with countable
dimension
$$\mathfrak{i}\colon\mo_{\aleph_0}(\uf{F}_2)\To\C{fp}(\C{vect}(\Real_+)).$$

\begin{prop}\label{diagiso2}
Let $V$ be a vector space with $\dim V\leq\aleph_0$, the following
diagram of $\C{vect}(\Real_+)$-modules
$$\xymatrix@C=40pt{\mathfrak{i}\otimes^2{V}\ar[d]_{\mathfrak{i}([-,-]\otimes\Z/2)}&&\\\mathfrak{i}((\Gamma
{V})\otimes\Z/2)\ar@{}[dr]|{\mathrm{push}}\ar@{->>}[d]_{\mathfrak{i}(\sigma\otimes\Z/2)}\ar[r]^<(.3){\mathfrak{i}(\tau\otimes\Z/2)}&
\mathfrak{i}\otimes^2{V}\ar@{->>}[r]^{\mathfrak{i}q}\ar@{->>}[d]_{\mathfrak{i}\bar{\sigma}}&\mathfrak{i}\wedge^2{V}\ar@{=}[d]\\
\mathfrak{i}{V}\ar@{^{(}->}[r]^<(.3){\mathfrak{i}\bar{\tau}}&\mathfrak{i}\hat{\otimes}^2{V}\ar@{->>}[r]^{\mathfrak{i}\bar{q}}&\mathfrak{i}\wedge^2{V}}$$
is naturally isomorphic to
$$\xymatrix@C=40pt{\otimes^2_{\Real_+}\mathfrak{i}{V}\ar[d]_{
[-,-]_{\Real_+}\otimes\Z/2}&&\\(\Gamma_{\Real_+}\mathfrak{i}{V})\otimes\Z/2\ar@{}[dr]|{\mathrm{push}}\ar@{->>}[d]_{\sigma_{\Real_+}\otimes\Z/2}\ar[r]^<(.3){\tau_{\Real_+}\otimes\Z/2}&
\otimes^2_{\Real_+}\mathfrak{i}{V}\ar@{->>}[r]^{q_{\Real_+}}\ar@{->>}[d]_{\bar{\sigma}_{\Real_+}}&\wedge^2_{\Real_+}\mathfrak{i}{V}\ar@{=}[d]\\
\mathfrak{i}{V}\ar[r]^<(.3){\bar{\tau}_{\Real_+}}&\hat{\otimes}^2_{\Real_+}\mathfrak{i}{V}\ar@{->>}[r]^{\bar{q}_{\Real_+}}&\wedge^2_{\Real_+}\mathfrak{i}{V}}$$
\end{prop}

\begin{proof}
As in the proof of Proposition \ref{diagiso1} it will de enough to
construct a natural isomorphism
$$(\Gamma_{\Real_+}\mathfrak{i}{V})\otimes\Z/2\simeq\mathfrak{i}((\Gamma
{V})\otimes\Z/2).$$

An explicit finite presentation
$$\uf{F}_2\grupo{D}_\delta\st{\rho}\To\uf{F}_2\grupo{C}_\gamma\st{p}\twoheadrightarrow\mathfrak{i}{V}$$
can be obtained from the proof of \cite{rca} 4.3. Here $p$ is
determined by a vector space homomorphism
$p_0\colon\uf{F}_2\grupo{C}\r V$. By right exactness the following
sequence is exact
\begin{small}
$$(\Gamma_{\Real_+}\uf{F}_2\grupo{D}_\delta)\otimes\Z/2\oplus\uf{F}_2\grupo{D}_\delta\otimes_{\Real_+}
\uf{F}_2\grupo{C}_\gamma\st{\xi}\To
(\Gamma_{\Real_+}\uf{F}_2\grupo{
C}_\gamma)\otimes\Z/2\st{(\Gamma_{\Real_+}p)\otimes\Z/2}\twoheadrightarrow(\Gamma_{\Real_+}\mathfrak{i}{V})\otimes\Z/2,$$
\end{small}
where
$\xi=((\Gamma_{\Real_+}\rho)\otimes\Z/2,([-,-]_{\Real_+}\otimes\Z/2)(\rho\otimes_{\Real_+}1))$.

Now the reader can check that the morphism
$p'\colon(\Gamma_{\Real_+}\uf{F}_2\grupo{
C}_\gamma)\otimes\Z/2\rightarrow\mathfrak{i}((\Gamma
{V})\otimes\Z/2)$ induced by the vector space homomorphism
$(\Gamma
p_0)\otimes\Z/2\colon(\Gamma\uf{F}_2\grupo{C})\otimes\Z/2\r(\Gamma
{V})\otimes\Z/2$ fits into an exact sequence
$$(\Gamma_{\Real_+}\uf{F}_2\grupo{D}_\delta)\otimes\Z/2\oplus\uf{F}_2\grupo{D}_\delta\otimes_{\Real_+}
\uf{F}_2\grupo{C}_\gamma\st{\xi}\To
(\Gamma_{\Real_+}\uf{F}_2\grupo{
C}_\gamma)\otimes\Z/2\st{p'}\twoheadrightarrow\mathfrak{i}((\Gamma
{V})\otimes\Z/2).$$ This determines the desired isomorphism.
\end{proof}

\begin{prop}\label{splitcinf}
The morphism
$\bar{\tau}_{\Real_+}\colon\um{C}_\infty\otimes\Z/2\r\hat{\otimes}^2_{\Real_+}\um{C}_\infty$
is indeed a split monomorphism.
\end{prop}

\begin{proof}
Consider the free $\Real_+$-controlled $\uf{F}_2$-module
$\uf{F}_2\grupo{A}_\alpha$ with $$A=\set{a_{mn}\,;\,0<n\leq m}$$ and
$\alpha(a_{mn})=m$, and its endomorphism $\varphi$ defined as

$$\varphi(a_{mn})=\left\{%
\begin{array}{ll}
    a_{mm}, & \hbox{if $n=m$;} \\
    a_{mn}+a_{m-1,n}, & \hbox{if $n<m$.} \\
\end{array}%
\right.    $$ By \cite{rca} 7.21 $\um{C}_\infty=\coker\varphi$.
Moreover by right exactness we have that
$\hat{\otimes}^2_{\Real_+}\um{C}_\infty$ is the cokernel of the
controlled homomorphism
$$(\hat{\otimes}^2_{\Real_+}\varphi,(\hat{\otimes}^2_{\Real_+}(\varphi|1))i_{12})\colon\uf{F}_2\grupo{\hat{\otimes}^2A}_{\hat{\otimes}^2\alpha}\oplus\uf{F}_2\grupo{A\times
A}_{\alpha\otimes\alpha}
\To\uf{F}_2\grupo{\hat{\otimes}^2A}_{\hat{\otimes}^2\alpha},$$ see
Proposition \ref{redten}. We
consider the following total order in $A$ for the definition of $\hat{\otimes}^2A$ $$a_{mn}\preceq a_{pq}\Leftrightarrow\left\{%
\begin{array}{ll}
    m<p, \\
    \hbox{or} \\
    m=p \hbox{ and } n\leq q. \\
\end{array}%
\right.    $$

There is a commutative diagram of solid arrows
$$\xymatrix{\uf{F}_2\grupo{A}_\alpha\ar[d]_\varphi\ar[r]^<(.13){(\bar{\tau}_{\Real_+},0)}&\uf{F}_2\grupo{\hat{\otimes}^2A}_{\hat{\otimes}^2\alpha}\oplus\uf{F}_2\grupo{A\times
A}_{\alpha\otimes\alpha}\ar[d]^{(\hat{\otimes}^2_{\Real_+}\varphi,(\hat{\otimes}^2_{\Real_+}(\varphi|1))i_{12})}\ar@{-->}@/^15pt/[l]^{s_1}\\
\uf{F}_2\grupo{A}_\alpha\ar[r]^<(.33){\bar{\tau}_{\Real_+}}&\uf{F}_2\grupo{\hat{\otimes}^2A}_{\hat{\otimes}^2\alpha}
\ar@{-->}@/^10pt/[l]^{s_0}}$$ One can check that the controlled
homomorphisms $s_1$ and $s_0$ given by $$s_1(a_{mn}\hat{\otimes}a_{pq})=\left\{%
\begin{array}{ll}
    a_{mn}, & \hbox{if $p=m\geq n=q$;} \\
    0, & \hbox{otherwise;} \\
\end{array}%
\right.    $$
$$s_0(a_{mn}\hat{\otimes}a_{pq})=s_1(a_{mn},a_{pq})=\left\{%
\begin{array}{ll}
    a_{mn}, & \hbox{if $p\geq m\geq n=q$;} \\
    0, & \hbox{otherwise;} \\
\end{array}%
\right.    $$ are retractions of the horizontal solid arrows and
satisfy $$\varphi
s_1=s_0(\hat{\otimes}^2_{\Real_+}\varphi,(\hat{\otimes}^2_{\Real_+}(\varphi|1))i_{12}),$$
hence we obtain the desired retraction
$\hat{\otimes}^2_{\Real_+}\um{C}_\infty\twoheadrightarrow\um{C}_\infty\otimes\Z/2$
by taking cokernels.
\end{proof}

The following corollary follows from the fact that $\um{C}$ is a
retract of $\um{C}_\infty$, see Theorem \ref{clasclas}.

\begin{cor}\label{splitc}
The morphism
$\bar{\tau}_{\Real_+}\colon\um{C}\otimes\Z/2\r\hat{\otimes}^2_{\Real_+}\um{C}$
is a split monomorphism.
\end{cor}

The following proposition is a consequence of the right exactness of
$-\otimes\Z/2$, Propositions \ref{quadpropri} and \ref{redten} and
Lemma \ref{conZ2}.

\begin{prop}\label{splitfree}
The morphism
$\bar{\tau}_{T}\colon\uf{F}_2\grupo{A}_\alpha\r\hat{\otimes}^2_T\uf{F}_2\grupo{A}_\alpha$
is indeed a split monomorphism for any free $T$-controlled
$\uf{F}_2$-module $\uf{F}_2\grupo{A}_\alpha$.
\end{prop}

\begin{prop}\label{splitmod2}
Given a tree $T$ with $\card\F(T)\leq 3$ the morphism
$\bar{\tau}_T\colon\um{M}\r\hat{\otimes}^2_T\um{M}$ is a split
monomorphism for any f. p. $\C{vect}(T)$-module $\um{M}$.
\end{prop}

\begin{proof}
It is enough to prove the proposition for elementary
$\C{vect}(T)$-modules. For trees with less than four ends a
complete list of elementary f. p. $\C{vect}(T)$-modules can be
found in Theorem \ref{clasclas}, see also Proposition \ref{losn}.

If $T$ has one end the proposition follows from \ref{splitfree} and
\cite{rca} 7.7 for the elementary modules $\um{A}$ and $\um{R}$,
from Proposition \ref{diagiso2} and \cite{rca} 7.11 for $\um{B}$ and
$\um{B}_\infty$, from Proposition \ref{splitcinf} for
$\um{C}_\infty$, and from Corollary \ref{splitc} for $\um{C}$.

For trees with $n= 2$ or $3$ ends by Proposition \ref{redten} it is
only left to prove this proposition for elementary f. p. modules
coming from rigid $n$-subspaces, but this is easy to check by using
Proposition \ref{diagiso1} and \ref{losn}.
\end{proof}

\begin{rem}\label{herculillo}\renewcommand{\theequation}{\alph{equation}}\setcounter{equation}{0}
Proposition \ref{splitmod2} does hot hold for trees with more than
$6$ ends. Indeed by Propositions \ref{redten} and \ref{diagiso1} it
will be enough to show a $7$-subspace $\ul{V}$ such that
$\bar{\tau}\colon\ul{V}\hookrightarrow\hat{\otimes}^2\ul{V}$ does
not split. We have found the following example,
\begin{eqnarray*}
  V_0&=&\uf{F}_2\grupo{x,y,z}\\
  V_1 &=& \uf{F}_2\grupo{x,y}, \\
  V_2 &=& \uf{F}_2\grupo{x,z}, \\
  V_3 &=& \uf{F}_2\grupo{y,z}, \\
  V_4 &=& \uf{F}_2\grupo{x+y,z}, \\
  V_5 &=& \uf{F}_2\grupo{x+z,y}, \\
  V_6 &=& \uf{F}_2\grupo{x,y+z}, \\
  V_7 &=& \uf{F}_2\grupo{x+y,x+z}.
\end{eqnarray*}
The $7$-subspace $\hat{\otimes}^2\ul{V}$ is given by
\begin{eqnarray*}
\hat{\otimes}^2V_0&=&\uf{F}_2\grupo{x\hat{\otimes} x,
x\hat{\otimes} y,x\hat{\otimes} z,y\hat{\otimes} y, y\hat{\otimes}
z,
z\hat{\otimes} z}\\
 \hat{\otimes}^2 V_1 &=& \uf{F}_2\grupo{x\hat{\otimes} x, x\hat{\otimes} y,y\hat{\otimes}y}, \\
\hat{\otimes}^2  V_2 &=& \uf{F}_2\grupo{x\hat{\otimes}x, x\hat{\otimes} z,z\hat{\otimes}z}, \\
\hat{\otimes}^2  V_3 &=& \uf{F}_2\grupo{y\hat{\otimes}y, y\hat{\otimes} z,z\hat{\otimes}z}, \\
\hat{\otimes}^2  V_4 &=& \uf{F}_2\grupo{x\hat{\otimes}x+y\hat{\otimes}y,x\hat{\otimes}z+y\hat{\otimes}z,z\hat{\otimes}z}, \\
\hat{\otimes}^2  V_5 &=& \uf{F}_2\grupo{x\hat{\otimes}x+z\hat{\otimes}z,x\hat{\otimes}y+y\hat{\otimes}z,y\hat{\otimes}y}, \\
\hat{\otimes}^2  V_6 &=& \uf{F}_2\grupo{x\hat{\otimes}x,x\hat{\otimes}y+x\hat{\otimes}z,y\hat{\otimes}y+z\hat{\otimes}z}, \\
\hat{\otimes}^2  V_7 &=&
\uf{F}_2\grupo{x\hat{\otimes}x+x\hat{\otimes}y+x\hat{\otimes}z+y\hat{\otimes}z,x\hat{\otimes}x+y\hat{\otimes}y,x\hat{\otimes}x+z\hat{\otimes}z}.
\end{eqnarray*}
If there where a retraction
$s\colon\hat{\otimes}^2\ul{V}\twoheadrightarrow\ul{V}$ of
$\bar{\tau}$ it should satisfy
\begin{eqnarray}
  s(x\hat{\otimes}x) &=& x, \\
\nonumber  s(y\hat{\otimes}y) &=& y, \\
\nonumber s(z\hat{\otimes}z) &=& z, \\
  s(x\hat{\otimes}y) &\in& V_1, \\
  s(x\hat{\otimes}z) &\in& V_2, \\
  s(y\hat{\otimes}z) &\in& V_3, \\
  s(x\hat{\otimes}z+y\hat{\otimes}z) &\in& V_4, \\
  s(x\hat{\otimes}y+y\hat{\otimes}z) &\in& V_5,\\
  s(x\hat{\otimes}y+x\hat{\otimes}z) &\in& V_6,\\
    s(x\hat{\otimes}x+x\hat{\otimes}y+x\hat{\otimes}z+y\hat{\otimes}z) &\in&
    V_7.
\end{eqnarray}
By using (b), (c) and (d) we see that there would exist
$a,b,c,d,e,f\in\uf{F}_2$ such that
\begin{eqnarray*}
  s(x\hat{\otimes}y) &=& ax+by,\\
  s(x\hat{\otimes}z) &=& cx+dz,\\
  s(y\hat{\otimes}z) &=& ey+fz.
\end{eqnarray*}
By combining this with (e), (f) and (g) we obtain that the
following equalities should hold
\begin{eqnarray*}
  e &=& c, \\
  f &=& a, \\
  d &=& b.
\end{eqnarray*}
Finally (h) together with (a) and the previous equalities would
imply that
$$(1+a+c)x+(b+c)y+(a+b)z\in V_7=\set{0,x+y,x+z,y+z}.$$ That is,
at least one of the following systems of linear equations with
three indeterminacies $(a,b,c)$ should have a solution over
$\uf{F}_2$,
$$\begin{array}{lcr}
\left\{%
\begin{array}{rcc}
     1+a+c&=&0,\\
     b+c&=&0,\\
     a+b&=&0,
\end{array}%
\right.&\;\;\;\;&
\left\{%
\begin{array}{rcc}
     1+a+c&=&1,\\
     b+c&=&1,\\
     a+b&=&0,
\end{array}%
\right.
\end{array}$$
$$\begin{array}{lcr}\left\{%
\begin{array}{rcc}
     1+a+c&=&1,\\
     b+c&=&0,\\
     a+b&=&1,
\end{array}%
\right.&\;\;\;\;&
\left\{%
\begin{array}{rcc}
     1+a+c&=&0,\\
     b+c&=&1,\\
     a+b&=&1.
\end{array}%
\right. \end{array}$$ It is easy to see that all these systems are
incompatible, therefore we would reach a contradiction, so the
retraction can not exist.

We do not know whether Proposition \ref{splitmod2} is true
for trees with $4$, $5$ or $6$ ends. One could try to check it for
trees with $4$ ends by using Theorem \ref{clasclas} and Nazarova's
classification of finite-dimensional $4$-subspaces in
\cite{nazarova}. We have not carried out such a computation
because the infinite list of indecomposable $4$-subspaces over
$\uf{F}_2$ is quite complicate to describe and we have already
obtained very satisfactory results for trees with less than $4$
ends. For trees $T$ with $5$ or $6$ ends there is not a direct
method to check wether Proposition \ref{splitmod2} is true since
the category of finitely presented $\C{vect}(T)$-modules has wild
representation type, see Theorem \ref{rep}.
\end{rem}

We conclude this section with explicit computations of values of
the controlled quadratic functors over elementary f. p.
$\C{vect}(T)$-modules.

\begin{prop}\label{calculazo}
We have the following identifications
\begin{enumerate}
\item $\wedge^2_{T_1}\um{B}=0$,
\item $\wedge^2_{T_1}\um{B}_\infty=\um{B}_\infty$,
\item $\wedge^2_{T_2}\uf{M}\ul{V}^{(2,1)}=0$,
\item $\wedge^2_{T_3}\uf{M}\ul{V}^{(3,i)}=0$ for $i=1,2,3$ and
$4$,
\item $\wedge^2_{T_3}\uf{M}\ul{V}^{(3,5)}=\um{A}$,
\item $\hat{\otimes}^2_{T_1}\um{B}=\um{B}$,
\item $\hat{\otimes}^2_{T_1}\um{B}_\infty=\um{B}_\infty$,
\item $\hat{\otimes}^2_{T_2}\uf{M}\ul{V}^{(2,1)}=\uf{M}\ul{V}^{(2,1)}$,
\item $\hat{\otimes}^2_{T_3}\uf{M}\ul{V}^{(3,i)}=\uf{M}\ul{V}^{(3,i)}$ for $i=1,2,3$ and
$4$.
\end{enumerate}
\end{prop}

\begin{proof}
The identifications (1), (2), (6) and (7) follow from Proposition
\ref{diagiso2} and \cite{rca} 7.11; (3), (4), (8) and (9) can be
easily derived from Proposition \ref{diagiso1}. Moreover, if
$\ul{W}$ is the $3$-subspace with $W_0=\uf{F}_2$ and $W_i=0$ for
$1\leq i\leq 3$ one readily checks by using Proposition
\ref{diagiso1} that
$\wedge^2_{T_3}\uf{M}\ul{V}^{(3,5)}=\uf{M}\ul{W}$. The finite
presentation of $\uf{M}\ul{W}$ as a $\C{vect}(T_3)$-module
constructed in the proof of \cite{rca} 9.1 shows that $\uf{M}\ul{W}$
is a free $T_3$-controlled $\uf{F}_2$-module with only one
generator, hence by \cite{rca} 7.7 $\uf{M}\ul{W}=\um{A}$ and (5)
holds.
\end{proof}



\appendix
\chapter{Proof of Proposition
{\ref{extigrande}}}\label{proofofextigrande}

For an arbitrary tree $T$ we will consider the following
commutative diagram of functors
$$\xymatrix{\C{M}_R(T)\ar[rr]^{Y}\ar[rd]_{J}&&\mo(\C{M}_R(T))\\&\C{M}^b_R(T)\ar[ru]_{Y'}&}$$
Here $Y$ is the Yoneda full inclusion in (\ref{yoneda}), $J$ is
the inclusion of the full subcategory, and $Y'$ is the faithful
extension of $Y$ along $J$ in (\ref{yonedabig}).

The functors $J$ and $Y$ are natural with respect to ``change of
tree'' functors induced by a proper map $f\colon T\rightarrow T'$,
i. e. the following two diagrams commute
\begin{equation}\label{2diag}
\begin{array}{cc}
\xymatrix{\C{M}_R(T)\ar[r]^{\uf{F}^f}\ar@{^{(}->}[d]_{J}&
\C{M}_R(T')\ar@{^{(}->}[d]^{J}\\
\C{M}^b_R(T)\ar[r]^{\uf{F}^f}&\C{M}^b_R(T')}
&\xymatrix{\C{M}_R(T)\ar[r]^{\uf{F}^f}\ar@{^{(}->}[d]_{Y}&
\C{M}_R(T')\ar@{^{(}->}[d]^{Y}\\
\mo(\C{M}_R(T))\ar[r]^{\uf{F}^f_*}&\mo(\C{M}_R(T'))}
\end{array}
\end{equation}
The first one commutes by the very definition of $\uf{F}^f$ in
Section \ref{freecontrol}, the second one is a particular example
of (\ref{yonedacom}). Proposition \ref{extigrande} establishes
that $Y'$ is also natural with respect to ``change of tree''
functors.

In this appendix we will make an extensive use of comma
categories, see \cite{cwm} II.6 for basic definitions. In
particular given a functor $F\colon\C{C}\rightarrow\C{D}$ we shall
consider the comma categories $c\downarrow F$ and $F\downarrow d$,
where $c$ and $d$ are objects of $\C{C}$ and $\C{D}$,
respectively. There is a canonical functor $P\colon F\downarrow
d\rightarrow \C{C}$. If we have a commutative diagram of functors
$$\xymatrix{\C{C}\ar[d]_F\ar[r]^M&\C{C}'\ar[d]^G\\\C{D}\ar[r]^N&\C{D}'}$$
there is also an induced functor $M\downarrow N\colon F\downarrow
d\rightarrow G\downarrow N(d)$ which satisfies
\begin{equation}\label{masuna}
P(M\downarrow N)=MP\colon F\downarrow d\rightarrow\C{D}.
\end{equation}

It is well-known that for any $\C{M}_R(T)$-module $\um{M}$
\begin{equation}\label{natural-es}
\colim [YP\colon
Y\downarrow\um{M}\rightarrow\mo(\C{M}_R(T))]=\um{M}.
\end{equation}

The crucial step in the proof of Proposition \ref{extigrande} is
the next lemma.

\begin{lem}\label{final-es}
For any object $R\grupo{A}_\alpha$ in $\C{M}^b_R(T)$ the functor
$$\uf{F}^f\downarrow \uf{F}^f\colon J\downarrow
R\grupo{A}_\alpha\rightarrow J\downarrow
\uf{F}^fR\grupo{A}_\alpha$$ is final.
\end{lem}

See the first definition of \cite{cwm} IX.3 for a definition of
\emph{final} functor. In the proof of Lemma \ref{final-es} we
shall use the following technical result.

\begin{lem}\label{tech}
If  $R\grupo{A}_\alpha$ is an object of $\C{M}^b_R(T)$ and
$R\grupo{B'}_{\beta'}$ is an object of $\C{M}_R(T')$, given a
controlled homomorphism $\varphi\colon
R\grupo{B'}_{\beta'}\r\uf{F}^fR\grupo{A}_\alpha$ there exists a
subset $\bar{A}\subset A$ such that if $\bar{\alpha}$ is the
restriction of $\alpha$ to $\bar{A}$ then
$R\grupo{\bar{A}}_{\bar{\alpha}}$ is an object of $\C{M}_R(T)$ and
if $i\colon R\grupo{\bar{A}}_{\bar{\alpha}}\hookrightarrow
R\grupo{A}_\alpha$ is the inclusion, which is obviously splitting,
$\varphi$ factors through the monomorphism $\uf{F}^fi$.
\end{lem}

\begin{proof}
If we define $\bar{A}\subset A$ as the least subset such that
$\varphi(B')\subset R\grupo{\bar{A}}$ it is enough to check that the
object $R\grupo{\bar{A}}_{\bar{\alpha}}$ belongs to $\C{M}_R(T)$.
Suppose by the contrary that there exists a compact subset $K\subset
T$ such that $\bar{\alpha}^{-1}(K)$ contains an infinite countable
subset $\set{a_n}_{n\geq 0}\subset A$, then we can take a sequence
$\set{b_n}_{n\geq 0}\subset B'$ such that $a_n$ appears in the
linear expansion of $\varphi(b_n)$ with non-trivial coefficient
$(n\geq 0)$. Since $\beta'\colon B'\r T'$ is proper we can suppose
without loss of generality, taking subsequences if necessary, that
$\lim\limits_{n\r\infty}\beta'(b_n)=\varepsilon\in\F(T')$. By using
the definition of controlled homomorphism, see Section
\ref{freecontrol}, it is easy to see that this implies that
$\lim\limits_{n\r\infty}f\alpha(a_n)=\varepsilon\in\F(T')$, but for
all $n\geq 0$ $f\alpha(a_n)\subset f(K)\in T'$ compact, and we reach
a contradiction.
\end{proof}

\begin{proof}[Proof of Lemma \ref{final-es}]
An object of $J\downarrow \uf{F}^fR\grupo{A}_\alpha$ is a morphism
$$\varphi\colon R\grupo{B'}_{\beta'}\r\uf{F}^fR\grupo{A}_\alpha$$ as
in the statement of Lemma \ref{tech}. We have to check that the
category $\varphi\downarrow(\uf{F}^f\downarrow\uf{F}^f)$ is
non-empty and connected. It is non-empty since by Lemma \ref{tech}
there is a commutative diagram
$$\xymatrix{R\grupo{B'}_{\beta'}\ar[rd]_\varphi\ar[rr]&&\uf{F}^fR\grupo{\bar{A}}_{\bar{\alpha}}
\ar@{^{(}->}[ld]^{\uf{F}^fi}\\
&\uf{F}^fR\grupo{A}_\alpha&}$$ where $i\colon
R\grupo{\bar{A}}_{\bar{\alpha}}\hookrightarrow R\grupo{A}_\alpha$
is an object of $J\downarrow R\grupo{A}_\alpha.$

Let us check now that
$\varphi\downarrow(\uf{F}^f\downarrow\uf{F}^f)$ is connected. Two
objects $a$ and $b$ in this category are exactly two commutative
triangles of controlled homomorphisms as follows
$$\xymatrix{&R\grupo{B'}_{\beta'}\ar[ld]_{\tau_1}\ar[d]^\varphi\ar[rd]^{\tau_2}&\\
\uf{F}^fR\grupo{B_1}_{\beta_1}\ar[r]_{{\uf{F}}^f\psi_1}&\uf{F}^fR\grupo{A}_\alpha&
\uf{F}^fR\grupo{B_2}_{\beta_2}\ar[l]^{{\uf{F}}^f\psi_2}}$$ Here
$R\grupo{B_1}_{\beta_1}$ and $R\grupo{B_2}_{\beta_2}$ are not big.
Since $\uf{F}^f$ is additive this diagram can be extended
commutatively in the following
way $$\xymatrix{&&&&\\&\uf{F}^fR\grupo{A}_\alpha&&\uf{F}^fR\grupo{B_1}_{\beta_1}\ar[ll]^{{\uf{F}}^f\psi_1}&\\
&&R\grupo{B'}_{\beta'}\ar[lu]_\varphi\ar[ru]^{\tau_1}
\ar[rd]^<(.25){\scriptsize{\left(\begin{array}{c}\tau_1\\\tau_2\end{array}\right)}}
\ar[ld]_{\tau_2}&&\\
&\uf{F}^fR\grupo{B_2}_{\beta_2}\ar[uu]_{{\uf{F}}^f\psi_2}&&
{\begin{array}{c}\uf{F}^f(R\grupo{B_1}_{\beta_1}\oplus
R\grupo{B_2}_{\beta_2})\\\parallel\\\uf{F}^fR\grupo{B_1}_{\beta_1}\oplus
\uf{F}^fR\grupo{B_2}_{\beta_2}\end{array}}\ar[uu]^{\uf{F}^fp_1}\ar[ll]_<(.3){\uf{F}^fp_2}
\ar`[ru]`[uuu]^{{\uf{F}}^f(\psi_1,\psi_2)}`[ll][lluu]
\ar`[dl]`[lll]^{{\uf{F}}^f(\psi_1,\psi_2)}`[uu][lluu]&\\&&&&}$$
This diagram represents two morphisms $a\leftarrow c\rightarrow b$
in $\varphi\downarrow(\uf{F}^f\downarrow\uf{F}^f)$, where $c$ is
the object given by $\varphi$,
$\left(\begin{array}{c}\tau_1\\\tau_2\end{array}\right)$, and any
of the outer arrows in the diagram. Hence
$\varphi\downarrow(\uf{F}^f\downarrow\uf{F}^f)$ is connected and
the proof is finished.
\end{proof}

Now we are ready to prove Proposition \ref{extigrande}.

\begin{proof}[Proof of Proposition \ref{extigrande}]
Let $R\grupo{A}_\alpha$ be a possibly big $T$-controlled
$R$-module. There is a chain of equalities and natural
isomorphisms
\begin{eqnarray*}
  Y'\uf{F}^fR\grupo{A}_\alpha &\st{(\mathrm{a})}=& \colim [YP\colon Y\downarrow Y'\uf{F}^f R\grupo{A}_\alpha\rightarrow\mo(\C{M}_R(T'))] \\
   &\st{(\mathrm{b})}=& \colim [YP(1\downarrow Y')=YP\colon J\downarrow \uf{F}^f R\grupo{A}_\alpha\rightarrow\mo(\C{M}_R(T'))] \\
   &\st{(\mathrm{c})}\simeq& \colim [YP(\uf{F}^f\downarrow \uf{F}^f)=Y\uf{F}^fP\colon J\downarrow R\grupo{A}_\alpha\rightarrow\mo(\C{M}_R(T'))] \\
   &\st{(\mathrm{d})}=& \colim [\uf{F}^f_*YP\colon J\downarrow R\grupo{A}_\alpha\rightarrow\mo(\C{M}_R(T'))] \\
   &\st{(\mathrm{e})}=& \uf{F}^f_* \colim [YP=YP(1\downarrow Y')\colon J\downarrow R\grupo{A}_\alpha\rightarrow\mo(\C{M}_R(T))] \\
   &\st{(\mathrm{f})}=& \uf{F}^f_* \colim [YP\colon Y\downarrow Y' R\grupo{A}_\alpha\rightarrow\mo(\C{M}_R(T))] \\
   &\st{(\mathrm{g})}=& \uf{F}^f_*R\grupo{A}_\alpha.
\end{eqnarray*}
Here (a) is a particular case of (\ref{natural-es}); (b) holds since
$(1\downarrow Y')\colon J\downarrow \uf{F}^f
R\grupo{A}_\alpha\rightarrow Y\downarrow Y'\uf{F}^f
R\grupo{A}_\alpha$ is obviously an isomorphism of categories by
definition of $Y'$, here $1$ is the identity functor and we also use
(\ref{masuna}); for (c) we use Lemma \ref{final-es}, \cite{cwm} IX.3
Theorem 1, and (\ref{masuna}); in (d) we use (\ref{masuna}) and the
commutativity of the second diagram in (\ref{2diag}); in (e) we use
that $\uf{F}^f_*$ is a left-adjoint and hence preserves colimits, as
well as (\ref{masuna}); (f) is analogous to (b); and finally (g) is
again a special case of (\ref{natural-es}).
\end{proof}

\chapter{Proof of Proposition
\ref{esexacto}}\label{proofofesexacto}

\begin{proof}[Proof of Proposition \ref{esexacto}]
One can suppose without loss of generality that
$\F(f)\colon\F(T)\subset\F(T')$ is in fact an inclusion and $f\colon
T\subset T'$ is an inclusion of a subtree. Moreover, since
$\uf{F}^f$ is full and faithful in this case, up to equivalence of
categories, we can regard $\C{M}_R(T)$ as the full subcategory of
$\C{M}_R(T')$ formed by those free $T'$-controlled $R$-modules with
support contained in $\F(T)$. Let $R\grupo{{T'}^0}_{\delta'}$ be the
canonical free $T'$-controlled $R$-module with support $\F(T)$,
$R\grupo{T^0}_\delta$ the free $T'$-controlled $R$-module with $T^0$
the vertex set of $T$ and $\delta\colon T^0\subset T'$ the
inclusion, and $\um{M}$ the left-$R(\F(T))$-right-$R(\F(T'))$-module
$\um{M}=\hom_{\C{M}_R(T')}(R\grupo{{T'}^0}_{\delta'},R\grupo{{T}^0}_{\delta})$.
By (\ref{conmotra}) there is a diagram of functors which commutes up
to natural equivalence
\begin{equation*}
\xymatrix@C=70pt{\mo(\C{M}_R(T))\ar[r]^{\uf{F}^f_*}\ar[d]_{ev_{R\grupo{T^0}_\delta}}^\sim&\mo(\C{M}_R(T'))\ar[d]^{ev_{R\grupo{{T'}^0}_{\delta'}}}_\sim\\
\mo(R(\F(T)))\ar[r]_{-\mathop{\otimes}\limits_{R(\F(T))}\um{M}}&
\mo(R(\F(T')))}
\end{equation*}
Also recall that the vertical arrows are equivalences of
categories, see (\ref{expli}), hence we only need to prove that
$\um{M}$ is flat as a left-$R(\F(T))$-module. We shall use the
planarity test (\cite{rmkt} 3.2.9) which ensures that it is enough
to check that for any finitely generated right-ideal
$\mathfrak{a}\subset R(\F(T))$ the epimorphism induced by
multiplication
$$m\colon\mathfrak{a}\mathop{\otimes}\limits_{R(\F(T))}\um{M}\twoheadrightarrow\mathfrak{a}\um{M},$$
is in fact an isomorphism.

By Proposition \ref{solo1} he ideal $\mathfrak{a}$ is necessarily
principal. Let $\psi\colon R\grupo{T^0}_\delta\rightarrow
R\grupo{T^0}_\delta$ be a generator. An arbitrary element in
$\mathfrak{a}\mathop{\otimes}\limits_{R(\F(T))}\um{M}$ has the
form $\psi\otimes \varphi$ for some controlled homomorphism
$\varphi\colon R\grupo{{T'}^0}_{\delta'}\rightarrow
R\grupo{T^0}_{\delta}$. Suppose that $m(\psi\otimes
\varphi)=\psi\varphi=0$. Consider the set $$A=\set{t\in
{T'}^0\,;\,\varphi(t)\neq 0}.$$ Since $\varphi$ is controlled is
easy to see that $\delta'(A)'\subset\F(T)$, in particular if
$\alpha$ is the restriction of $\delta'$ to $A$ then
$R\grupo{A}_\alpha$ belongs to $\C{M}_R(T)$ so it is a retract of
$R\grupo{T^0}_\delta$,
$$\xymatrix{R\grupo{T^0}_\delta\ar@<.5ex>[r]^r&R\grupo{A}_\alpha\ar@<.5ex>[l]^i},\;\; ri=1.$$
There is also a retraction diagram
$$\xymatrix{R\grupo{{T'}^0}_{\delta'}\ar@<.5ex>[r]^p&R\grupo{A}_\alpha\ar@<.5ex>[l]^j},\;\; pj=1,$$
where $p$ vanishes over the elements in ${T'}^0-A$ and $j$ is
induced by the inclusion $A\subset {T'}^0$.

Clearly $\varphi$ coincides with the following composite
$$R\grupo{{T'}^0}_{\delta'}\st{p}\twoheadrightarrow R\grupo{A}_\alpha\st{i}\hookrightarrow
R\grupo{T^0}_\delta\st{r}\twoheadrightarrow
R\grupo{A}_\alpha\st{j}\hookrightarrow
R\grupo{{T'}^0}_{\delta'}\st{\varphi}\rightarrow
R\grupo{T^0}_{\delta},$$ so
\begin{eqnarray*}
  \psi\otimes\fy &=& \psi\otimes\varphi jrip \\
   &=& \psi\varphi jr\otimes ip \\
   &=& 0.
\end{eqnarray*} Therefore $m$ is injective and the proof is
finished.
\end{proof}

\backmatter

\bibliographystyle{amsalpha}
\bibliography{Fernando}
\end{document}